\documentclass[10pt]{amsart}

\usepackage{amsmath}
\usepackage{amssymb}
\usepackage{amscd}
\usepackage{latexsym}
\usepackage[latin1]{inputenc}

\usepackage{graphicx}
\usepackage{changebar,color}

\newtheorem{theorem}{Theorem}[section]
\newtheorem{corollary}[theorem]{Corollary}
\newtheorem{lemma}[theorem]{Lemma}
\newtheorem{conjecture}[theorem]{Conjecture}

\newtheorem{proposition}[theorem]{Proposition}

\newcommand{\Proof}[1]{\noindent {\bf Proof{#1}.}}

\newcommand{\BOX}[1]{\\ \hspace*{\fill} ${\square}_{\tiny{\mbox{#1}}}$}

\numberwithin{equation}{section}

\newcommand{\Hom}{\operatorname{Hom}}

\def\A{\mathbb{A}}

\def\N{\mathbb{N}}
\def\Z{\mathbb{Z}}
\def\PP{\mathbb{P}}
\def\Q{\mathbb{Q}}
\def\R{\mathbb{R}}
\def\C{\mathbb{C}}
\def\E{\mathbb{E}}

\def\G{\mathbb{G}}
\def\H{\mathbb{H}}

\newcommand{\Qb}{{\overline {\mathbb Q}}}
\newcommand{\Gm}{{{\mathbb G}_m}}

\def\bG{\mathbf{G}}

\newcommand{\cf}{{\em cf. }}

\newcommand{\cA}{{\mathcal  A}}
\newcommand{\cB}{{\mathcal  B}}

\newcommand{\cC}{{\mathcal  C}}

\newcommand{\cE}{{\mathcal  E}}
\newcommand{\cF}{{\mathcal  F}}

\newcommand{\cH}{{\mathcal  H}}
\newcommand{\cI}{{\mathcal  I}}
\newcommand{\cJ}{{\mathcal  J}}

\newcommand{\cL}{{\mathcal L}}
\newcommand{\cM}{{\mathcal M}}
\newcommand{\cN}{{\mathcal N}}

\newcommand{\cO}{{\mathcal O}}
\newcommand{\cOan}{{{\mathcal O}^{\rm{an}}}}
\newcommand{\cP}{{\mathcal P}}
\newcommand{\cU}{{\mathcal  U}}
\newcommand{\cV}{{\mathcal  V}}
\newcommand{\cX}{{\mathcal X}}
\newcommand{\cY}{{\mathcal Y}}

\newcommand{\OK}{{\mathcal O}_K}

\newcommand{\CdRB}{{\mathcal C}_{\rm dRB}}

\newcommand{\Lie}{{\rm Lie\,}}
\newcommand{\Lieb}{{\rm Lie}}

\newcommand{\Spec}{{\rm Spec\, }}

\newcommand{\pr}{{\rm pr}}

\newcommand{\Jac}{{\rm Jac}}
\newcommand{\Pic}{{\rm Pic}}
\newcommand{\mon}{{\rm mon}}

\newcommand{\rB}{{\rm B}}

\newcommand{\oli}{\overline}

\newcommand{\Ext}{{\rm Ext}}

\newcommand{\an}{{\rm an}}
\newcommand{\dR}{{\rm dR}}
\newcommand{\dRB}{{\rm dRB}}
\newcommand{\topo}{{\rm top}}

\newcommand{\LiePer}{{\rm LiePer}}
\newcommand{\Per}{{\rm Per}\,}
\newcommand{\Perb}{{\rm Per}}

\newcommand{\bfe}{{\mathbf e}}

\newcommand{\bfE}{{\mathbf E}}

\newcommand{\HdR}{H_{\rm dR}}

\newcommand{\gp}{{\rm gp}}

\newcommand{\lrasim}{\stackrel{\sim}{\longrightarrow}}
\newcommand{\lra}{\longrightarrow}
\newcommand{\lmt}{\longmapsto}

\newcommand{\MIC}{\mathbf{MIC}}
\newcommand{\Rep}{\mathbf{Rep}}

\newcommand{\hA}{\hat{A}}
\newcommand{\hB}{\hat{B}}
\newcommand{\ha}{\hat{a}}

\newcommand{\hlra}{{\lhook\joinrel\longrightarrow}}

\newcommand{\cDgp}{\text{c$D$-gp}}
\newcommand{\cgp}{\rm{c-gp}}
\newcommand{\ggp}{\rm{gp}}

\newcommand{\mic}{\text{mic}}

\newcommand{\bfGm}{{\mathbf{G}_{m}}}            

\begin{document}

\date{\today}

\title[Algebraization, Transcendence, and $D$-groups]{Algebraization, Transcendence, and $D$-group schemes} 
\author{Jean-Beno\^{\i}t Bost}
\address{J.-B. Bost, D{é}partement de Math{é}matiques, Universit{é}
Paris-Sud,
B{â}timent 425, 91405 Orsay cedex, France}
\email{jean-benoit.bost@math.u-psud.fr}

\begin{abstract} We present a conjecture in Diophantine geometry concerning the construction of line bundles over smooth projective varieties over $\Qb.$ This conjecture, closely related to the Grothendieck Period Conjecture for cycles of codimension 1, is also motivated by classical algebraization results in analytic and formal geometry and in transcendence theory. Its formulation involves  the consideration of $D$-group schemes attached to abelian schemes over algebraic curves over $\Qb.$
We also derive the Grothendieck Period Conjecture for cycles of codimension 1 in abelian varieties over $\Qb$ from a classical transcendence theorem \emph{\`a la} Schneider-Lang.

\end{abstract}

%
%
\maketitle

\setcounter{tocdepth}{2}
{
\tableofcontents
}

\setcounter{section}{-1}
\section{Foreword}

My aim, in this largely expository article, is to present a conjecture in Diophantine geometry, concerning the construction of line bundles over smooth projective varieties over $\Qb.$ This conjecture is motivated by the classical Grothendieck Period Conjecture (\cf Section \ref{GPC}) and by the philosophy, already advocated in diverse  places (see for instance \cite{Bost01}, \cite{Chambert01}, \cite{BostChambert-Loir07}, \cite{Gasbarri10}), that various results in Diophantine approximation and transcendence theory are arithmetic counterparts, valid in varieties over number fields, or rather in their model of finite type over $\Z$, of geometric algebraicity criteria, concerning  formal objects inside algebraic varieties over some (algebraically closed) field $k$.

Most of the presently known results in transcendence appear actually to be analogues of geometric algebraicity criteria concerning germs $\widehat{V}$  of formal subvarieties along a projective subvariety $Y$ of some ambient variety $X$ over $k$ --- 
by such a $\widehat{V}$, we mean a smooth formal subscheme  $\widehat{V}$ of the completion $\widehat{X}_{Y}$ admitting $Y$ as scheme of definition. (Any such $\widehat{V}$ may be written as the limit 
$$\widehat{V} = \lim_{\stackrel{\lra}{i}} V_{i}$$
of the successive infinitesimal neighbourhoods $V_{i},$ $i \in \N$, of $Y$ in $\widehat{V}$, which are closed subschemes of $X$, of support $\vert V_{i}\vert = Y.$ ) 
These criteria assert that, if $Y$ is smooth, \emph{of dimension at least one}, and if the normal bundle $N_{Y}\widehat{V}$ of $Y$ in $\widehat{V}$ satisfies some suitable positivity condition, then $\widehat{V}$ is algebraic --- roughly speaking, this means that $\widehat{V}$ is a ``branch'' along $Y$ of some  subvariety $W$ of $X$ containing $Y$.

When the base field $k$ is the field $\C$ of complex numbers, that kind of result may be stated in the following terms, which avoid an explicit appeal to formal geometry and so may  look more familiar. In the situation when $k =\C,$ any germ of $\C$-analytic submanifold $\cV$ of $X$ along $Y$  defines a smooth formal germ $\widehat{V}:= \widehat{\cV}_{Y}$ along $Y$ (namely, the limit $\lim_{i} \cV_{i}$ of the successive infinitesimal neighbourhood of $Y$ inside $\cV$; these are projective analytic subspaces in $X,$ which may be identified to projective subschemes over $\C$). Then the above-mentioned algebraicity criteria  assert  that, \emph{when the normal bundle of $Y$ in $\cV$} satisfies a suitable positivity condition, for instance when it \emph{is ample, then $\cV$ is contained in some algebraic subvariety $W$ of $X$ of the same} (complex) \emph{dimension as} $\cV$. That type of geometric result goes back to Andreotti \cite{Andreotti63}.

In transcendence theory, one deals with algebraicity criteria concerning smooth formal germs of subvarieties $\widehat{V}$ through some $K$-rational point $P$  in a variety $X$ over a number field $K$. According to a viewpoint that goes back to Kronecker, it is appropriate to consider a model $\cX$ of $X$ of finite type over the ring of integers $\OK$ of $K$ (hence over $\Z$), in which $P$ extends to a point $\cP$ in $\cX(\OK)$. The algebraicity criteria established in transcendence turn out to deal with a formal germ in the completion $\widehat{\cX}_{\cP}$ along the ``arithmetic curve'' $\cP \simeq \Spec \OK.$ In this Kroneckerian perspective, transcendence results are indeed algebraicity criteria concerning formal germs along \emph{curves}, analogue to the geometric algebraicity criteria \emph{\`a la} Andreotti.

It turns out that, in the context of analytic and formal geometry, algebraicity criteria have been established that concern, not only subvarieties, but also coherent sheaves (for examples, line bundles or vector bundles), notably  by  Grothendieck  (\cite{GrothendieckFGA}, \cite{GrothendieckSGA2}) in the context of formal geometry. In their most basic geometric version, for instance, the algebraization results in    \cite{GrothendieckSGA2}  (also presented in \cite{Hartshorne70})   deal with germs of formal (or analytic) vector bundles along suitable ample projective subvarieties $Y$ of some algebraic variety $X$ over some base field $k$. Their validity requires $Y$ to be \emph{of dimension at least two}. The Kroneckerian viewpoint mentioned above --- in which the arithmetic counterpart of a surface over some base field is an ``arithmetic surface'', that is an integral model of a curve over a number field  --- leads one to expect  that one could formulate, and possibly  establish, some significant arithmetic algebraization criterion, concerning \emph{formal line or vector bundles over the completion $\widehat{X}_{Y}$  of some algebraic variety $X$ over a number field along some projective curve $Y$.}

In this article, I present a conjectural transcendence statement of this kind (Conjecture \ref{Main} \emph{infra}), the validity of which would actually imply some new cases of  the classical Grothendieck Period Conjecture.

 An interesting feature of this conjectural statement is that it introduces differential algebraic groups in a classical Diophantine context, concerning algebraic varieties over number fields. Recall that the role of differential algebra in Diophantine geometry over function fields is well established since the work of Manin on algebraic curves over function fields, culminating with his proof of the geometric Mordell conjecture (\cite{Manin58}, \cite{Manin61}, \cite{Manin63}), and has more recently considerably expanded, in a series of works initiated by the contributions of Buium (\cite{Buium92Annals}, \cite{Buium93}, \cite{Buium93eff}, \cite{BuiumVoloch93}) and Hrushovski  (\cite{Hrushovski96}), which make conspicuous the role of differential algebraic groups in the Diophantine geometry of abelian varieties over function fields\footnote{We refer the reader to \cite{Buium92}, \cite{Buium94}, and \cite{Pillay97BAMS}, \cite{Bouscaren98}, \cite{Marker2000} for more systematic presentations, surveys, and additional references.}.   
The occurrence of nonlinear differential algebraic groups over curves over number fields in  Conjecture \ref{Main}, which reflects the two-dimensional nature of the problem at hand, has appeared to me worthy of attention, and I took the opportunity of the Ol\'eron conference to present it to experts in model theory and differential algebra gathered at the occasion of Anand Pillay's 60th birthday.

Actually, although the content of this work has presently  no explicit link with model theory, it turns out to involve several of the mathematical themes so successfully explored by Anand Pillay during the recent years, notably the interplay between the analytic geometry of compact complex manifolds and algebraic geometry, and the study of algebraic $D$-groups, especially in relation to abelian varieties and their universal vector extensions. This article is dedicated to him, as a token of appreciation and confidence in his mathematical vision.

This paper, like my oral presentation in Ol\'eron, is to a large extent expository: I seriously attempted to  discuss the classical facts relevant to the formulation of Conjecture \ref{Main} in a form accessible to mathematicians of diverse backgrounds (with possibly a limited success, notably in the last sections of this article). Especially I tried to avoid any real knowledge of formal geometry, by putting forward the analytic variants of diverse  results usually formulated in terms of  formal geometry, or by translating statements in formal geometry into equivalent statements involving systems of successive thickenings, to stay in the realm of algebraic geometry. I also tried to present various themes from some unconventional point of view, for instance in emphasizing the role of moduli spaces of vector bundles with integrable connections.  

However, besides Conjecture \ref{Main} itself, I also included some original content, notably in Part 5  a proof of the Grothendieck Period Conjecture in codimension 1 for abelian varieties. Readers interested in this result may only  read Parts 4 and 5, independently of the rest of the article. 

I heartily thank   Daniel Bertrand for generously sharing his insights of transcendence theory and differential algebraic groups over the years, and for helpful remarks on a preliminary version. I am grateful to the referee for useful comments and to J.P. Serre for his remarks on Section 2.1. I also thank Zo\'e Chatzidakis for her gentle insistence that I transform my oral presentation in Ol\'eron into some written contribution, and the Centro di Ricerca Matematica Ennio di Giorgi (Pisa) for its hospitality during the completion  of this article.

During the preparation of this article, the author has been  partially supported by the ANR project {MODIG}\footnote{ANR-09-BLAN-0047.} and by the Institut Universitaire de France.

%
%
%
\section{Algebraization of analytic objects I}\label{sec:AlgAna}
\subsection{Algebraization of compact Riemann surfaces and of projective analytic sets}\label{PRChow} Algebraization of analytic objects (such as varieties and their morphisms, vector bundles, coherent sheaves, etc.) is a central theme in the development of algebraic and analytic geometry at least since the 1830s. Already recognizable in the pioneering work of Abel and Jacobi on elliptic functions and elliptic curves, it appears in a form familiar to modern mathematicians in the work of Puiseux  and Riemann. 

For instance, in the first part of  his memoir on abelian functions \cite{Riemann57IV} --- devoted to a systematic study of what today would be called ``compact Riemann surfaces realized as a finite covering of the projective complex line $\PP^1(\C)$'' ---  Riemann establishes the \emph{algebraicity} of any pair $(C,\nu)$ where $C$ is \emph{a compact connected Riemann surface} and  $\nu : C \longrightarrow \PP^1(\C)$ \emph{a ramified analytic covering} (or equivalently, a nonconstant $\C$-analytic map). 

Namely, he proves that, for any such pair $(C, \nu)$, there exists an irreducible polynomial $P$ in $\C[X,Y]$ (of positive degree in $Y$), and an isomorphism from $C$ to the compact Riemann surface associated to the plane algebraic curve of equation $P(X,Y) = 0$ such that, through this isomorphism, the map $\nu$ (seen as a meromorphic function on $C$) gets identified with the meromorphic function defined by the first coordinate $X$. To achieve this, Riemann constructs a suitable meromorphic function on $C$ (which ultimately will become the second coordinate $Y$) by appealing to the Dirichlet principle.

An important step in the development of algebraization theorems  has been  the theorem of Chow (\cite{Chow49}), which asserts that  \emph{any closed $\C$-analytic subset $X$ of the projective space $\PP^N(\C)$ is algebraic.} In other words, there exists a finite family $(P_{\alpha})_{1\leq \alpha \leq A}$ of homogeneous polynomials in $\C[X_{0},\cdots,X_{N}]$ such that, for any point $(x_{0}:\cdots:x_{N})$ in $\PP^N(\C),$
$$ (x_{0}:\cdots:x_{N}) \in X  \Longleftrightarrow \mbox{for $\alpha = 1,\ldots, A,$ } P_{\alpha}(x_{0},\ldots,x_{N})= 0.$$

The statement of Chow's theorem  clearly did not come as a surprise at the time of the publication of \cite{Chow49} (see for instance H. Cartan's summary of \cite{Chow49} in \emph{Mathematical Reviews}). A significant point in \cite{Chow49} is the formal rigour of its proofs --- based on some algebraicity criterion formulated in terms of intersections numbers with algebraic subvarieties of $\PP^N(\C)$ ---  which links the theme of algebraization of analytic objects to the development of rigorous  foundations for algebraic topology and geometry, in the line of earlier works by Lefschetz, van der Waerden,  and Chevalley.

\subsection{Algebraization of line bundles over complex projective varieties}\label{subsec:algline} Actually, more than forty years before Chow's work, a remarkable variation on this theme of algebraization was initiated by Poincar\'e and Lefschetz during their investigation of \emph{algebraic cycles on complex surfaces} by means of the so-called normal functions. Motivated by techniques and problems of the Italian school of algebraic geometry and by Picard's contributions to the theory of algebraic surfaces, they basically established the following theorem, when $\dim X = 2$ : 

\emph{Let $X$ be a smooth closed $\C$ analytic subvariety of $\PP^N(\C)$ \emph{(necessarily algebraic, according to Chow's theorem).} Then any analytic line bundle $L$ over $X$ is algebraic.}

This result was extended by Hodge (\cite{Hodge41}, p. 214-216)  to higher-dimensional  smooth projective varieties. Kodaira and Spencer \cite{KodairaSpencer53II} gave a new ``modern'' proof of this theorem in 1953, in what probably constitutes the first application of sheaf theory and cohomological techniques to projective complex varieties.


 Let us formulate a few comments on the content of  the Poincar\'e-Lefschetz-Hodge theorem.  

We shall denote $\cO^\an_{X}$ and $\cC_{X}$ (resp. $\cO_{X}$) the sheaf of analytic and complex-valued continuous functions (resp. of regular functions) on $X$ equipped with the usual ``analytic'' topology (resp. with the Zariski topology). 

Recall that, for any analytic line bundle $L$ over $X$, there exist an open covering $\cU := (U_{\alpha})_{\alpha \in A}$ of $X$ (in the analytic topology) and, for every $\alpha \in A,$ an analytic trivialisation of $L$ over $U_{\alpha}$ :
$$s_{\alpha} : \cO^\an_{U_{\alpha}} \stackrel{\sim}{\longrightarrow} L_{U_{\alpha}}.$$
By comparing the trivialisations --- namely by introducing the functions $\phi_{\alpha \beta}$ in $\cO^\an_X(U_{\alpha}\cap U_{\beta})^\ast$ defined by
$$s_{\alpha} = \phi_{\alpha \beta} s_{\beta} \mbox{ over } U_{\alpha}\cap U_{\beta}$$
--- one defines a 1-cocycle $(\phi_{\alpha \beta})$ in $Z^1(\cU, \cO^{\an\ast}_X)$. The class of this cocycle in $H^1(X,\cO^{\an\ast}_{X})$ determines the isomorphism class of $L$, and any cohomology class in $H^1(X,\cO^{\an\ast}_X)$ arises through this construction from a suitable analytic line bundle $L$.

The line bundle $L$ is \emph{algebraic} precisely when the above covering $\cU=:= (U_{\alpha})_{\alpha \in A}$ and trivialisations $(s_{\alpha})_{\alpha \in A}$, may be chosen in such a way that every $U_{\alpha}$ is \emph{Zariski} open in $X$ and every  function $\phi_{\alpha \beta} s_{\beta}$ is \emph{regular}\footnote{that is, given on the Zariski-open set  $U_{\alpha}\cap U_{\beta}$ by the quotient of two (nonvanishing over $U_{\alpha}\cap U_{\beta}$) homogeneous polynomials of the same degree on $\C^{N+1}.$} over $U_{\alpha}\cap U_{\beta}$; then $(\phi_{\alpha \beta})$ defines a 1-cocycle in $Z^1(\cU, \cO^\ast).$

The above formulation  of the theorem of Poincar\'e-Lefschetz-Hodge, in terms of algebraicity of analytic line bundles, is basically its ``modern'' formulation by Kodaira and Spencer. Let us recall how it translates into its ``classical'' formulation \`a la Lefschetz-Hodge, involving (co)homology classes of divisors. The following arguments, now classical,  appear in \cite{KodairaSpencer53I}.

Consider the short exact sequences of sheaves of abelian groups over $X$ defined by the ``exponential'' map $\bfe : = \exp (2\pi i .)$ :
$$0 \longrightarrow \Z_{X}   \longrightarrow \cC_{X} \stackrel{\bfe}{\longrightarrow} \cC_{X}^\ast \longrightarrow 0 $$
and
$$0 \longrightarrow \Z_{X}   \longrightarrow \cO_{X}^\an \stackrel{\bfe}{\longrightarrow} \cO_{X}^{\an\ast} \longrightarrow 0.$$
The abelian group of isomorphism classes of topological (resp. analytic line) bundles over $X$ is naturally identified with $H^1(X, \cC_{X}^\ast)$ (resp. $H^1(X,\cO_{X}^{\an\ast})$). The long exact sequences of cohomology groups associated to the above short exact sequences of sheaves fit into a commutative diagram :
\begin{equation} \label{diagChern}
\begin{CD}
H^1(X, \cC_{X})@>{\bfe}>> H^1(X,\cC_{X}^\ast) @>{\delta}>>  H^2(X,  \Z) @>>>   H^2(X,\cC_{X}) \\
@.                                                   @AAA                               @AAA                          @AAA \\
  @.                                                 H^1(X,\cO^{\an\ast}_X)  @>{\delta^\an}>> H^2(X, \Z)     @>>>   H^2(X, \cO_{X}).
\end{CD}
\end{equation}

The exactness of the first line and the vanishing of $H^1(X,\cC_{X})$ and $H^2(X,\cC_{X})$ define an isomorphism 
\begin{equation}\label{deftopChern}
c_{1, \topo}:= \delta : H^1(X,\cC_{X}^\ast)        \stackrel{\sim}{\longrightarrow}    H^2(X, \Z),
\end{equation}
which maps the isomorphism class of some topological line bundle $L$ to its so-called \emph{first Chern class}. The exactness of the second line in (\ref{diagChern}) precisely asserts that a class $\alpha$ in $H^2(X, \Z)$ belongs to the image of $\delta^\an$ --- or equivalently, is the first Chern class $c_{1}(L)$ of some \emph{analytic} line bundle --- if and only if $\alpha$ belongs to the kernel
$$\ker (H^2(X, \Z) \longrightarrow H^2(X,\cO_{X}^\an))$$
of the map induced by the inclusion of sheaves  
$\Z_{X} \hlra \cO_{X}^\an,$
or equivalently, if the real cohomology class $\alpha_{\R}$ in $H^2(X, \R)$ belongs to
$$\ker (H^2(X, \R) \longrightarrow H^2(X,\cO^\an_X)).$$
In the classical notation of Hodge theory, this is precisely the space $H^2(X, \R) \cap H^{1,1}(X)$ of real 2-cohomology classes on $X$ of type $(1,1).$  In the case of surfaces, considered by Lefschetz, this space may be defined by the classical  vanishing condition 
$$\int_X \alpha \wedge \omega = 0$$
of the integrals  along $\alpha$ of the global regular algebraic 2-forms $\omega$ on $X$.

Besides, an \emph{algebraic} line bundle $L$ may be described in terms of the divisor $D$ of some nonzero rational section $s$ : the section $s$ establishes an isomorphism from $L$ to the line bundle $\cO(D)$, and the class $c_{1}(L)=c_{1}(\cO(D))$ coincides with the class $[D]$ in $H^2(X, \Z)$ Poincar\'e dual to the divisor $D$, seen as a codimension 1 algebraic cycle on $X$.

Taking the above facts into account, Kodaira-Spencer's version of the theorem of Poincar\'e-Lefschetz-Hodge admits the following consequence, which is actually its original version due to Lefschetz and Hodge\footnote{Conversely, to recover Kodaira-Spencer's version from the Lefschetz-Hodge's, one needs to know that any topologically trivial analytic line bundle over $X$ is algebraic :  this follows from the algebraicity of the Albanese variety and of the Albanese morphism of $X$, and from the algebraicity of analytic line bundles over complex abelian varieties. But for the algebraicity of the Albanese morphism, itself a consequence of Chow's theorem (\cf 2.3.1 \emph{infra}), these  results are actually consequences of Hodge theory and of Lefschetz's work on complex abelian varieties.}   : \emph{a class $\alpha$ in $H^2(X,\Z)$ is algebraic --- namely, the class $[D]$ of some algebraic cycle $D$ of codimension 1 on $X$ --- if and only if $\alpha_{\R}$ is of type  $(1,1)$.}

\subsection{GAGA} The diverse algebraicity statements in the previous sections appear today as special instances of Serre's GAGA Theorem (1956, \cite{Serre56}).

To formulate Serre's results, consider a complex algebraic variety $X$. From any  algebraic coherent sheaf $F$ over $X$ equipped with the Zariski topology --- for example, an algebraic vector bundle $E$ over $X$, defined by some 1-cocycle $(\phi_{\alpha \beta}) \in Z^1((U_{\alpha}), GL_{N}(\cO_{X}))$, attached to some Zariski-open covering $(U_{\alpha})$ of $X$, with values in invertible matrices of regular functions --- we deduce an analytic coherent sheaf $F^{\rm an}$ on $X$ equipped with the analytic topology --- for instance, $E^{\rm an}$ is the analytic vector bundle defined by the cocycle $(\phi_{\alpha \beta})$ seen as as an analytic cocycle (that is, as an element of $Z^1((U_{\alpha}), GL_{N}(\cO_{X}^\an))$. This is a straightforward consequence of the facts that the analytic topology of $X$ is finer than its Zariski topology,  and that, for every Zariski open subset $U$ of $X$, $\cO_{X}(U)$ is a subring of $\cO_{X}^\an(U).$  

These facts also imply the existence of canonical ``analytification maps'' between cohomology groups :
\begin{equation}
\label{anmaps}
H^i(X,F) \longrightarrow H^i(X^{\rm an}, F^{\rm an}).
\end{equation}
Here $X$ (resp. $X^\an$) denotes the variety $X$ equipped with the Zariski topology (resp. the underlying analytic space, which topologically is the set of complex points of $X$ equipped with the usual ``analytic'' topology).

Serre's GAGA Theorem is the conjunction of the following two statements :

\noindent \textbf{GAGA Comparison Theorem.} \emph{For any projective complex variety $X$ and any coherent algebraic sheaf $F$ on $X$, the ``analytification maps'' (\ref{anmaps}) are isomorphisms:}
\begin{equation}
\label{anmapsiso}
H^i(X,F) \stackrel{\sim}{\longrightarrow} H^i(X^{\rm an}, F^{\rm an}).
\end{equation}

\noindent \textbf{GAGA Existence Theorem.} \emph{For any projective complex variety $X$ and for any analytic coherent sheaf $\cF$ on $X^{\rm an}$, there exists some algebraic coherent sheaf $F$ over $X$ \emph{(unique up to unique isomorphism)} such that $\cF$
is isomorphic to $F^{\rm an}$} (as analytic coherent sheaf over $X^\an$).

Let us stress that the projectivity assumption in the GAGA Theorem is essential (see Section \ref{subscec:alganstructures} for a discussion of counterexamples in the quasi-projective situation).

The Poincar\'e-Lefschetz-Hodge Theorem is nothing but the special case of the GAGA Existence Theorem concerning line bundles over smooth varieties. 

Chow's Theorem also follows from the GAGA Existence Theorem --- with the notation of paragraph (\ref{PRChow}), it follows from this theorem applied to $\cO^{\rm an}_{X}$, seen as a coherent analytic sheaf over $\PP^N(\C)^{\rm an}.$ Observe also that conversely, by considering graphs, Chow's theorem implies the comparison isomorphism (\ref{anmapsiso}) when $i = 0$ and $F$ is a vector bundle.

Serre's proof of GAGA Theorems is the archetype of ``modern cohomological proofs'' and, beside its considerable importance in itself, has also played an important role as a model for the development of cohomological techniques in algebraic and formal geometry.

To establish the GAGA Comparison Theorem, using that $X$ may be embedded into some  projective space $\PP^N_{\C}$, one reduces to the special case $X = \PP^N_{\C}$. In that case, Serre's proof relies on some ``algebraic d\'evissage of $F$'' by means of a left resolution by algebraic coherent sheaves that are direct sums of line bundles of the form $\cO_{\PP^N}(k)$, $k \in \Z,$ combined with a direct computation of the algebraic and analytic cohomology groups in (\ref{anmapsiso}) when $F = \cO_{\PP^N}(k).$

The proof of the GAGA Existence Theorem may be seen as a deep amplification and simplification of Kodaira-Spencer's proof in \cite{KodairaSpencer53II}. Besides the Comparison Theorem previously established, it relies on the finite dimensionality of the analytic cohomology groups $H^i(X^{\rm an}, \cF)$ attached to an arbitrary analytic coherent sheaf $\cF$ on $X$. This result, of analytic nature,  was established by Cartan and Serre (\cite{CartanSerre53}) with $X^{\rm an}$ an arbitrary compact complex analytic space. Actually only the degree $i=1$ case of the finiteness theorem of Cartan-Serre is used in  the proof of the Existence Theorem. When $X$ is smooth and $\cF$ is a line bundle, it was established by Kodaira and Spencer as a consequence of the description of $H^i(X^{\rm an}, \cF)$ by means of harmonic forms and of the fact that elliptic differential operators on compact manifolds are Fredholm.

\section{Algebraization of analytic objects II : comments and applications}\label{sec:AlgAnaII}

\subsection{Un peu d'histoire}\label{subsec:histoire}

I would like to stress that the content of the previous sections provides a very fragmentary image of the history of algebraization theorems, a topic especially rich in results and techniques, where the evolution of ideas over the long term seems rather difficult to untangle.

To illustrate this last point, let me indicate that  algebraicity theorems \emph{\`a la} Chow may be derived from B\'ezout-type bounds on intersection multiplicities. That line of argument  appears for instance in Poincar\'e's survey article on abelian functions \cite{Poincare02}, when he proves that a compact complex torus imbedded in a complex projective space is actually algebraic (see \emph{loc. cit.}, Section 2, 53--56).  It constitutes the central point in Chow's proof in \cite{Chow49}, and more recently, plays a key role in the work of Hrushovski and Zilber on Zariski geometries (see \cite{HrushovskiZilber96}, section 7). The influence of Poincar\'e's work on \cite{Chow49} and \cite{HrushovskiZilber96} seems unclear, and \cite{Poincare02} could be a striking example of double \emph{plagiat par anticipation} by Poincar\'e.

Another approach due to Serre to Chow's Theorem --- which appears as an anonymous contribution in \cite{Anonymous56} --- consists in deriving it from the fact  that the transcendence degree over $\C$ of the field  $\cM(X)$ of meromorphic functions on some compact connected complex manifold $X$ is at most its (complex) dimension :
\begin{equation}\label{degtrleqdim}
{\rm degtr}_{\C} \cM(X) \leq \dim X.
\end{equation}

Indeed, if  $X$ is analytically embedded in $\PP^N(\C)$, its Zariski closure $\overline{X}^{\rm Zar}$ is irreducible, the field 
$\C(\overline{X}^{\rm Zar})$ of rational function on $\overline{X}^{\rm Zar}$ may be identified to a subfield of the field of meromorphic function $\cM(X)$, and the upper bound (\ref{degtrleqdim}) implies that the Zariski closure $\overline{X}^{\rm Zar}$ of $X$ in $\PP^N(\C)$ has dimension at most $\dim X$, hence equal to $\dim X$. Besides, the irreducibility of $\overline{X}^{\rm Zar}$ implies its connectedness and the connectedness of its subset $\overline{X}^{\rm Zar}_{\rm reg}$ of smooth points in the analytic topology. This connectedness is a GAGA-type statement which goes back to Puiseux \cite{Puiseux51}, Section I, in the case of plane curves; Puiseux's original proof actually extends to higher-dimensional varieties (see for instance \cite{Shafarevich77}, Section VII.2), and probably constitutes, with other arguments in \cite{Puiseux50} and \cite{Puiseux51}, the first proof of such results satisfactory according to modern standards. The connectedness of $\overline{X}^{\rm Zar}_{\rm reg}$ and its density in $\overline{X}^{\rm Zar}$ for the analytic topology,   together with the inclusion $X \subset \overline{X}^{\rm Zar}$ and the equality of dimension $\dim X = \dim \overline{X}^{\rm Zar}$,  imply the equality $X =\overline{X}^{\rm Zar},$ that is the algebraicity of $X$. 

In turn, proofs of  the upper bound (\ref{degtrleqdim}) appear to have a complicated history --- this bound seems to have been established for the first time in a completely satisfactory way by Serre (\cite{Serre53}, \S 3) and Thimm (\cite{Thimm54}). In \cite{Siegel55}, Siegel discusses the history of the question and gives an ingenious ``elementary" proof, directly influenced by Poincar\'e's article\footnote{Curiously enough, Siegel points out  the relation of Chow's paper with Poincar\'e's article, but does not seem aware that Chow's Theorem may be derived from (\ref{degtrleqdim}).}  \cite{Poincare02} and actually very close to the proof in \cite{Serre53}.  Conversely, as observed in \cite{Remmert56},  (\ref{degtrleqdim}) is an easy consequence of Chow's Theorem and Remmert proper image theorem. In turn, both these theorems may be derived from the fundamental extension theorems concerning complex  analytic sets, due to Thullen, Remmert, and Stein (see for instance \cite{Mumford76}, Section 4A,  or \cite{Gunning90}, Chapters K and M).

Concerning the history of the Poincar\'e-Lefschetz-Hodge theorem, I refer to the classical analysis by Zariski and to the additional comments by Mumford in \cite{Zariski71} Chapter VII\footnote{In a more mundane vein,  I would simply add that   an especially negative assessment by Lefschetz of the approach of Kodaira-Spencer \cite{KodairaSpencer53II} turns out to be well documented (see for instance \cite{SpencerVita04}, p. 21).}.

\subsection{Algebraic de Rham cohomology}\label{subsec:algebdeRham}

In this section, we apply the GAGA Comparison Theorem to the study of the algebraic de Rham cohomology, in the ``easy'' case of projective smooth varieties. The formalism below seems to appear in printed form in the famous letter of Grothendieck to Atiyah \cite{Grothendieck66}, although algebraic de Rham cohomology already occurs implicitly in diverse classical works on algebraic curves, surfaces, and abelian varieties. See \cite{Hartshorne75} for a systematic presentation of  the de Rham cohomology of algebraic varieties and for references.

\subsubsection{}\label{GAGAdR}  Let $X$ be a smooth projective complex algebraic variety. It is equipped with the algebraic de Rham complex
\begin{equation}\label{algdR}
\Omega^\bullet_{X/\C} : 0 \longrightarrow
\Omega^0_{X/\C} =\cO_{X} \stackrel{d}{\longrightarrow}
\Omega^1_{X/\C} \stackrel{d}{\longrightarrow}
\Omega^2_{X/\C} \stackrel{d}{\longrightarrow} \cdots
\end{equation}
and the hypercohomology groups of this complex of sheaves over $X$ equipped with the Zariski topology define the \emph{algebraic de Rham cohomology groups} of $X$ :
$$\HdR^i(X/\C) := \H^i (X, \Omega^\bullet_{X/\C}).$$

By ``analytification'', the algebraic  de Rham complex (\ref{algdR}) becomes the analytic de Rham complex of the $\C$-analytic manifold $X^{\an}$:
\begin{equation}\label{andR}
\Omega^\bullet_{X^{\an}} : 0 \longrightarrow
\Omega^0_{X^{\an}} = \cO^\an_{X^{\an}} \stackrel{d}{\longrightarrow}
\Omega^1_{X^{\an}} \stackrel{d}{\longrightarrow}
\Omega^2_{X^{\an}} \stackrel{d}{\longrightarrow} \cdots
\end{equation}
The hypercohomology groups of $\Omega^\bullet_{X^{\an}}$ define the 
\emph{analytic de Rham cohomology groups} of $X^{\an}$ $\H^i (X^{\an}; \Omega^\bullet_{X^{\an}}),$ and ``analytification'' defines canonical $\C$-linear maps:
\begin{equation}\label{analytifdR}
\H^i (X, \Omega^\bullet_{X/\C}) \longrightarrow \H^i (X^{\an}, \Omega^\bullet_{X^{\an}}).
\end{equation}
The algebraic (resp. analytic) de Rham cohomology groups are related to the algebraic (resp. analytic) ``Hodge cohomology groups''
$H^q(X,\Omega^p_{X/\C})$ (resp. $H^q(X^{\an},\Omega^p_{X^{\an}})$) by the usual spectral sequences
$$E_{1}^{p,q} = H^q(X,\Omega^p_{X/\C}) \Rightarrow \H^{p+q}(X, \Omega^\bullet_{X^{\an}})$$
$$\mbox{(resp. $E_{1}^{p,q} = H^q(X^{\an},\Omega^p_{X^{\an}}) \Rightarrow \H^{p+q}(X^{\an}, \Omega^\bullet_{X^{\an}})$)}.$$
The formation of these spectral sequences is compatible with analytification. Consequently, from the GAGA comparison isomorphisms
$$H^q(X,\Omega^p_{X/\C}) \stackrel{\sim}{\longrightarrow} H^q(X^{\an},\Omega^p_{X^{\an}}),$$
we deduce that the analytification maps (\ref{analytifdR}) from algebraic to analytic de Rham cohomology groups are isomorphisms.

Besides, according to the analytic Poincar\'e Lemma, the inclusion of the locally constant sheaf $\C_{X^{\an}}$ into $\cO^\an_{X^{\an}}$ defines a quasi-isomorphism of complex of sheaves on $X^{\an}$:
$$\C_{X^{\an}} \stackrel{q.i.}{\longrightarrow} \Omega^\bullet_{X^{\an}},$$
and consequently an isomorphism of (hyper)cohomology groups:
\begin{equation}\label{PdR}
H^i (X^{\an}, \C) \stackrel{\sim}{\longrightarrow} \H^i (X^{\an}, \Omega^\bullet_{X^{\an}}).
\end{equation}
The isomorphisms (\ref{analytifdR}) and (\ref{PdR}) define by composition an isomorphism of finite-dimensional $\C$-vector spaces:
\begin{equation}\label{dRB}
\begin{array}{rcl}
\HdR^i(X/\C)  & \longrightarrow   & H^i(X^{\an}, \C)   \\
\beta  & \longmapsto   & \beta^{\an}.     
\end{array}
\end{equation}

\subsubsection{}\label{subsub:algdRk} Observe that the definition of the algebraic de Rham cohomology makes sense for any smooth projective variety $X_{0}$ defined over an arbitrary base field $k$. Indeed we may consider the algebraic de Rham complex
\begin{equation}\label{algdRk}
\Omega^\bullet_{X_{0}/k} : 0 \longrightarrow
\Omega^0_{X_{0}/k} =\cO_{X_{0}} \stackrel{d}{\longrightarrow}
\Omega^1_{X_{0}/k} \stackrel{d}{\longrightarrow}
\Omega^2_{X_{0}/k} \stackrel{d}{\longrightarrow} \cdots
\end{equation}
and define 
$$\HdR^i(X_{0}/k) := \H^i (X_{0},  \Omega^\bullet_{X_{0}/k}).$$

These are finite-dimensional $k$-vector spaces, and when $k$ is a subfield of $\C$, this construction defines a natural ``form over $k$'' of the cohomology with complex coefficients $H^i(X^{\an}; \C)$ of the $\C$-analytic manifold $X^{\an}$ attached to complex algebraic variety $X := X_{0} \otimes_{k} \C$ deduced from $X_{0}$ by extending the base field from $k$ to $\C$. Indeed, by composing a straightforward base change isomorphism and the comparison isomorphism (\ref{dRB}), we obtain a canonical isomorphism
\begin{equation}\label{dRkB}
\HdR^i(X_{0}/k) \otimes_{k}\C \stackrel{\sim}{\longrightarrow} \HdR^i(X/\C)   \lrasim    H^i(X^{\an}, \C).
\end{equation}

\subsubsection{Example I. Smooth projective curves.}

Let $X_{0}$ be a smooth, projective, geometrically connected curve, of genus $g$, over $k$. Then $\HdR^i(X_{0}/k)$ vanishes if $i>2$ and is a canonically isomorphic to $k$ when $i =0$ or $2$. The first de Rham cohomology group $\HdR^1(X_{0}/k)$ is a $2g$-dimensional $k$-vector space. It may be identified with the quotient of the space of  meromorphic 1-forms over $X_{0}/k$ of the second kind (that is, with vanishing residues) by its subspace $dk(X_{0})$ formed by the differentials of rational functions $k(X_{0})$ over $X_{0}.$

For instance, when $k$ is a field of characteristic $\neq 2, 3$, if $X_{0}$ is an elliptic curve $E$ of plane equation 
$$y^2 = 4 x^3 -g_{2}x -g_{3},$$
then $\HdR^1(E/k_{0})$ is a $2$-dimensional $k$-vector space with basis $([\alpha],[\beta])$, where $\alpha : = dx/y$ and  $\beta := x.dx/y.$

\subsubsection{Example II. The first Chern class in algebraic de Rham cohomology.} The morphism of sheaves of abelian groups  over $X_{0}$
$$
\begin{array}{rccc}
  d\log : & \cO^\ast_{X_{0}}   & \lra &\Omega ^1_{X_{0}/k}  \\
& \phi  & \lmt   &    d\phi/\phi
\end{array}
$$
takes its values in the subsheaf $\Omega ^{1{\rm closed}}_{X_{0}/k}$ of closed 1-forms. Therefore it defines a morphism of complex of sheaves 
$$d\log :  \cO^\ast_{X_{0}} \lra \Omega^\bullet_{X_{0}/k}[1],$$
and finally of (hyper)cohomology groups
$$H^1(X_{0}, \cO^\ast_{X_{0}}) \lra \H^1(X_{0}, \Omega^\bullet_{X_{0}/k}[1]) = \H^2(X_{0}, \Omega^\bullet_{X_{0}/k}).$$

 The map so defined will be denoted:
 $$c_{1, \dR} : \Pic (X_{0}) := H^1(X_{0},  \cO^\ast_{X_{0}}) \lra H^2_{\dR}(X_{0}/k).$$
 It sends the class of the line bundle $L$ over $X_{0}$ defined by a cocycle 
 $(\phi_{\alpha \beta})$ in $Z^1(\cU, \cO^\ast_{X_{0}})$ to the class of the (hyper)cocycle
 $(d\phi_{\alpha \beta}/\phi_{\alpha \beta})$ in $Z^1(\cU, \Omega ^{1{\rm closed}}_{X_{0}/k})$, identified to a subspace of $Z^2(\cU, \Omega^\bullet_{X_{0}/k}).$
 
This construction of the first Chern class in algebraic de Rham cohomology is compatible with the topological first Chern class defined in (\ref{deftopChern}):

\begin{lemma}\label{compChern} Assume that $k$ is a subfield of $\C$, and consider a smooth projective variety $X_{0}$ over $k$,   the complex algebraic projective variety $X:= X_{0}\otimes_{k}\C$, and the associated $C$-analytic manifold $X^{\an}$, as in \ref{subsub:algdRk}.
Let $L$ be a line bundle over $X_{0}$, let $L_{\C}$ be the algebraic line bundle over $X$ deduced from $L$ by extension of scalars from $k$ to $\C$, and let $L^{\an}_{\C}$ be the associated analytic line bundle over $X^{\an}.$

The morphism
$$
\begin{array}{rcccl}
 \HdR^i(X_{0}/k) & \lra  & \HdR^i(X/\C)  & \lrasim & H^i(X^{\an}, \C) \\
\alpha  & \lmt  & \alpha_{\C}:= \alpha \otimes_{k}1_{\C} & \lmt & \alpha^{\an}_{\C}   
\end{array}
$$
maps $c_{1, \dR}(L)$ to $2\pi i\, c_{1,\topo}(L_{\C}^{\an}).$
\end{lemma}

To prove this Lemma, it is enough to consider the case $k = \C.$ Then it follows from the fact that the composite morphism of sheaves over $X^{\an}$
$$\cOan \stackrel{\bfe}{\lra} \cOan^\ast \stackrel{d\log}{\lra} \Omega^1_{X^{\an}}$$
is\footnote{The precise definition of the map $\alpha \mapsto \alpha^{\an}_{\C}$ actually involves  the specific sign conventions  used in homological algebra and sheaf cohomology.  The ``standard'' convention used in \cite{Deligne71} indeed introduces a minus sign in the above compatibility relation :  
$c_{1, \dR}(L)_{C}^{\an} = - 2\pi i\, c_{1,\topo}(L_{\C}^{\an}).$
 
In the sequel, we shall generally neglect these delicate problems of signs involved in various ``canonical'' isomorphisms and their compatibility --- although the important sign issue encountered in Section \ref{PrelAb} (see notably (\ref{sign}) and (\ref{altphi})) would plead for a more careful treatment, on the model of \cite{BerthelotBreenMessing82}, Section V.1
.} $2\pi i \, d.$

\subsubsection{Amplification: modules with integrable connections and de Rham cohomology}\label{deRhamcoeff}
In the last sections of this article, we shall use a generalization of the previous results, concerning cohomology with coefficients not only in $\C$, but in local systems of finite-dimensional $\C$-vector spaces. 

Let $(E,\nabla)$ be a ``module with integrable connection'' over $X$, namely a vector bundle $E$ over $X$ equipped with a connection
$$\nabla: E \lra E \otimes_{\cO_{X}}\Omega_{X/\C}^1$$
with vanishing curvature. Then $\nabla$ canonically extends to morphisms of sheaves over $X$
$$\nabla : E \otimes_{\cO_{X}}\Omega_{X/\C}^l \lra E \otimes_{\cO_{X}}\Omega_{X/\C}^{l+1}$$
which satisfy the Leibniz rule --- namely, for any sections $\omega$ of $\Omega^k_{X/\C}$ and $\alpha$ of $E \otimes_{\cO_{X}}\Omega_{X/\C}^\ast$,
$$\nabla(\omega \wedge \alpha) = d\omega \wedge \alpha + (-1)^k \omega \wedge \nabla \alpha$$
--- and the relation
$$\nabla \circ \nabla = 0.$$  
Consequently we may define:
\begin{equation}\label{coeffalg}
H^i_{\dR}(X/\C,(E,\nabla)):=\H^i(X, (\Omega^\bullet_{X/\C}\otimes_{\cO_{X}}E, \nabla)).
\end{equation}

By analytification, we obtain a complex of sheaves $(\Omega^\bullet_{X^\an}\otimes_{\cO^\an_{X}}E^\an, \nabla)$ on $X^\an$ from $(\Omega^\bullet_{X/\C}\otimes_{\cO_{X} }E, \nabla)$, and we may define
\begin{equation}\label{coeffan}
H^i_{\dR}(X^\an,(E^\an,\nabla)):=\H^i(X^\an, (\Omega^\bullet_{X^\an}\otimes_{\cO^\an_{X}}E^\an, \nabla)).
\end{equation}
An application of GAGA similar to the one in paragraph \ref{GAGAdR} shows that (\ref{coeffalg}) and (\ref{coeffan}) are finite-dimensional vector spaces and that the analytification morphisms
\begin{equation}\label{GAGAdRcoeff}
H^i_{\dR}(X/\C,(E,\nabla)) \lra H^i_{\dR}(X^\an,(E^\an,\nabla))
\end{equation}
are isomorphisms.

Besides, the ``analytic de Rham complex with coefficients'' $(\Omega^\bullet_{X^\an}\otimes_{\cO^\an_{X}}E^\an, \nabla)$ is a resolution of the local constant sheaf $E^h$ of finite-dimensional complex vector spaces (of dimension the rank of $E$) defined by the $\C$-analytic sections of $E^\an$ which are ``horizontal'', that is in the kernel of $\nabla.$ In other words, we have an ``analytic Poincar\'e lemma with coefficients'' over $X^\an$, 
$$E^h \stackrel{q.i}{\lra} (\Omega^\bullet_{X^\an}\otimes_{\cO^\an_{X}}E^\an, \nabla),$$
and consequently an isomorphism of (hyper)cohomology groups:
\begin{equation}\label{dRthcoeff}
H^i(X^\an, E^h) \lrasim H^i_{\dR}(X^\an,(E^\an,\nabla)).
\end{equation}

The isomorphisms (\ref{GAGAdRcoeff}) and (\ref{dRthcoeff}) define by composition an isomorphism
$$H^i_{\dR}(X/\C,(E,\nabla)) \lrasim H^i(X^\an, E^h).$$

When $X=X_{0}\times_{k}\C$ and $(E,\nabla)$ are defined over some subfield $k$ of $\C$, we may define
$$H^i_{\dR}(X_{0}/k,(E,\nabla)):=\H^i(X_{0}, (\Omega^\bullet_{X_{0}/k}\otimes_{\cO_{X_{0}}}E, \nabla)).$$
It is a finite-dimensional $k$-vector space, which defines a natural ``form over $k$'' of the cohomology $H^i(X^\an, E^h)$ with coefficients in the local system $E^h.$ 

\subsection{Algebraic and analytic structures, and moduli spaces of vector bundles with integrable connections}\label{subscec:alganstructures}

\subsubsection{} Applied to graphs of morphisms, Chow's Theorem shows that, for any two \emph{projective} complex varieties $X_{1}$ and $X_{2}$ (say smooth for simplicity), the analytification map defines a bijection :
$$
\begin{array}{rcl}
\left\{\substack{\mbox{morphisms $\phi : X_{1} \rightarrow X_{2}$}\\ \mbox{of complex algebraic varieties}}\right\} &
\lrasim &
\left\{\substack{\mbox{morphisms $\psi : X^{\an}_{1} \rightarrow X^{\an}_{2}$}\\ \mbox{of complex analytic manifolds}}\right\} \\
\phi & \lmt & \phi^{\an}.
\end{array}
$$

(See for instance \cite{Mumford76}, Section 4B,  for details.)

In particular, $X_{1}$ and $X_{2}$ are isomorphic as complex algebraic varieties if and only if $X^\an_{1}$ and $X^\an_{2}$ are isomorphic as complex analytic manifolds. Moreover, for any smooth  projective complex algebraic variety $X$, the algebraic variety structure  of $X$ is uniquely determined by the structure of $\C$-analytic manifold $X^\an$ it induces.

This does not hold anymore for general quasi-projective varieties. In this section, we want to discuss a remarkable families of counterexamples, namely of pairs $(X_{1}, X_{2})$ of smooth quasi-projective complex algebraic varieties such that $X_{1}^\an$ and $X_{2}^\an$ are ``naturally'' isomorphic complex manifolds, although $X_{1}$ and $X_{2}$ are not algebraically isomorphic. 

The GAGA Existence Theorem will actually play a crucial role in the construction of these counterexamples, which are built from  moduli spaces of vector bundles with integrable connections of a given rank $N$ on a smooth projective variety $M$, and from spaces of representations of degree $N$ of the fundamental group of $M^\an.$ When $N =1$, these spaces have been classically considered by Severi and Conforto, and then by  Rosenlicht and Serre, during the decades around 1950. For arbitrary $N \geq 1$, they have been investigated thoroughly by Simpson (\cite{Simpson94I}, \cite{Simpson94II}; see also \cite{LePotier91} for a survey).

\subsubsection{}\label{cor(i)(ii)} Let $M$ be a smooth connected projective complex algebraic variety, and let $o$ be a (complex) point of $X$. Choose a positive integer $N$, and consider the following kinds of data :

(i) 3-uples $(E, \nabla, \psi)$ consisting in \emph{a vector bundle $E$ of rank $N$ over $M,$ an integrable connection $\nabla$ on $E$, and a ``rigidification'' }$\psi$ of $E$ at $o$, namely an isomorphism of $\C$-vector spaces 
$$\psi : E_{o} \lrasim \C^N.$$

(ii) \emph{Representations of degree $N$
$$\rho : \Gamma \lra GL_{N}(\C)$$
of the fundamental group} $\Gamma := \pi_{1}(M^\an, o)$ of the complex analytic manifold $M^\an$ with base point $o$.

Observe that we may consider $\C$-analytic versions of data of type (i), namely:

$\mbox{(i)}^\an$ 3-uples $(E^\an, \nabla^\an, \psi)$ consisting in \emph{an analytic vector bundle $E$ of rank $N$ over $M^\an,$ an integrable analytic connection $\nabla^\an$ on $E^\an$, and a rigidification} $\psi$ of $E^\an$ at $o$.

The notion of isomorphisms between two data of type (i), or between two data  of type $\mbox{(i)}^\an$, is defined in the obvious manner as an isomorphism of (algebraic or analytic) vector bundles, compatible with the connections and rigidifications. Observe that, when such an isomorphism exists, it is actually unique.

Through analytification, any data $(E, \nabla, \psi)$ of type (i) determines a data $(E^\an, \nabla^\an, \psi)$ of type $\mbox{(i)}^\an$. Conversely, GAGA Theorems show that any data of type $\mbox{(i)}^\an$ may be obtained by analytification from some data of type (i), that is uniquely determined (up to unique algebraic isomorphism)\footnote{To ``algebraize'' an analytic connection $\nabla^\an$ over $E^\an$ by means of GAGA Comparison Theorem, identify (algebraic or analytic) connections with (algebraic or analytic) splittings of the Atiyah extension of $E$, $0\rightarrow \Omega^1_{M}\otimes E \rightarrow J_{M}^1E \rightarrow E \rightarrow 0$, defined by the vector bundle  $J_{M}^1E$ of 1-jets of $E$ over $M$.}.

In turn, to any data of type  $\mbox{(i)}^\an$ is associated its monodromy representation in the fiber $E_{0}$ of the flat vector bundle $(E^\an, \nabla^\an)$, which may be identified to a $GL_{N}(\C)$-representation by means of the rigidification $\psi$: 
$$\rho : \Gamma \lra GL(E_{o}) \xrightarrow{\psi . \psi^{-1}} 
GL_{N}(\C).$$

Conversely, we may introduce the universal covering $(\tilde{M},\tilde{o})$ of the pointed connected complex manifold $(M^\an, o)$ --- it is a $\Gamma$-covering of $M^\an$ --- and the trivial vector bundle $\tilde{E} := \tilde{M}\times \C^N$ of rank $N$ over $\tilde{M}$, equipped with the ``trivial'' integrable analytic connection $\tilde{\nabla} := d \otimes Id_{\C^N}.$ If $\rho : \Gamma \lra GL_{N}(\C)$ denotes an arbitrary representation, the action of $\Gamma$ on $\C^N$ defined by $\rho$ makes $(\tilde{E}, \tilde{\nabla})$ a $\Gamma$-equivariant analytic vector bundle with integrable connection, which moreover is naturally rigidified at $\tilde{o}$. 
This equivariant rigidified vector bundle with integrable connection over $(\tilde{M}, \tilde{o})$ descends to some rigidified vector bundle of rank $N$ with integrable connection  $(E^\an, \nabla^\an, \psi)$ on the pointed complex manifold $(M^\an, o)$.

These last two constructions are clearly inverse of each other and establish a natural bijection between (isomorphism classes) of data of type   $\mbox{(i)}^\an$ and representations of type (ii). Combined with the above GAGA correspondence between data of type (i) and  $\mbox{(i)}^\an$, this becomes a natural bijection between (isomorphism classes) of data of type (i) and representations of type (ii).

\subsubsection{}\label{MIC} The set of (isomorphism classes) of data of type (i) coincides with the set of complex points $\MIC_{N}(M,o)(\C)$ of some quasi-projective scheme $\MIC_{N}(M,o)$ over $\C$, which represents the functor which maps a $\C$-scheme (of finite type) $S$ to the isomorphism classes of ``data of type (i) over $S$'', defined as 3-uples $(E,\nabla, \psi)$ where $E$ denotes a locally free coherent sheaf of rank $N$ over $M \times S,$ $\nabla$ an integrable connection on $E$, relative to the projection $M \times S \rightarrow S,$ and $\psi$ a rigidification $E_{\vert o\times S} \lrasim \cO_{S}^{\oplus N}$.

At this level of generality, the existence of the quasi-projective scheme $\MIC_{N}(M,o)$ representing this functor is one of the main results of Simpson in \cite{Simpson94I, Simpson94II}, where it is denoted $\mathbf{R}_{\rm DR}(M,o,N)$. A central point in the construction of  $\MIC_{N}(M,o)$  is the fact that the vector bundles $E$ of rank $N$ over $M$ admitting an integrable connection $\nabla$ constitute a bounded family (see \cite{LePotier91}, Lemme 9, for a concise presentation of Simpson's argument in this specific situation).

The set of representations of type (ii) coincides with the set of complex points $\Rep_{N}(\Gamma)(\C)$ of the quasi-projective (actually affine) scheme $\Rep_{N}(\Gamma)$ over $\C$ which represents the functor which sends a $\C$-scheme of finite type $S$ to the set of representations 
$$\rho : \Gamma \lra GL_{N}(\Gamma(S, \cO_{S})).$$
The existence   of the  scheme    $\Rep_{N}(\Gamma)$ is a straightforward consequence of the existence of a finite presentation for the fundamental group $\Gamma$ (see for instance  \cite{Simpson94II}, Section 5, where this scheme is denoted $\mathbf{R}(\Gamma,N)$ or $\mathbf{R}_{\rm B}(M,o,N)$).

The bijection constructed in \ref{cor(i)(ii)}, by associating  the monodromy representation of its analytification to some data of type (i), defines a bijection:
\begin{equation}\label{anisopoints}
\MIC_{N}(M,o)(\C) \lrasim \Rep_{N}(\Gamma)(\C),
\end{equation}
 which turns out to be defined by a canonical isomorphism of $\C$-analytic spaces
 \begin{equation}\label{aniso}
\mon_{o} : \MIC_{N}(M,o)^\an \lrasim \Rep_{N}(\Gamma)^\an.
\end{equation}
(Compare \cite{Simpson94II}, Section 7. This formally expresses the fact that the construction in \ref{cor(i)(ii)} ``ana\-lytically depends on parameters'' in an arbitrary analytic space.)

\subsubsection{} However, in general, the analytic isomorphism (\ref{aniso}) is \emph{not} induced by an algebraic isomorphism from $\MIC_{N}(M,o)$ to $\Rep_{N}(\Gamma)$. 

This is already the case when $M$ is a smooth connected projective curve $C$ of positive genus $g$ and $N=1$.
Then $$\Pic^\natural (C) := \MIC_{1}({C},o)$$ may be identified with the \emph{universal vector extension}
$E(\Pic_0(C))$ of the connected Picard variety $\Pic_{0}(C)$ of $C$ (see for instance \cite{Messing73}, \cite{MazurMessing74}, \cite{BK09}). Actually, $\Pic^\natural (C)$ classifies pairs $(L,\nabla)$ consisting in a line bundle $L$ of degree 0 over $C$ and a (necessarily integrable) connection $\nabla$ over $L$. The tensor product of line bundles with connections induces a structure of algebraic groups on $\Pic^\natural (C)$. It  fits into the following exact sequence of connected commutative group schemes over $\C$, which displays it as a vector extension of $\Pic_{0}(C)$:
\begin{equation}\label{Picnaturalext}
\begin{array}{crcccclc}
  0 \lra   & \Omega^1(C)   & \lra    & \Pic^\natural(C) & \lra  & \Pic_{0}(C) & \lra 0 \\
  & \alpha &\lmt & [(\cO_{C}, d + \alpha)] & &   &   \\
  & & & [(L,\nabla)] & \lmt & [L] &   &   
\end{array}
\end{equation}

Besides, the representation space $\Rep_{1}(\pi_{1}(C^\an, o))$ may be identified with the torus
$$H^1(C^\an, \Z) \otimes_{\Z}\G_{m} \simeq \G_{m}^{2g},$$
and the monodromy isomorphism (\ref{aniso}) takes the form of an isomorphism of complex Lie groups :
$$\Pic^\natural(C)^\an \lrasim \C^{\ast 2g}.$$

However the description of  $\Pic^\natural(C)$ as a vector extension of an abelian variety easily implies that every morphism of algebraic variety from $\Pic^\natural(C)$ to $\G_{m}$ is constant. \emph{A fortiori}, the algebraic varieties $ \MIC_{1}({C},o) = \Pic^\natural(C)$ and $\Rep_{1}(\pi_{1}(C^\an, o)) \simeq \G_{m}^{2g}$ are not isomorphic\footnote{This occurence of commutative algebraic groups over $\C$ that are analytically,  but not algebraically,  isomorphic has been first pointed out  by Conforto; see \cite{Conforto48}, \cite{Conforto49II}, and  \cite{Severi61}, Appendice.}. 

\subsubsection{}\label{VarMic} For later reference, let us indicate diverse variants of the previous constructions.

First of all, for any base field of characteristic zero and any pointed connected smooth pointed variety $(M,o)$ over $k$, the construction of the quasi-projective scheme $\MIC_{N}(M,o)$ makes sense over $k$ : it classifies data of type (i) over varying $k$-schemes $S$. This follows from a straightforward generalization of the arguments in \cite{Simpson94I}, or (say, when $k$ is a subfield of $\C$) from a descent argument.

When $N=1,$ the tensor product of line bundles with (necessarily integrable) connections makes the quasi-projective scheme $\MIC_{1}(M,o)$ a group scheme, necessarily smooth over $k$. Moreover its connected component $\MIC_{1}(M,o)_{0}$ may be identified with the universal vector extension $E(\Pic_{0}(M))$ of the connected Picard variety $\Pic_{0}(M)$ of $M$. Indeed the obvious analogue of the short exact sequence (\ref{Picnaturalext}) still holds in this setting (see for instance \cite{BK09}, Appendix B).  

When $M$ is the abelian variety $\hat{A}$ dual to some abelian variety $A$ over $k$, this construction identifies the universal vector extension $E(A)$ of $A$ to the $k$-algebraic group
$$\Pic^\natural (\hat{A}) := \MIC_{1}(\hat{A}, 0_{\hat{A}}),$$
which classifies line bundles with (necessarily integrable) connections over $A$, and the short exact sequence  (\ref{Picnaturalext}) becomes the extension defining $E(A)$:
$$0 \lra \E_{\hA} := (\Lie \hA)^\vee \lra E(A) \stackrel{p_{A}}\lra A \lra 0.$$ 

Second, it is convenient to have at one's disposal diverse generalizations of the moduli spaces $\MIC_{N}(M,o)$. For instance, if $(M, o, o')$ denotes a connected smooth projective variety over $k$, endowed with two (possibly equal)  ``base points'' $o$ and $o'$ in $M(k)$, we may construct a quasi-projective scheme $\MIC_{N}(M,o,o')$ that classifies vector bundles $E$ of rank $N$ over $M$, equipped with an integrable connection $\nabla$ and with rigidifications $\psi : E_{o} \lrasim k^N$ and $\psi': E_{o'} \lrasim k^N$ at $o$ and $o'$  (\cf \cite{Simpson94I}, Remark p. 109). Thanks to the morphism 
$$\digamma : \MIC_{N}(M,o,o') \lra \MIC_{N}(M,o)$$
defined by forgetting the rigidifications $\psi'$ at $o'$ and to the action by composition of $GL_{N,k}$ on these rigidifications, $\MIC_{N}(M,o,o')$ becomes a $GL_{N,k}$-torsor over $\MIC_{N}(M,o)$. When $N=1,$ the tensor product again makes  $\MIC_{N}(M,o,o')$ a commutative algebraic group over $k$, and the above structure of $GL_{N,k}$-torsor becomes an extension of commutative algebraic groups:
\begin{equation}
0 \lra \G_{m,k} \lra \MIC_{1}(M,o,o') \lra \MIC_{1}(M,o) \lra 0.
\end{equation}

When $M=\hat{A}$ as above, $o = 0_{\hat{A}},$ and $o'$ is a point $P$ in $\hA(k)$ parameterizing some line bundle $L$ over $A$ (equipped with a rigidification $\epsilon : k\simeq L_{0_{A}}$ and algebraically equivalent to zero), one gets an extension
\begin{equation}\label{extEA}
0 \lra \G_{m,k} \lra \MIC_{1}(\hA,0_{A},P) \lra E(A) \lra 0
\end{equation}
 which may be described as follows. The $\G_{m,k}$-torsor $L^\times$ over $A$, deduced from the total space of $L$ by deleting its zero section may be endowed with a unique structure of $k$-algebraic group which makes the diagram
 \begin{equation}\label{extA}
0 \lra \G_{m,k} \stackrel{\epsilon}{\lra} L^\times \lra A \lra 0
\end{equation}
a short exact sequence of commutative algebraic groups over $k$, and the extension  (\ref{extEA}) coincides with the pullback of the extension (\ref{extA}) by $p_{A}: E(A) \lra A.$ 

\section{Algebraization of formal objects}

\subsection{A Theorem of Grauert-Grothendieck} Since the work of Zariski on ``holomorphic functions'' (\cite{Zariski51}) and its amplification in Grothendieck's new foundations  of algebraic geometry (\cite{GrothendieckFGA}), \emph{formal schemes} and coherent sheaves over them play a central role in modern algebraic geometry. Grothendieck notably established some comparison and existence theorems that relate algebraic and formal geometry over a suitable complete ``adic'' base ring
(\cf \cite{GrothendieckFGA}, \cite{EGAIII1}, \cite{Illusie05}). In SGA2 (\cite{GrothendieckSGA2}), motivated by some earlier work of Grauert, he also used formal geometry to investigate the classical Lefschetz theorems comparing the geometry of projective varieties and of their hyperplane sections.

In the sequel, we shall be concerned by the algebraization theorems of ``Lefschetz type'' established in SGA2 rather than by the earlier ``fundamental'' comparison and existence theorems discussed in \cite{GrothendieckFGA}, \cite{EGAIII1} and \cite{Illusie05}. 

For the sake of simplicity, we first state a (weaker) analytic version of these theorems of Lefschetz type in a special simple case.

\begin{theorem}[Grauert, Grothendieck, \cite{GrothendieckSGA2}]\label{GrGr}
Let $X\hookrightarrow \PP^N_{\C}$ be a smooth projective complex variety of dimension $d$, and let $Y:= X \cap \PP^{N-1}_{\C}$ be a hyperplane section of $X$ of dimension $d-1$.

{\bf Gr1.} If $d \geq 2,$ then for every algebraic vector bundle $E$ over $X$, the restriction map
$$\Gamma(X,E) \longrightarrow \left\{\text{germs of analytic sections of $E$ along $Y$}\right\}$$
is an isomorphism.

{\bf Gr2.} If $d \geq 3,$ any germ of analytic vector bundle $\cE$ on some analytic neighbourhood of $Y$ in $X$ ``extends'' to some coherent sheaf $E$ over $X$. 
\end{theorem}

Observe that, like GAGA, this theorem decomposes into two parts: a ``comparison theorem'' {\bf Gr1}, and an ``existence theorem'' {\bf Gr2}.

Observe also that, according to Serre's GAGA, the vector bundle $E$ in {\bf Gr1} and its space of global sections $\Gamma(X,E)$ may be equivalently taken in the algebraic or in the analytic category. The same remark applies to the coherent sheaf $E$ the existence of which is asserted in {\bf Gr2}. Accordingly, when the conclusion of {\bf Gr2} holds, we shall say that $\cE$ is \emph{algebraizable}.

Let us emphasize that the assumptions on the dimension $d$ are crucial in Theorem \ref{GrGr}.

Indeed {\bf Gr1} trivially fails for $X=\PP^1$, $Y=\{\mbox{point}\}$, and $E=\cO_{X}.$ 

The existence theorem {\bf Gr2} already fails for line bundles when $X$ is the projective plane $\PP_{\C}^2$ and $Y= \PP_{\C}^1$ a projective line in $X$. This follows from Proposition \ref{gerlineb} below, which is a simple consequence of {\bf Gr1}. 

Let $X_{\infty}$ denote a projective line  in $X$ distinct from $Y$, and let us consider the affine plane $\A_{\C}^2 := X \setminus X_{\infty}$ and the affine line $\A_{\C}^1:= \A_{\C}^2 \cap Y$. Choose affine coordinates $(x,y)$ on $\A_{\C}^2$ such that $\A_{\C}^1 = (x=0)$. For any converging power series $f$ in $\C\{T\}$, the equation 
$$y = f(x)$$ defines a germ $T_{f}$ of smooth analytic curve in $X=\PP_{\C}^2$ transverse to $Y= \PP_{\C}^1.$

\begin{proposition}\label{gerlineb} The germ of analytic line bundle $\cOan(T_{f})$ along $\PP_{\C}^1$ in $\PP_{\C}^2$ is algebraizable if and only if the series $f$ belongs to $\C T + \C.$ 
\end{proposition}

Observe also that Theorem \ref{GrGr} admits striking elementary geometric applications. For instance, it implies that \emph{any germ of analytic hypersurface 
 along $\PP^2_{\C}$ in $\PP^3_{\C}$ extends to a global algebraic hypersurface, defined by the vanishing of some homogeneous polynomial in $\C[X_{0},X_{1},X_{2},X_{3}]$.}

\subsection{Formal geometry} In SGA2, Theorem \ref{GrGr} is stated and proved in a more general formulation, in which (i) concerns \emph{formal} sections and vector bundles instead of analytic germs of sections and vector bundles, (ii) makes sense over an arbitrary base field $k$ --- indeed over an arbitrary Noetherian base $S$ --- instead of $\C$, and (iii) holds under  regularity assumptions weaker than the smoothness of $X$, formulated in terms of depth. In this paragraph, we want to give some indication of the generalizations (i) and (ii), while keeping minimal the prerequisites from formal geometry.  

Recall (see for instance \cite{Illusie05}) that, for any Noetherian scheme $X$  and any closed subscheme $Y$ in $X$, a coherent formal sheaf $\cE$ over the formal scheme $\widehat{X}_{Y}$, completion of $X$ along $Y$, ``is'' nothing else than the data of a system $(\cE_{n})_{n \in \N}$ of coherent sheaves on the successive infinitesimal neigbourhoods $Y_{n}$ of $Y$ in $X$ ($Y_{0}:=Y;$ $Y_{n}$ is defined by the $n+1$-th power $\cI_{Y}^{n+1}$ of the ideal sheaf $\cI_{Y}$ of $Y$ in $\cO_{X}$), equipped with isomorphisms
\begin{equation}\label{forsyst}
\cE_{n+1\vert Y_{n}} \lrasim \cE_{n}.
\end{equation}
The coherent formal sheaf   $\cE$ is locally free --- and then called a \emph{vector bundle} --- if and only if, for every $n$, $\cE_{n}$ is a locally free coherent sheaf of $\cO_{Y_{n}}$-modules.

By definition, the space of sections of $\cE$ over $\widehat{X}_{Y}$  ``is'' precisely the projective limit
$$\Gamma(\widehat{X}_{Y}, \cE) := \lim_{\stackrel{\longleftarrow}{n}}\Gamma(Y_{n}, \cE_{n}),$$ defined by means of the isomorphisms (\ref{forsyst}) and of the induced projective system of spaces of sections:
$$ \Gamma(Y_{n+1}, \cE_{n+1}) \stackrel{._{\vert Y_{n}}}{\lra} \Gamma(Y_{n}, \cE_{n+1 \vert Y_{n}}) \lrasim \Gamma(Y_{n}, \cE_{n}).$$

A coherent sheaf $E$ over $X$ defines a formal coherent sheaf  $E_{\vert \widehat{X}_{Y}} := (E_{\vert Y_{n}})$ over $\widehat{X}_{Y}$. A formal coherent sheaf on $\widehat{X}_{Y}$ will be called \emph{algebraizable} if, up to isomorphism, it is of the form $E_{\vert \widehat{X}_{Y}}$ for some coherent sheaf $E$ over $X$.

Using these definitions, we may state a generalized version of Theorem \ref{GrGr} valid for a smooth projective scheme over an arbitrary base field $k$.

\begin{theorem}\label{GrGrfor}
Let $X\hookrightarrow \PP^N_{k}$ be a smooth projective scheme over $k$, of pure dimension $d$, and let $Y := X \cap \PP^{N-1}_{k}$ be some hyperplane section, of dimension $d-1$.

{\bf Gr1.} If $d \geq 2,$ then for any  vector bundle $E$ over $X$, the restriction map
$$\Gamma(X,E) \lra \Gamma(\widehat{X}_{Y}, E_{\vert \widehat{X}_{Y}}) := \lim_{\stackrel{\longleftarrow}{n}}\Gamma(Y_{n}, \cE_{\vert Y_{n}})$$
is an isomorphism.

{\bf Gr2.} If $d \geq 3,$ then any vector bundle $\cE$ over $\widehat{X}_{Y}$ is algebraizable.

\end{theorem} 

Like the proof of Serre's GAGA and of Grothendieck's Comparison and Existence Theorems in \cite{GrothendieckFGA}, \cite{EGAIII1}, \cite{Illusie05}, the proofs in SGA2 are cohomological. For instance, a key point in the proof of {\bf Gr2} is that, since $d\geq 3,$ the Cartier divisor $Y$ has depth $\geq 2$ and the ampleness of $\cO_{X}(Y)_{\vert Y}$ implies that, for every vector bundle $E_{0}$ over $Y$, the cohomology group $H^1(Y,E_{0}\otimes \cO_{X}(-Y)^{\otimes n}_{\vert Y})$ vanishes for $n$ a sufficiently large positive integer  (Lemma of Enriques-Severi-Zariski). This implies that, for any vector bundle $\cE = (E_{n})$ over  $\widehat{X}_{Y}$, the system $(H^1(Y_{n},E_{n}))$ is essentially constant, and consequently
$$H^1(\widehat{X}_{Y}, \cE) = \lim_{\stackrel{\longleftarrow}{n}}H^1(Y_{n}, \cE_{n})$$
is a \emph{finite-dimensional} $k$-vector space. The finite dimensionality of a first cohomology group plays  the same role here
as in the proofs of the Poincar\'e-Lefschetz-Hodge Theorem by  Kodaira-Spencer, and of the GAGA Existence Theorem by Serre. 

Let us also indicate that the results in SGA2 have been extended in diverse directions by Mich\`ele Raynaud (\cite{RaynaudMe75}) and Faltings (\cite{Faltings79}), and that, besides the original cohomological proofs, it is possible to give more ``classical" proofs of Theorems \ref{GrGr} and \ref{GrGrfor},  based on Theorem \ref{thAndreotti} \emph{infra} and its formal variant, which ultimately rely on the use of ``auxiliary polynomials," familiar in Diophantine approximation and transcendence.

\subsection{A Theorem of Andreotti and Hartshorne} Let us mention that diverse algebraization results concerning formal meromorphic functions along subvarieties have also been established, notably by Hironaka-Matsumura (\cite{HironakaMatsumura68}), Faltings (\cite{Faltings80}, \cite{Faltings81}), and Chow (\cite{Chow86}). 

We want to discuss briefly an algebraization result, concerning formal germs along curves, that is related both to the results in \emph{loc. cit.} and to the Grauert-Grothendieck Theorems \ref{GrGr} and \ref{GrGrfor}. For the sake of simplicity, we state it in the analytic framework, in which situation it goes back to Andreotti \cite{Andreotti63} :

\begin{theorem}\label{thAndreotti}
Let $C \hookrightarrow \PP^N_{\C}$ be a smooth connected projective complex algebraic curve, and let $\cV$ be a germ of smooth $\C$-analytic submanifold along $C$ in $\PP^N(\C).$ 

If the normal bundle $N_{C}\cV$ to $C$ in $\cV$ is ample, then $\cV$ is algebraic. 
\end{theorem}

Observe that the normal bundle $N_{C}\cV$ is an analytic vector bundle over $C$, which by GAGA defines an algebraic vector bundle over $C$. When $\dim \cV = 2$, it is a line bundle, and its ampleness is equivalent to the positivity of its degree $\deg_{C}N_{C}\cV.$

In Theorem \ref{thAndreotti}, the algebraicity of $\cV$ precisely means that the dimension $\dim \oli{\cV}^{\rm Zar}$ of its Zariski closure  $\oli{\cV}^{\rm Zar}$ in $\PP^N_{\C}$, which is at least equal to the complex dimension $\dim \cV$ of the complex manifold $\cV$, actually coincides with $\dim \cV$. This is equivalent to the fact that the germ $\cV$ is a ``branch'' along $C$ of some (irreducible) algebraic subset of $\PP^N_{\C}$ containing $C$.

Here again, Theorem \ref{thAndreotti} admits  a formal generalization, valid over any base field, where $\cV$ is a smooth formal subscheme containing $C$ of the formal completion of $\PP^N_{k}$ along a smooth projective $k$-curve. It may also be extended to higher-dimensional situations : the curve $C$ may be replaced by any smooth projective subvariety $Y$, of dimension at least $1$.  This condition is similar to the dimension condition in the assertions {\bf Gr1} in Theorems \ref{GrGr} and \ref{GrGrfor}. Actually  {\bf Gr1} may be derived from Theorem \ref{thAndreotti} and its higher-dimensional and formal generalization  by considering the graphs of analytic or formal sections (see \cite{BostChambert-Loir07}).

In its analytic (resp. formal) form, Theorem \ref{thAndreotti} is a direct consequence --- by the ``anonymous'' argument recalled in Section \ref{subsec:histoire} --- of a result of Andreotti \cite{Andreotti63} (resp. of Hartshorne \cite{Hartshorne68}) which asserts that the field of meromorphic functions (resp. of formal meromorphic functions) on $\cV$ is a field of transcendence degree at most $\dim \cV$ over $\C$ (resp. over $k$). 

Theorem \ref{thAndreotti} may also established by directly estimating the Hilbert function of the Zariski closure of $\cV$, with no recourse to the (formal) meromorphic functions (\cf \cite{Bost01}, Section 3.3, and \cite{Bost06}). This type of argument may be seen as a geometric counterpart of the use of auxiliary polynomials in Diophantine approximation and transcendence proofs.

Algebraization criteria in the style of Theorem \ref{thAndreotti} have been recently reconsidered in \cite{BogomolovMcQuillan01} and \cite{Bost01} in relation to algebraicity properties of leaves of algebraic foliations; see \cite{KebekusSolaToma07} for geometric applications and references, and \cite{Bost04} for similar geometric applications to groups schemes over projective curves.

\subsection{Algebraization over function fields}  The above algebraization theorems, concerning formal ``objects'' over projective varieties on some base field $k$ may be used to derive algebraization theorems over projective varieties on function fields of the form $k(C)$, where $C$ denotes some projective variety over $k$.

We illustrate this general principle by formulating an application of Theorem \ref{thAndreotti} to the algebraicity of formal germs in varieties over the function field $\C(C)$ defined by some smooth projective complex curve $C$. The details of its derivation, which is straightforward, will be left to the reader, as well as the derivation from the formal variant of Theorem \ref{thAndreotti} of a similar algebraicity criterion for formal germs in varieties over a general function field $k(C)$. 

Let $C$ be a smooth projective complex curve and let $\pi : \cX \rightarrow C$ be a projective complex variety fibered over $C$. (In other words, $\pi$ is a flat surjective morphism of complex schemes.)

Let $K:= \C(C)$ be the function field of $C$, and let  $X := \cX_{K}$ be the generic fiber of $\pi$. It is a projective $K$-variety, and conversely, any projective $K$-variety may be realized as the generic fiber of a suitable projective model $\cX$ fibered over $C$ as above.

Let $P$ be a $K$-point of $X$. By the projectivity of $\pi,$ it extends to a section $\cP$ of $\pi$ over $C$.

Consider a smooth formal germ of a subvariety through $P$ in $X$, 
$$\widehat{V} := \lim_{\stackrel{\lra}{i}} V_{i},$$
namely a smooth formal subscheme of the completion $\widehat{X}_{P}$. Here again it is said to be algebraic when its Zariski closure $\oli{\widehat{V}}^{{\rm Zar}_{X}}$ in the $K$-scheme $X$ has the same dimension as $\widehat{V}.$

The $V_{i}$'s are zero-dimensional subschemes of $X =\cX_{K}$ supported by $P$. Their closures in $\cX$ 
$$\cV_{i}:= \oli{V_{i}}^{{\rm Zar}_{\cX}}$$ 
are one-dimensional subschemes of $\cX$ with support $\cP$, and constitute an inductive system
$$\cV_{0} = \cP \hlra \cV_{1} \hlra \cV_{2} \hlra \ldots \hlra \cV_{i} \hlra \cV_{i+1} \hlra \ldots$$
In general this system $(\cV_{i})_{i \in \N}$ does \emph{not} define a formal subscheme of the completion $\hat{\cX}_{\cP}$ smooth over $C$. However it is the case when there exists a germ $\cV$ of analytic submanifold of $\cX^\an$ along $\cP$ that ``extends'' $(\cV_{i})_{i \in \N}$ in the sense that $\cV_{i}$ is the $i$th infinitesimal neighbourhood of $\cP$ in $\cV$.  

\begin{corollary}\label{corthA} With the above notation, if $\widehat{V}$ extends to a germ $\cV$ of a smooth analytic submanifold of $\cX^\an$ along $\cP$ and if the normal bundle $N_{\cP}\cV$ to $\cP$ in $\cV$ is ample, then $\widehat{V}$ is algebraic. 
\end{corollary}

A generalization of this corollary, formulated in terms of formal geometry only, holds when the base field $\C$ is replaced by an arbitrary base field $k$. Namely, \emph{$\widehat{V}$ is algebraic when it extends to a formal subscheme $\hat{\cV}$ of $\hat{\cX}_{\cP}$ smooth over the base curve $C$ and when the normal bundle $N_{\cP}\hat{\cV}$ is ample.}

\section{Algebraization and transcendence}

Various classical results in transcendance theory and Diophantine approximation may be rephrased in geometric terms as algebraization results, asserting the algebraicity of certain formal or analytic subvarieties inside algebraic varieties defined over number fields, provided suitable arithmetic and analytic conditions are satisfied (see for instance \cite{Bost01}, \cite{Chambert01}, \cite{Bost06}, \cite{Gasbarri10}).

In this article, we are concerned with transcendence results of ``Schneider-Lang type'', in the line of the classical theorems of Schneider about the transcendence of values of abelian functions (\cite{Schneider41}, \cite{Schneider57}) and of their modern amplification by Lang (\cite{Lang62, Lang65, Lang66}). We shall content ourselves with   two instances of these transcendence theorems, whose proofs involve only elementary analytic techniques. We refer the reader to \cite{Bombieri70}, \cite{Waldschmidt79}, \cite{Demailly82}, \cite{Gasbarri10}, \cite{Herblot12} for more general higher-dimensional situations and references to related work.

In the sequel, $\Qb$ will denote the algebraic closure of $\Q$ in $\C$ --- or equivalently, an algebraic closure of $\Q$ equipped with  some preferred  embedding in $\C$.

\subsection{Algebraicity of leaves of rank one algebraic foliations}  Let $K$ be a number field, embedded in $\C$, and $X$ a smooth quasi-projective variety over $K$, and let $L \hookrightarrow T_{X/K}$ be a sub-vector bundle of rank 1 of its tangent bundle.

By base field extension from $K$ to $\C$ and analytification, we obtain a complex analytic manifold $X_{\C}^\an$ and an analytic sub-vector bundle $L^\an_{C} \hookrightarrow T_{X_{C}^\an}.$ Since $L^\an_{\C}$ has rank 1, it is integrable (in other words, its sheaf of sections is stable under Lie bracket), and defines a $\C$-analytic foliation of $X_{\C}^\an$. Consider some analytic leaf $\cF$ of this foliation --- it is a connected Riemann surface, equipped with an injective analytic immersion into $X_{\C}^\an$ --- and assume that, for some closed discrete subset $\Delta$ of $\C$, we are given a nonconstant analytic map:
$$f : \C \setminus \Delta \lra \cF.$$

The map $f$ defines an analytic map from $\C \setminus \Delta$ into the quasi-projective complex variety $X_{\C}^\an \hookrightarrow \PP^N(\C)$. As such, it is said to be meromorphic on $\C$ when it extends to an analytic map, which we will still denote $f$, from $\C$ to $\PP^N(\C).$ When this holds, it is said to be \emph{of order} $\leq \rho$ for some $\rho \in \R_{+}$ when, for every $\epsilon >0,$ it admits an analytic lift\footnote{In other words, for every $t \in \C$, $f(t)=(F_{0}(t):\cdots:F_{N}(t)).$}
$$F = (F_{0},\ldots,F_{N}) : \C \lra \C^{N+1}\setminus\{0\}$$
such that
$$\log^+ \max_{0\leq i \leq N} \vert F_{i}(t)\vert = O(\vert t \vert^{\rho + \epsilon}) \mbox{ when $\vert t \vert \rightarrow + \infty$.}$$

Here is a first instance of a transcendence theorem \emph{\`a la} Schneider--Lang (see for instance \cite{Herblot12}, notably Section 6,  for a proof and for a discussion of earlier variants):

\begin{theorem}\label{SL1} Let $K, X, \cF, \Delta$, and $f$ be as above. If

\noindent (1) $f$ is meromorphic of finite order $\leq \rho$, and

\noindent (2) there exists a subset $A$ of $\C\setminus \Delta$ such that $f(A) \subset X(K)$, whose cardinality $\vert A \vert$ satisfies
$$\vert A \vert > 2 \rho [K:\Q],$$
then $\cF$ is algebraic.
\end{theorem}

Here the algebraicity of $\cF$ precisely means that the Riemann surface $\cF$, injectively immersed in $X^\an_{\C}$ is actually a (necessarily closed and smooth) complex algebraic curve in $X_{\C}$. It is equivalent to the algebraicity of the formal germ $\widehat{\cF}_{f(z)}$ of $\cF$ through $f(z)$, for any $z \in A.$ The formal germ $\widehat{\cF}_{f(z)} \hookrightarrow \widehat{X}_{\C,f(z)}$ is indeed defined\footnote{In other words, it is deduced by extension of scalars from $K$ to $\C$ from a formal germ in the formal completion $\widehat{X}_{f(z)}$ of $X$ at the $K$-rational point $f(z)$.}  over $K$, and consequently its Zariski closure in $X_{\C}$  is also. Finally, when conditions (1) and (2) hold, $\cF$ is the set of complex points of some smooth closed $K$-curve in $X.$

Classically a transcendence theorem \emph{\`a la} Schneider--Lang like Theorem \ref{SL1} is rather expressed in the following contrapositive formulation:  \emph{if $f$ is meromorphic of finite order $\rho$ and if $\cF$ is not algebraic, then the cardinality of the subset $f^{-1}(X(K))$ of $\C\setminus \Delta$ is at most $2\rho [K:\Q].$}

A simple but nontrivial instance of Theorem \ref{SL1} arises when
$$X := \A^1 \times \G_{m},$$
$$L:=(\partial/\partial x + y \,\partial/\partial y) \cO_{X}$$
(where $x$ and $y$ denote the standard coordinates on $\A^1 \times \G_{m} \hookrightarrow \A^2$), and $\cF$ is the image of 
$$
\begin{array}{crcl}
f : & \C  & \lra   & X^\an_{\C}  \\
& t  & \lmt   & (t,e^t).  
\end{array}
$$
Clearly $f$ is of order $\leq 1$  and $\cF$ is not algebraic, and Theorem \ref{SL1} asserts that, for any number field $K$ in $\C$, the intersection $f^{-1}(X(K))$ is finite, of cardinality $\leq 2 [K:\Q].$ Besides, if for  some $z$ in $K$, $f(z)$ belongs to $X(K)$, then for any $n\in \Z,$ $f(nz)$ belongs to $X(K)$. Consequently in this case Theorem \ref{SL1} boils down to the \emph{Theorem of Hermite-Lindemann}, which asserts that \emph{for any non-zero complex number $z$, $(z, e^z)$ does not belong to $\Qb^2$.} 

\subsection{Algebraic Lie subalgebras} Let $G$ be a (quasi-projective) algebraic group over $\Qb$, and let $\Lie G$ denote its Lie algebra. Observe that
$$\Lie G_{\C} := \Lie G \otimes_{\Qb}\C \simeq \Lie (G_{\C})$$
may be identified with the Lie algebra of the complex Lie group $G^\an_{\C}.$ In particular, we may consider the exponential map of this Lie group:
$$\exp_{G_{\C}} : \Lie G_{\C} \lra G^\an_{\C}.$$
It is a $\C$-analytic map, \'etale at $0$, and of finite order.

We may also consider the formal variant of this exponential map:
$$\widehat{\exp}_{G}: (\Lie G)^{\wedge}_{0} \lrasim \widehat{G}_{e},$$ 
which is an isomorphism between the formal completion of $\Lie G$ at $0$ --- defined as the formal spectrum of the completion of the symmetric algebra ${\rm Sym}^\bullet(\Lie G)^\vee$,
$$(\Lie G)^{\wedge}_{0} := {\rm Spf}[{\rm Sym}^\bullet(\Lie G)^\vee]^\wedge$$
--- and the formal completion $\widehat{G}_{e}$ of $G$ at its unit element $e$.

A $\Qb$-Lie subalgebra $V$ of $\Lie G$ will be called \emph{algebraic} when the formal subgroup $\widehat{\exp}_{G} V_{0}^\wedge$ that it defines may be algebraized, or equivalently, when  \emph{there exists a $\Qb$-algebraic subgroup $H$ of $G$ such that $V = \Lie H$.}

Transcendence techniques \emph{\`a la} Schneider-Lang may be used to derive ``arithmetic criteria'' for a Lie subalgebra of $\Lie G$ to be algebraic. For instance, when $G$ is commutative --- so that any $\Qb$-vector subspace of $\Lie G$ is a Lie subalgebra --- they lead to the following result, which appears as a vast generalization of Schneider's original result in \cite{Schneider41} (see \cite{Lang66b}, IV, \S 4, Th. 2, when $G$ is a linear group or an abelian variety, and  \cite{Waldschmidt79}, Th. 5.2.1, for a general commutative algebraic group $G$):
 
\begin{theorem}\label{SL2} For any commutative algebraic group $G$ over $\Qb$ and any $\Qb$-vector subspace $V$ of $\Lie G$, the following two conditions are equivalent :

\noindent \emph{(1)} $V$ is an  algebraic Lie subalgebra;

\noindent \emph{(2)} there exists a family $(w_{i})_{i \in I}$ of element of $V_{\C}$ such that, for any $i \in I,$
$$\exp_{G_{\C}} w_{i} \in G(\Qb),$$
which generates the $\C$-vector space $V_{\C}.$

\end{theorem}

The direct implication ${\rm (1)} \Rightarrow {\rm (2)}$ is straightforward. The converse implication ${\rm (2)} \Rightarrow {\rm (1)}$ is a transcendence statement. Consider for instance the case where $G = \G_{m}^2$. Then the (connected) algebraic subgroup of $G$ are defined by monomial equations, and consequently the  algebraic Lie subalgebras $V$ of 
$$\Lie G = \Lie \G_{m}^2= \Qb.x\partial/\partial x  \oplus \Qb.y\partial/\partial y$$
are precisely the $\Qb$-vector subspaces of $\Lie G$ which are $\Q$-rational in the basis $(x\partial/\partial x, y\partial/\partial y)$. Therefore Theorem \ref{SL2} for $G=\G_{m}^2$ becomes the \emph{Theorem of Gelfond-Schneider}, which asserts that \emph{ for any $\alpha$ in $\Qb^\ast$ and any \emph{non-zero} complex number $\log \alpha$ such that $\exp (\log \alpha) = \alpha$, and for any $\beta$ in $\Qb \setminus \Q,$ $\alpha^\beta := \exp (\beta \log \alpha)$ does not belong to $\Qb$.}

Observe also that, when $\dim V =1,$ Theorem \ref{SL2} follows from Theorem \ref{SL1} applied to the translation invariant sub-vector bundle $L$ in $T_{G/\Qb}$ such that $L_{e}= V$. (Choose $K$ large enough to have $G$ and $V$ defined over $K$.) In general, Theorem \ref{SL2} may be seen as an algebraic integrability criterion for translation-invariant algebraic foliations on the algebraic groups $G$.

Let me point out that Theorem \ref{SL2} is now subsumed in stronger transcendence results on commutative algebraic groups, such as the theorems of Baker on linear forms in logarithms and the analytic subgroup theorem of W\"ustholz. The reader may find a recent survey of these results  in the monograph \cite{BakerWuestholz07}.

\subsection{Morphisms of commutative algebraic groups}\label{Morag} In the sequel, we shall use a corollary of Theorem \ref{SL2} which describes morphisms of connected commutative algebraic groups over $\Qb$ in terms of Lie theoretic data. This type of consequence was already pointed out by Bertrand in \cite{Bertrand83},   Section 5, Prop. 2B, where Theorem \ref{SL2} is applied in a similar way to investigate the ring of endomorphisms of a commutative algebraic group.

If $G$ is a connected commutative algebraic group over $\C$, we may introduce its group of ``periods'' 
$$\Per G := \ker \exp_{G},$$
defined as the kernel of its exponential map. It is a discrete subgroup of its Lie algebra $\Lie G$, and fits into an exact sequence of commutative complex Lie groups
$$0 \lra \Per{G}\, \hlra \Lie G \xrightarrow{\exp_{G}} G^\an \lra 0.$$ 

We shall say that \emph{$G$ satisfies Condition} $\mathbf{LP}$ when the group of periods $\Per G$ generates $\Lie G$ as a complex vector space.  

Observe that this condition is preserved by isogenies, and by forming quotients and products, and is satisfied by the multiplicative group $\G_{m\C}$, complex abelian varieties, and universal vector extensions. Actually, a connected commutative algebraic group $G$ over $\C$
satisfies Condition $\mathbf{LP}$ precisely when $G$ is ``almost semi-abelian" or ``anti-additive" in the sense of \cite{BertrandPillay10}, Section 3.1, namely when the torsion points of $G(\C)$ are Zariski dense in $G$, or equivalently when there is no nontrivialmorphism of algebraic groups from $G$ to the additive group $\G_{a\C}$ (\cf \emph{loc. cit.}, Appendix I).  In particular condition $\mathbf{LP}$ is a purely algebraic condition, invariant under the automorphisms of the field $\C$.  

\begin{corollary}\label{CorSL2} Let $G_{1}$ and $G_{2}$ be connected commutative algebraic groups over $\Qb$. 

1) For any $\phi$ in the $\Z$-module $\Hom_{\ggp/\Qb}(G_{1}, G_{2})$ of morphisms of  algebraic groups over $\Qb$ from $G_{1}$ to $G_{2}$, the $\Qb$-linear map
$$\Lie \phi := D\phi (e) : \Lie G_{1} \lra \Lie G_{2}$$
satisfies
$$(\Lie \phi)_{\C}(\Per G_{1\C}) \subset \Per G_{2\C}.$$

The  map
\begin{equation}\label{LieGamma}
\Lie : \Hom_{\ggp/\Qb}(G_{1}, G_{2}) \lra \{ \psi \in \Hom_{\Qb}(\Lie G_{1},\Lie G_{2})\vert \psi_{\C}(\Per G_{1\C}) \subset \Per G_{2\C} \}
\end{equation} so defined
 is an injective morphism of $\Z$-modules.

2) When the group $G_{1\C}$ satisfies condition 
$\mathbf{LP}$, then the morphism (\ref{LieGamma}) is bijective.
\end{corollary}

\Proof{}
 Assertion 1) follows from identification of   $(\Lie \phi)_{\C}$  with the differential $\Lie \phi_{\C}:= D\phi_{\C}(e)$ of the complexification $\phi_{\C}: G_{1\C}\rightarrow G_{2\C}$ of the morphism of $\Qb$-algebraic groups $\phi,$ together with the commutativity of the diagram:
 $$
 \begin{CD}
 \Lie G_{1\C} @>{\Lie \phi_{\C}}>> \Lie G_{2\C} \\
@V{\exp_{G_{1\C}}}VV                @VV{\exp_{G_{2\C}}}V \\
G_{1\C}^\an @>{\phi_{\C}}>>       G_{2\C}^\an.
\end{CD}
$$

To prove 2), assume that condition  $\mathbf{LP}$ is satisfied by $G_{1\C},$ and consider some $\Qb$-linear map
$$\psi : \Lie G_{1} \lra \Lie G_{2}$$
such that $\psi_{\C}(\Per G_{1\C}) \subset \Per G_{2\C}$. We need to establish the existence of a morphism of $\Qb$-algebraic groups 
$\phi: G_{1} \lra G_{2}$
such that 
\begin{equation}\label{psiphi}
\psi = \Lie \phi.
\end{equation}

To achieve this, we will apply Theorem \ref{SL2} to the group $G := G_{1} \times G_{2}$, and to the subspace  $V$ of 
$$\Lie G = \Lie G_{1} \oplus \Lie G_{2}$$
defined as the graph of $\psi$.  

Indeed, as $G$ is commutative, $V$ is  a Lie subalgebra of $\Lie G$. Moreover the complex vector space $V_{\C}$ is the graph of $\psi_{\C}$ and therefore contains 
$$\widetilde{\Per G_{1\C}} := \{ (\gamma, \psi_{\C}(\gamma)), \gamma \in \Per G_{1\C} \},$$ 
which is included in $\Per G_{1\C} \times \Per G_{2\C} = \Per G_{\C}$. Besides, the condition $\mathbf{LP}$   on $G_{1\C}$ shows that $\widetilde{\Per G_{1\C}}$ generates this $\C$-vector space. According to Theorem \ref{SL2}, $V$ is  algebraic and is the Lie algebra of some connected $\Qb$-algebraic subgroup $H$ of $G$.

The first projection $p :={\pr}_{1\vert H} : H \lra G_{1}$ is \'etale. Moreover $H_{\C}^\an$ is the image by $\exp_{G_{\C}}$ of $V_{\C}$. This immediately implies that  ${p}_{\C} : H_{\C} \lra G_{1\C}$ is injective, and finally that $p$ is an isomorphism. In other words, $H$ is the graph of some morphism $\phi$ of algebraic groups from $G_{1}$ to $G_{2}$. Clearly it satisfies (\ref{psiphi}). \BOX{}

\subsection{Transcendence theorems and the analogy between number fields and functions fields}\label{subsec:analogy} 

Theorems \ref{SL1} and \ref{SL2} may be seen as arithmetic counterparts of algebraization theorems such as Andreotti's Theorem \ref{thAndreotti}, or  $\mathbf{Gr1}$ in Theorems \ref{GrGr} and \ref{GrGrfor}, or more specifically, of their consequences concerning algebraization over function fields, such as Corollary \ref{corthA} and its formal variant. The role of the function field $\C(C)$ or $k(C)$ is now played by $\Qb$ or by a number field $K$ over which the geometric data $X$ and $L$, or $G$ and $V$, are defined. 

Observe that the so-called Kronecker dimension of $K$ --- namely the Krull dimension of $\Spec \OK$ --- is \emph{one}, and that the algebraization Theorems \ref{SL1} and \ref{SL2}, which are algebraicity criteria for smooth formal germs of subvarieties through $K$-rational \emph{points}, isomorphic to $\Spec K$, are indeed algebraization theorems concerning smooth formal germs along some \emph{arithmetic curves} $\Spec \OK$ in some integral model of the given $K$-variety.  

The classical proofs of Theorems \ref{SL1} and \ref{SL2} may be understood in a way that makes this geometric analogy precise. This geometric approach even suggests the formulation and the proof of new  transcendence theorems, as demonstrated by the  recent works of Gasbarri \cite{Gasbarri10} and Herblot \cite{Herblot12} who have established sophisticated generalizations of previously known transcendence theorems \emph{\`a la} Schneider-Lang. I might also refer the reader to  \cite{Chambert01} and \cite{Bost06} for discussions of this geometric approach and of some of its applications in the framework of  Diophantine results \emph{\`a la} Chudnovsky (\cite{ChudnovskysGroth85}, \cite{ChudnovskysAcad85}) instead of Schneider-Lang. The arithmetic counterparts of the ampleness conditions   in the geometric theorem of Andreotti-Hartshorne  and $\mathbf{Gr1}$ appear more clearly in this somewhat simpler framework.

At the present stage, in this analogy, there is no known counterpart in transcendence theory of the general Existence Theorems, such as $\mathbf{Gr2}$ in Theorems \ref{GrGr} and \ref{GrGrfor}. This absence appear especially regrettable when one considers the important geometric applications of these theorems: we have discussed at length several consequences of GAGA Existence Theorem in Sections \ref{subsec:algline}, \ref{subsec:algebdeRham}, and \ref{subscec:alganstructures};  as demonstrated in \cite{GrothendieckSGA2}, $\mathbf{Gr2}$ is the key to a modern approach to ``Lefschetz-type theorems'' which compare invariants, such as their fundamental group or their Picard group, of projective varieties to the ones of their hyperplane section.

The dimension condition 
$$\dim Y \geq 2$$
in $\mathbf{Gr2}$ leads one, in a Kroneckerian perspective, to expect a suitable arithmetic counterpart of  $\mathbf{Gr2}$ to be an algebraization criterion concerning formal line or vector bundles over the completion $\widehat{X}_{Y}$ of some algebraic variety $X$ over a number field $K$, along a smooth projective embedded curve $Y$ over $K$, or if one prefers, over the completion $\widehat{\cX}_{\cY}$ of some scheme of finite type $\cX$ over $\Spec \OK$ along a projective arithmetic surface $\cY.$ 

In the spirit of transcendence theorems \emph{\`a la} Schneider--Lang like Theorems \ref{SL1} and \ref{SL2}, this criterion would also require some ``differential algebraic'' conditions (comparable to the occurrence of algebraic foliations in these theorems) and some ``analytic control'' on the considered formal vector bundles.

The remainder of this article is devoted to presenting such a criterion, in a conjectural form, and its relation to Grothendieck Period Conjecture in codimension 1. 

The proof of this last conjecture for abelian varieties may actually be derived from Theorem \ref{SL2} and its Corollary \ref{CorSL2}. As it provides a further  illustration of the ``concrete geometric content'' of transcendence theorems \emph{\`a la} Schneider--Lang, we begin by a discussion of this material in Part \ref{GPCAb}. Then, in Sections \ref{basicD} to \ref{ExtD}, we review the formalism of $D$-group schemes and of their extensions that will be used in the last part to formulate our conjectural algebraization criterion.  

\section{The Grothendieck Period  Conjecture for cycles of codimension 1  in abelian varieties}\label{GPCAb}

\subsection{Grothendieck's conjecture $GPC^1(X)$}\label{GPC}

Let $X$ be a smooth projective algebraic variety over $\Qb,$ and let $X_{\C}$  denote the smooth complex projective variety $X \otimes_{\Qb}\C,$  and $X^\an$ the corresponding compact complex manifold. 

As discussed in Section \ref{subsec:algebdeRham}, the Picard groups of $X$, $X_{\C}$, and $X_{\C}^\an$ --- which classify the algebraic lines bundles over $X$ and $X_{\C}$, and the analytic line bundles over $X_{\C}^\an$ --- fit into the following commutative diagram:
$$\label{bigdiag}
 \begin{CD}
 \Pic(X)   @>{c_{1\dR/\Qb}}>> \HdR^2(X/\Qb) \\
 @VVV                                  @VV{.\otimes_{\Qb}1_{\C}}V \\ 
 \Pic(X_{\C}) @>{c_{1\dR/\C}}>>             \HdR^2(X_{\C}/\C) \\ 
 @VV{.^\an}V                                                @VV{.^\an}V \\
 \Pic (X_{\C}^\an)  @>{c_{1\dR}^\an}>>                                  \HdR^2(X^\an_{\C}/\C) \\
 @VV{c_{1\topo}}V                                       @VV{\text{de Rham isomorphism}}V \\
 H^2(X^\an_{\C}, \Z) @>{2\pi i (. \otimes_{\Z}1_\C)}>> H^2(X_{\C}^\an,\C).
\end{CD}
$$

The upper vertical arrows are induced by the field extension $\Qb \hookrightarrow \C$. The map $\Pic(X) \lra \Pic(X_{\C})$ maps the class of some line bundle $L$ over $X$ to the class of the line bundle  $L_{\C}$ over $X_{\C}$, and is injective, but not surjective when the connected Picard variety $\Pic_{0}(X/\Qb)$ has positive dimension\footnote{that is,  when the ``irregularity" $h^{1,0} (X) = h^{0,1} (X)$ of $X$ is positive.}. However, since any line bundle over $X_{\C}$ is algebraically equivalent to some line bundle defined over $\Qb$, the images of  $\Pic(X)$ and $\Pic(X_{\C})$ by the first Chern class coincide. The map  $\HdR^2(X/\Qb) \lra \HdR^2(X_{\C}/\C)$ induces an isomorphism $\HdR^2(X/\Qb) \otimes_{\Qb}\C \lrasim \HdR^2(X_{\C}/\C)$. The image in $\HdR^2(X_{\C}/\C)$ of an element $\alpha$ in $\HdR^2(X/\Qb)$ will be denoted $\alpha \otimes_{\Qb}1_{\C}$.

The two middle vertical arrows  $.^\an$, defined by analytification, are isomorphisms according to GAGA. The analytification isomorphism $\HdR^2(X_{\C}/\C) \lrasim \HdR^2(X^\an_{\C}/\C)$ will be noted as an equality. 

The image of some class $\beta \in H^2(X^\an_{\C}, \Z)$ by the natural map $H^2(X^\an_{\C}, \Z) \lra H^2(X^\an_{\C}, \C)$ (defined by extending the coefficients from $\Z$ to $\C$) will be denoted $\beta \otimes_{\Z} 1_\C$, and the image of some class $\gamma$ in $\HdR^2(X^\an_{\C}/\C)$ by the de Rham isomorphism will be denoted $\gamma^{\rm B}$.

We may  define the subgroup
$H^2_{\rm Gr}(X)$ of ``Grothendieck's classes'' in $H^2_{\dR}(X/\Qb) \oplus H^2(X_{\C}^\an, \Z)$
by the condition that, for any $\alpha \in \HdR^2(X/\Qb)$ and any $\beta \in H^2(X_{\C}^\an, \Z)$:
\begin{equation}\label{defGr} (\alpha,\beta) \in H^2_{\rm Gr}(X) \Longleftrightarrow (\alpha \otimes_{\Qb} 1_{\C})^{\rB}  = 2\pi i \, \beta\otimes_{\Z}1_{\C}.
\end{equation}
The commutativity of the diagram above 
shows that the algebraic and topological first Chern classes define a morphism of abelian groups:
$$
\begin{array}{crcl}
c_{1\dRB}: & \Pic(X)  & \lra  & H^2_{\rm Gr}(X)  \\
& [L]  & \lmt  & (c_{1\dR}(L), c_{1\topo}(L^\an_{\C})).
  \end{array}
$$

The classical Grothendieck Period Conjecture\footnote{This conjecture is mentioned briefly  in \cite{Grothendieck66} (note (10) p.102) and with more details in \cite{Lang66b} (Historical Note of Chapter IV). We refer the reader to \cite{AndreMotives04}, Section 7.5 and Chapitre 23 for a ``modern'' presentation and for variants and generalizations.} 
 leads one to conjecture that \emph{the morphism $c_{1\dRB}$ is onto}, namely that \emph{a class $\gamma$ in $H^2(X^\an_{\C}, \Z)$ such that $2 \pi i. \gamma \otimes_{\Z}1_{\C}$ is $\Qb$-rational in 
$$H^2(X^\an_{\C},\C) \simeq \HdR^2(X/\Qb) \otimes_{\Qb} \C$$
is algebraic} in the sense of Section \ref{subsec:algline}.

This conjectural assertion assertion may be called \emph{the Grothendieck Period Conjecture in codimension 1} for the smooth projective variety $X$ over $\Qb$ and will be denoted $GPC^1(X)$ in the sequel. 



Conjecture $GPC^1(X)$ admits a $\Q$-rational version, \emph{a priori} weaker, that asserts the surjectivity of the map
$$c_{1\dRB\Q}:  \Pic(X)_{\Q}   \lra   H^2_{\rm Gr}(X)_{\Q}$$
deduced from $c_{1\dRB}$ by tensoring with $\Q$. (The tensor product $H^2_{\rm Gr}(X)_{\Q}:= H^2_{\rm Gr}(X)\otimes{\Q}$ may be identified with the $\Q$-vector subspace of 
$H^2_{\dR}(X/\Qb) \oplus H^2(X_{\C}^\an, \Q)$  defined by the right-hand side of (\ref{defGr}), with $.\otimes_{\Z}.$ replaced by $.\otimes_{\Q}.$) A special feature of the codimension 1 case of the Grothendieck Period Conjecture is that this rational version of the conjecture --- which is the one that appears in \emph{loc. cit.} --- actually implies the above ``integral'' version. Indeed, for any positive integer $n$, a class $\gamma$ in $H^2(X_{\C}^\an, \Z)$ is algebraic if $n\gamma$ is algebraic.

More generally, for any positive integer $k$, we may consider the Grothendieck Period Conjecture in codimension $k$, $GPC^k(X)$ : it asserts that any class $\gamma$ in $H^{2k}(X^\an_{\C}, \Q)$ such that $(2\pi i)^k \gamma \otimes_{\Q}1_{\C}$ is $\Qb$-rational in $H^{2k}(X^\an_{\C},\C) \simeq \HdR^{2k}(X/\Qb) \otimes_{\Qb} \C$ is algebraic. See \cite{AndreMotives04}, Section 7.5, for  a discussion of the  close relation between the original version of the Grothendieck Period Conjecture and the  fullness conjecture for the ``de Rham--Betti realization'', namely  the conjunction of  Conjectures $GPC^k(X)$ for all smooth projective varieties $X$ over $\Qb$ and all integers $k$\footnote{Notably the original Grothendieck Period Conjecture for a given smooth projective variety $X$ over $\Qb$ should imply the conjunction of Conjectures $GPC^k(X^n)$ for all positive integers $k$ and $n$.}.  To my knowledge, the known results concerning these conjectures may be summarized as follows : 

(i)  the original Grothendieck Period Conjecture is known to be valid  for a motive  in the Tannakian category generated by the Tate motive (transcendence of $\pi$) or for an elliptic curve with complex multiplication (Chudnovsky); 

(ii)   the fullness of the de Rham-Betti realization is known for $H^1$ (\cf \cite{AndreMotives04}, 7.2.3, where it is derived from the transcendence results in \cite{Wuestholz84}; this fullness is basically  the content of Theorem \ref{HomAB} \emph{infra}, and as shown in the next paragraphs, may be derived from Schneider-Lang's Theorem \ref{SL2} and its Corollary \ref{CorSL2}).

In the next sections, we shall establish the validity of Grothendieck's Period Conjecture in codimension  1 for abelian varieties:

\begin{theorem}\label{GPCA}
For any abelian variety $A$ over $\Qb,$ $GPC^1(A)$ holds.
\end{theorem}

The proof of Theorem \ref{GPCA} will be based on the ``transcendental'' characterization of  algebraic Lie subalgebras in Theorem \ref{SL2}, via its Corollary \ref{CorSL2} applied to universal vector extensions of abelian varieties, and on the identification of the N\'eron-Severi group of an abelian variety with the group of symmetric morphisms from the abelian variety to its dual (compare \cite{Bost06}, Theorem 6.4). We present the details of this proof in Section \ref{GPCAb}. As a preliminary, in Section \ref{PrelAb} we recall classical facts concerning abelian varieties, their duality, and their universal vector extensions, and in Section \ref{CDRB} we introduce the elementary, but convenient, formalism of the category $\CdRB$ of the ``de Rham--Betti realisations'' (in the spirit of the realisation categories \emph{\`a la} Deligne--Jannsen \cite{Jannsen90}; see also \cite{AndreMotives04}, Section 7.5.).

\subsection{Abelian varieties, duality, and universal extensions}\label{PrelAb}

In this section, we work over an algebraically closed field $k$ of characteristic zero.

\subsubsection{Dual abelian varieties and de Rham (co)homology}\label{subsubPoinc} 
If $A$ is an abelian variety over $k$, we shall denote $\hA := \Pic_{0}(A/k)$ the dual abelian variety. The group $\hA(k)$ of its $k$-rational points may be identified with the subgroup $\Pic^{0}(A)$ of $\Pic(A)$ of isomorphism classes of line bundles algebraically equivalent to zero, or equivalently, with the kernel of 
$$c_{1\dR}: \Pic(A) \lra \HdR^2(A/k).$$

To any morphism $\phi : A \lra B$ of abelian varieties over $k$ is attached the dual morphism $\hat{\phi} :\hB \lra \hA.$ It maps the class of some line bundle $L$ over $B$ algebraically equivalent to zero to the class of $\phi^\ast(L).$ This construction is additive and (contravariantly) functorial. 

Let $\cP_{A}$ denote the Poincar\'e line bundle over $A \times \hA$. Its restriction to $0_{A}\times \hA$ is trivial, and for any $\ha \in \hA(k),$ the isomorphism class of its restriction to $A \times \ha$ is precisely $\ha$ itself, and these properties characterize $\cP_{A}$ up to isomorphism. By mapping a point $a$ in $A(k)$ to the class $\iota_{A}(a)$ of $\cP_{A\vert a \times \hA}$, ones defines a canonical isomorphism 
$$\iota_{A} : A \lrasim \hat{\hA},$$
which is sometimes written as an equality.

Recall that the following ``biduality'' properties are satisfied (compare \cite{BerthelotBreenMessing82}, Section V.1, or \cite{Coleman91}, Section 1). For any $\phi : A \lra B$ as above, $\hat{\hat{\phi}} :\hat{\hA} \lra \hat{\hat{B}}$ and $\phi$ (or more exactly $\iota_{B}\circ \phi \circ \iota_{A}$) coincide. Moreover, under the composite isomorphism
$$
\begin{array}{rclc}
 A \times \hA & \stackrel{\sigma}{\lra} & \hA \times A & \xrightarrow{Id_{\hA}\times\iota_{A}} \hA \times \hat{\hA}  \\
  (a,\ha) & \lmt   & (\ha,a) &   
\end{array}
$$
the Poincar\'e bundle $\cP_{A}$ of $A$ becomes the Poincar\'e bundle $\cP_{\hA}$ of $\hA$:
\begin{equation}\label{Poincdual}
((Id_{\hA}\times \iota_{A})\circ \sigma)^\ast \cP_{\hA} \lrasim \cP_{A}.
\end{equation}

Moreover $c_{1\dR}(\cP_{A})$ belongs to the K\"unneth component $H^1_{\dR}(A/k) \otimes H^1_{\dR}(\hA/k)$ of $H^2(A\times \hA/k).$ If we define
$$H_{1\dR}(A/k) := H^1_{\dR}(A/k)^\vee = \Hom_{k}(H^1_{\dR}(A/k),k),$$
then $c_{1\dR}(\cP_{A})$ defines an element $\varpi_{A}$ in
$$H_{1\dR}(A/k)^\vee \otimes_{k} \HdR^1(\hA/k) \simeq \Hom_{k}(H_{1\dR}(A/k),\HdR^1(\hA/k))$$
which actually is an isomorphism:
$$\varpi_{A} : H_{1\dR}(A/k) \lrasim \HdR^1(\hA/k) = H_{1\dR}(\hA/k)^\vee.$$

The duality isomorphism $\varpi_{A}$ satisfies the following functoriality property. 

Let $\phi: A \lra B$ be a morphism of abelian varieties over $k$. It induces a $k$-linear map between  de Rham cohomology groups:
$$H^1_{\dR}(\phi) :=\phi^\ast : H^1_{\dR}(B/k) \lra H^1_{\dR}(A/k),$$
and then by duality, between homology groups:
$$H_{1\dR}(\phi) := H^1_{\dR}(\phi)^t : H_{1\dR}(A/k) \lra H_{1\dR}(B/k).$$
Then the dual morphism of abelian varieties
$$\hat{\phi}: \hB \lra \hA$$
satisfies
\begin{equation}\label{dualdual}
H_{1\dR}(\hat{\phi}) = \varpi_{A}^{\vee -1} \circ H_{1}(\phi)^\vee \circ \varpi_{B}^\vee.
\end{equation}
This follows from the isomorphism of line bundles over $A \times \hat{B}$:
$$(Id_{A}\times \hat{\phi})^\ast \cP_{A} \simeq (\phi \times Id_{\hB})^\ast \cP_{B},$$
and from the implied equality between first Chern classes.

Observe however that the isomorphism 
$$\varpi_{\hA}: H_{1\dR}(\hA/k) \lrasim H^1_{\dR}(\hat{\hA}/k) \simeq H_{1\dR}(\hat{\hA}/k)^\vee$$
differs by a sign from the transpose of $\varpi_{A}$:
\begin{equation}\label{sign}
\varpi_{\hA} = - H_{1\dR}(\iota_{A})^\vee \circ \varpi_{A}^\vee.
\end{equation}
This follows from the equality of first Chern classes implied by the isomorphism (\ref{Poincdual}), and from the fact that switching the factors $A \simeq \hat{\hA}$ and $\hA$ introduces a sign in the K\"unneth morphism 
$$H^1_{\dR}(A/k) \otimes_{k} H^1_{\dR}(\hA/k) \hlra H^2_{\dR}(A\times\hA/k).$$

\subsubsection{N\'eron-Severi groups and symmetric morphisms}\label{NSsym} To any line bundle $L$ over $A$ is attached a morphism of abelian varieties over $k$,
$$\phi_{L}: A \lra \hat{A},$$
that is defined by 
$$\phi_{L}(a) := [\tau_{a}^\ast L \otimes L^\vee]$$
for any $a \in A(k)$, where $\tau_{a}$ denotes the translation by $a$ on $A$. Moreover $\phi_{L}$ is zero if and only if $L$ is algebraically equivalent to zero, and, for any two line bundles $L_{1}$ and $L_{2}$ on $A,$ $\phi_{L_{1}\otimes L_{2}}=\phi_{L_{1}}+\phi_{L_{2}}.$ Consequently this construction induces an injective morphism of $\Z$-modules:
$$
\begin{array}{rcl}
 NS(A) := \Pic(A)/\Pic_{0}(A) & \lra  & \Hom_{\gp/k}(A,\hA)  \\
{[L]} & \lmt   & \phi_{L} . 
\end{array}
$$ 
Its image is the subgroup $\Hom_{\gp/k}(A,\hA)^{\text{sym}}$ of \emph{symmetric} morphisms, namely the subgroup of morphisms $\phi : A \lra \hA$ such that
\begin{equation}\label{symphi}
\hat{\phi} \circ \iota_{A} = \phi.
\end{equation}

This actually holds for abelian schemes over an arbitrary base, as established by Nishi and Oda (\cf \cite{Oda69}, p. 77, note $(^2)$).

Observe that, at the level of de Rham (co)homology groups, the symmetry condition (\ref{symphi}) translates into a \emph{skew-symmetry} condition on 
$$\varpi_{A}^\vee \circ H_{1}(\phi) : H_{1\dR}(A/k) \lra H_{1\dR}(A/k)^\vee.$$
Indeed   the ``duality'' formulas (\ref{dualdual}) and (\ref{sign}) imply the relation:
\begin{equation}\label{altphi}\varpi_{A}^\vee \circ H_{1}(\hat{\phi}\circ \iota_{A}) = - (\varpi_{A}^\vee \circ H_{1}(\phi))^\vee.
\end{equation}

In particular, when the base field $k$ is $\C$, the above identification of $NS(A)$ with $\Hom_{\gp/k}(A,\hA)^{\text{sym}}$ is basically the classical theory of Riemann forms attached to line bundles over complex abelian varieties.

\subsubsection{Universal vector extensions} (\cf \cite{Rosenlicht58}, \cite{Serre59}, \cite{Messing73}, \cite{MazurMessing74}, \cite{Coleman91}, \cite{BK09}). 

For any abelian variety $A$ over $k$, we shall denote $\E_{A}$ the $k$-vector space
$$\Gamma(A,\Omega^1_{A/k}) \simeq \Omega^1_{A/k, 0_{A}} \simeq (\Lie A)^\vee.$$
Observe that we have a canonical identification
$$\E_{\hA} \simeq (\Lie \hA)^\vee \simeq H^1(A,\cO_{A})^\vee.$$

Let $V$ a finite-dimensional $k$-vector space, and let $V^{\gp}$ denote the associated $k$-vector group (namely the commutative algebraic group over $K$, such that the group $V^{\gp}(k)$ ``is'' the additive group $(V,+)$). Recall that any extension of commutative algebraic groups over $k$
\begin{equation}\label{VGA}
0 \lra V^{\gp} \lra G \lra A \lra 0
\end{equation}
of some abelian variety $A$ over $k$ by $V^\gp$ determines an $\cO_{A}\otimes_{k}V$-torsor over $A$, and that this construction defines a canonical isomorphism\footnote{Where $\Ext^1_{\cgp/k}$ and $\Ext^1_{\cO_{A}-\text{mod}}$ stand for ``group of 1-extensions'' in the category of commutative algebraic groups over $k$, and of sheaves of $\cO_{A}$-modules, respectively.}
\begin{equation}\label{canext}
\Ext^1_{\cgp/k}(A, V^\gp) \lrasim \Ext^1_{\cO_{A}-\text{mod}}(\cO_{A}, \cO_{A}\otimes_{k}V) \simeq H^1(A,\cO_{A}) \otimes_{k} V \simeq \Hom_{k}(\E_{\hA},V).
\end{equation}
Moreover an extension (\ref{VGA}) of commutative algebraic groups of an abelian variety by a vector group admits no nontrivial automorphism. Consequently the isomorphism (\ref{canext}) with $V= \E_{\hA}$ 
shows that, to the element
$Id_{\E_{\hA}}$
is canonically associated a vector extension of $A$ by the vector group defined by $\E_{\hA}$, which we shall denote 
\begin{equation}\label{defEA}
0 \lra \E_{\hA} \hlra E(A) \stackrel{p_{A}}{\lra} A \lra 0.
\end{equation}
It is the \emph{universal vector extension} of $A$ : any vector extension (\ref{VGA}) may be realized uniquely as a pushout of (\ref{defEA}), namely, as the pushout by its ``classifying element'' in the right-hand side of (\ref{canext}).

\subsubsection{The functor $E$}\label{functE} Let $\phi : A \lra B$ be a morphism of abelian varieties over $k$. We may consider the pullback by $\phi$ of the universal vector extension of $B$, and use the universal property of the universal vector extension of $A$. We thus get the existence and unicity of a morphism $E(\phi)$ of $k$-algebraic groups which makes the following diagram commutative:
$$\begin{CD}
E(A) @>{E(\phi)}>> E(B) \\
@VV{p_{A}}V                      @VV{p_{B}}V \\
A  @>{\phi}>>  B.
\end{CD}
$$ 
The construction of $E(\phi)$ is clearly additive and functorial in $\phi.$  Moreover it is easily seen to be fully faithful: 

\begin{lemma} For any two abelian varieties $A$ and $B$ over $k$, the morphism of $\Z$-modules
\begin{equation}\label{Eff}
\begin{array}{rcl}
\Hom_{\gp/k} (A, B)  & \lra   & \Hom_{\gp/k} (E(A), E(B))   \\
 \phi  & \lmt  &   E(\phi).
\end{array}
\end{equation}
is an isomorphism.
\end{lemma}

\subsubsection{Biduality and universal vector extensions}\label{biduve} We shall also use that the biduality isomorphism
$$\iota_{A} : A(k) \lrasim \hat{\hA}(k) = \ker c_{1\dR}: H^1(\hA, \cO_{\hA}^\ast) \lra H^1_{\dR}(\hA, \Omega^\bullet_{\hA/k})
$$
may be lifted to an isomorphism
$$\iota_{E(A)} : E(A)(k) \lrasim H^1(\hA, \Omega^\times_{\hA/k}),$$
where $\Omega^\times_{\hA/k}$ denotes the complex
$$\cO^\ast_{\hA} \stackrel{d\log}{\lra} \Omega^1_{\hA/k} \stackrel{d}{\lra} \Omega^2_{\hA/k}\stackrel{d}{\lra} \cdots,
$$
which makes commutative the following diagram with exact lines\footnote{Recall that $\sigma^{\geq 1} \Omega^\bullet_{\hA/k}$ denotes the ``stupid'' truncation $0\rightarrow \Omega^1_{\hA/k} \rightarrow \Omega^2_{\hA/k}\rightarrow \cdots$ of $\Omega^\bullet_{\hA/k}$.}:
\begin{equation}\label{horrible}
\begin{CD}
0 @>>> \E_{\hA}      @>>>     E(A)(k) @>{p_{A}}>> A(k) @>>> 0 \\
@.        @V{\simeq}VV          @V{\simeq}V{\iota_{E(A)}}V                            @V{\simeq}V{\iota_{A}}V   @. \\
0 @>>> H^1(\hA, \sigma^{\geq 1} \Omega^\bullet_{\hA/k}) @>>> H^1(\hA,  \Omega^\times_{\hA/k}) 
@>>> \hat{\hA}(k) 
@>>> 0. \\
\end{CD}
\end{equation}
(For constructing the second line, recall that  $F^1H^2_{\dR}(\hA/k) 
:= H^2(\hA, \sigma^{\geq 1} \Omega^\bullet_{\hA/k})
$ injects into $H^2_{\dR}(\hA/k)$, and that $c_{1\dR}:  H^1(\hA, \cO_{\hA}^\ast) \rightarrow H^1_{\dR}(\hA, \Omega^\bullet_{\hA/k})$ coincides with $d\log: H^1(\hA, \cO_{\hA}^\ast) \rightarrow F^1H^2_{\dR}(\hA/k).$) 

Moreover the ``infinitesimal'' version\footnote{Both the above isomorphism $\iota_{E(A)}$ at the level of $k$-points and this infinitesimal version are special instances  of a canonical isomorphism $\iota_{E(A)}$ of fpqc $k$-sheaves; \cf \cite{MazurMessing74}, \cite{BK09}.} of $\iota_{E(A)}$ defines an isomorphism
$$I_{A}:= \Lie \iota_{E(A)} : \Lie E(A) \lra H^1(\hA, \Omega^\bullet_{\hA/k}) = H^1_{\dR}(\hA/k),$$
and the infinitesimal version of (\ref{horrible}) is an isomorphism of exact sequences of finite-dimensional $k$-vector spaces:
\begin{equation}\label{moinshorrible}
\begin{CD}
0 @>>> \E_{\hA} @>>> \Lie E(A) @>{\Lie p_{A}}>> \Lie A @>>> 0 \\
@.          @V{=}VV           @V{\simeq}V{I_{A}}V                       @V{\simeq}V{\Lie \iota_{A}}V   @. \\
0 @>>> \E_{\hA} @>>>  H^1_{\dR}(\hA/k) @>>> H^1(\hA, \cO_{\hA}) @>>> 0.\\
\end{CD}
\end{equation}
(The second line defines the Hodge filtration on $H^1_{\dR}(\hA/k)$.) 

Finally we get an isomorphism of $k$-vector spaces
$$J_{A}:=\varpi_{A}^{-1}\circ I_{A} : \Lie E(A) \lrasim H_{1\dR}(A/k).$$
It is easily checked to be functorial. Namely, for any morphism $\phi: A \lra B$ of abelian varieties over $k$, the diagram
\begin{equation}
\begin{CD}
\Lie E(A) @> {\Lie E(\phi)}>> \Lie E(B) \\
@V{\simeq}V{J_{A}}V                      @V{\simeq}V{J_{B}}V \\
H_{1\dR}(A/k)  @>{H_{1\dR}(\phi)}>> H_{1\dR}(B/k)
\end{CD}
\end{equation}
 is commutative.
 
\subsection{The category $\CdRB$}\label{CDRB}

\subsubsection{Definitions} We define an additive category $\CdRB$ --- where $\mathcal{C}$ stands for ``category'' or ``comparison'' and $\dRB$ stands for ``de Rham -- Betti'' --- in the following way.

Its objects are triples 
$$M = (M_{\dR}, M_{\rB}, c_{M}),$$
where $M_{\dR}$ is a finite-dimensional $\Qb$-vector space, $M_{\rB}$ a free $\Z$-module of  finite rank, and $c_{M}$ an isomorphism of $\C$-vector spaces:
$$c_{M}: M_{\dR}\otimes_{\Qb}\C \lrasim M_{\rB}\otimes_{\Z}\C.$$

In other terms, an object $M$ of $\CdRB$ may be seen as the data of the finite-dimensional $\C$-vector space
$$M_{\C} := M_{\dR}\otimes_{\Qb}\C \simeq M_{\rB}\otimes_{\Z}\C,$$
together with a ``$\Qb$-form'' $M_{\dR}$ and a ``$\Z$-form'' $M_{\rB}$ of $M_{\C}$.

If $M$ and $N$ are objects in $\CdRB$, the additive group of morphisms from $M$ to $N$ in $\CdRB$ is the subgroup $\Hom_{\dRB}(M,N)$ in
$\Hom_{\Qb}(M_{\dR},N_{\dR}) \oplus \Hom_{\Z}(M_{\rB},N_{\rB})$ consisting of pairs of maps $\phi = (\phi_{dR},\phi_{\rB})$ such that the following diagram is commutative:
$$\begin{CD}
M_{\dR}\otimes_{\Qb}\C @>{\phi_{\dR}\otimes_{\Qb}Id_{\C}}>> N_{\dR}\otimes_{\Qb}\C \\
@V{\simeq}V{c_{M}}V                                                                                  @V{\simeq}V{c_{N}}V \\
M_{\rB}\otimes_{\Z}\C  @>{\phi_{\rB}\otimes_{\Z}Id_{\C}}>> N_{\rB}\otimes_{\Z}\C.            
\end{CD}
$$

These morphisms may be identified with the $\C$-linear maps $\phi_{\C}$ from $M_{\C}$ to $N_{\C}$ which are compatible both to their $\Qb$-forms and their $\Z$-forms.
The composition of these morphisms is the obvious one, defined by the composition of the ``de Rham'', ``Betti'', and ``complex'' realizations $\phi_{\dR},$ $\phi_{\rB}$, and $\phi_{\C}$ respectively.
 
The category $\CdRB$ is endowed with an internal tensor product, defined by $$M\otimes N := (M_{\dR}\otimes_{\Qb}N_{\dR},M_{\rB}\otimes_{\Z}N_{\rB}, c_{M}\otimes_{\C}c_{N}),$$ and with an internal duality functor, defined by $$M^\vee := (\Hom_{\Qb}(M_{\dR},\Qb), \Hom_{\Z}(M_{\rB}, \Z), c^t),$$ and $$\phi^\vee := (\phi_{\dR}^{t},\phi_{\rB}^t) = (. \circ \phi_{\dR}, . \circ \phi_{\rB}).$$

For any integer $k$, we denote $\Z(k)$ the object of $\CdRB$ defined by $\Z(k)_{\Qb}= \Qb$ and $\Z(k)_{\rB} = (2\pi i)^k \Z$ in $\Z(k)_{\C}= \C$. Observe that $\Z(0)$ and the obvious isomorphism $\Z(0)\otimes\Z(0) \lrasim \Z(0)$, mapping $1\otimes 1$ to $1$, define a unit object of $\CdRB$, which, endowed with $\otimes$ and $.^\vee$ becomes a rigid tensor category. In particular,  for any two objects $M$ and $N$ of $\CdRB$, we have a natural isomorphism:
\begin{equation}\label{dualhom}
\begin{array}{rcl}
  \Hom_{\dRB}(M,N) & \lrasim  & \Hom_{\dR}(\Z(0), M^\vee\otimes N)  \\
  (\phi_{\dR}, \phi_{\rB}) & \lmt   & (1\mapsto \phi_{\dR}, 1 \mapsto \phi_{\rB}).  
\end{array}
\end{equation}

Moreover, for every integer $k$, we get an identification
\begin{equation}\label{homtwist}
\Hom_{\dRB}(\Z(0), M\otimes \Z(k)) \lrasim M_{\dR} \cap (2\pi i)^k M_{\rB},
\end{equation}
where the intersection is taken in $M_{\C}$, by mapping a morphism $\phi:\Z(0) \lra M\otimes\Z(k)$ to $\phi_{\C}(1)$. 

\subsubsection{Examples, I:  The (co)homology of smooth projective varieties over $\Qb$.} 

For any smooth projective variety $X$ over $\Qb$ and for any integer $i\geq 0,$ the algebraic de Rham cohomology of $X$ and the Betti cohomology of  $X^\an_{\C}$ determine an object $H^i_{\dRB}(X)$ in $\CdRB$ defined as follows:
$$H^i_{\dRB}:= (H^i_{\dR}(X/\Qb), H^i_{\rB}(X^\an_{\C},\Z)/\text{torsion}, c),$$
where $c$ denotes the composition of the comparison isomorphism defined by the base change isomorphism, analytification, and the de Rham isomorphism
$$H^i_{\dR}(X/\Qb)\otimes_{\Qb}{\C} \lrasim H^i_{\dR}(X_{\C}/\C)  \lrasim H^i_{\dR}(X^\an_{\C}) \lrasim H^i (X_{\C}^\an, \C)$$
and of the inverse of the isomorphism defined by extension of coefficients
$$ (H^i(X^\an_{\C},\Z)/\text{torsion}) \otimes_{\Z} \C \simeq  H^i(X^\an_{\C},\Z)\otimes_{\Z}\C \lrasim H^i(X^\an_{\C},\C).$$ 

To a morphism 
$$\phi : X \longrightarrow Y$$
of smooth projective varieties over $\Qb$ is attached a morphism in ``de Rham--Betti cohomology'':
$$H^i_{\dRB}(\phi):=(H^i_{\dR}(\phi),H^i_{B}(\phi))$$
defined by the ``pullback'' morphisms
$$H^i_{\dR}(\phi) := \phi^\ast : \HdR^i(Y/\Qb) \lra \HdR^i(X/\Qb)$$
and 
$$H^i_{\rB}(\phi) := \phi_{\C}^{\an \ast} : H^i(Y^\an_{\C},\Z)/\text{torsion} \lra H^i(X^\an_{\C},\Z)/\text{torsion}$$
in algebraic de Rham and Betti cohomology. This construction is clearly functorial.

Observe that, as an instance of (\ref{homtwist}), we have a natural identification:
\begin{equation}\label{HGrHom}
H^2_{\rm Gr}(X) \simeq \Hom_{\dRB}(\Z(0), H^2_{\dRB}(X)\otimes \Z(1)).
\end{equation}

We shall also define the de Rham--Betti \emph{homology} functor by duality in $\CdRB$:
$$H_{i \dRB}(X) := H^i_{\dRB}(X)^\vee \mbox{ and } H_{i \dRB}(\phi) := H^i_{\dRB}(\phi)^\vee.$$
Observe that $H_{i \dRB}(X)_{\rB}$ and $H_{i \dRB}(X)_{\C}$ may be identified with the Betti homology groups $H_{i}(X_{\C}^\an, \Z)$ modulo torsion and $H_{i}(X_{\C}^\an, \C)$ of $X_{\C}^\an$.

\subsubsection{Examples,  II: The homology of abelian varieties.}\label{ExamHom}
Let $A$ be an abelian variety of dimension $g$ over $\Qb$, and $E(A)$ its universal vector extension.

Consider the exponential map of the associated complex Lie group:
$$\exp_{E(A)_{\C}} : \Lie E(A)_{\C} \lra E(A)_{\C}^\an.$$
Its kernel, the  group of periods $\Per E(A)_{\C}$ of $E(A)_{\C}$, 
is a free $\Z$-module of rank $2g$, and the inclusion  
$\Per E(A)_{\C} \hookrightarrow \Lie E(A)_{\C}$ extends to an isomorphism 
\begin{equation}\label{perlie}
\Per E(A)_{\C} \otimes_{\Z} \C \lrasim \Lie E(A)_{\C}.
\end{equation}
Consequently we may attach the following object of $\CdRB$ to the abelian variety $A$:
$$\LiePer E(A) := (\Lie E(A), \Per E(A)_{\C}, c),$$
where $c$ denotes the inverse of the isomorphism (\ref{perlie}).

As recalled in \ref{biduve} above, the construction of $E(A)$ as the moduli space of line bundles with (integrable) connections over the dual abelian variety $\hat{A}$ provides a canonical isomorphism of $\Qb$-vector spaces:
$$ I_{A} : \Lie E(A) \lrasim H^1_{\dR}(\hA/\Qb).$$
Moreover the isomorphism of complex vector spaces
$$\begin{CD} \Lie E(A)_{\C} @>{I_{A,\C}= I_{A_{\C}}}>> \HdR^1(\hA/\Qb)\otimes_{\Qb}\C \simeq \HdR^1(\hA_{\C}/\C) @>{\text{GAGA + de Rham}}>> H^1(\hA^\an_{\C}, \C)
\end{CD}$$
 maps $\Per E(A)_{\C}$ onto $H^1(\hA_{\C}^\an, 2\pi i \Z).$ This follows from the description of $E(A)^\an_{\C}$ as $H^1(\hA^\an_{\C}, \Omega^\times_{\hA^\an_{\C}}),$ where  $\Omega^\times_{\hA^\an_{\C}}$ denotes the complex $\cO^\an_{\hA^\an_{\C}} \stackrel{d\log}{\lra} \Omega^1_{\hA^\an_{\C}}\stackrel{d}{\lra} \Omega^1_{\hA^\an_{\C}}\stackrel{d}{\lra} \cdots $.
 
 In other words, $I_{A}$ defines an isomorphism in $\CdRB$:
 $$I_{A,\dRB}: \LiePer E(A) \lrasim H^1_{\dRB}(\hA) \otimes \Z(1).$$

Besides, the isomorphism $\varpi_{A,\dR}$ constructed in paragraph \ref{subsubPoinc} above admits an obvious analogue $\varpi_{A_{\C},\rB}$ involving the Betti (co)homology of $A^\an_{\C}$ and $\hA^\an_{\C},$ which are defined by means of $c_{1\rB}(\cP_{A_{\C}}).$  Up to a factor $2\pi i$ coming from the relation 
$$c_{1\dR}(\cP_{A})_{\C} = 2 \pi i \; c_{1\rB}(\cP_{A_{\C}}),$$
it is compatible with the isomorphism   $\varpi_{A,\dR}$ in algebraic de Rham (co)homology. In other words, they define an isomorphism in $\CdRB$:
$$\varpi_{A,\dRB} := (\varpi_{A,\dR}, \varpi_{A_{\C},\rB}) : H_{1,\dRB}(A) \lrasim H^1_{\dRB}(\hA)\otimes \Z(1).$$

Finally, we get a canonical isomorphism in $\CdRB$:
\begin{equation}\label{JA}
J_{A,\dRB}:= \varpi_{A,\dRB}^{-1} \circ I_{A,\dRB} : \LiePer E(A) \lrasim H_{1,\dRB}(A).
\end{equation}
This construction is easily seen to be functorial in $A$. Namely, for any morphism $\phi: A \lra B$ of abelian varieties over $\Qb,$
$$\LiePer E(\phi) := (\Lie E(\phi), \Lie E(\phi)_{\C \vert \Per E(A)_{\C}})$$
is an element of $\Hom_{\dRB}(\LiePer E(A), \LiePer E(B)),$ and the following diagram commutes in $\CdRB$:
$$\begin{CD}
\LiePer E(A) @>{\LiePer E(\phi)}>> \LiePer E(B) \\
@V{\simeq}V{J_{A,\dRB}}V                             @V{\simeq}V{J_{B,\dRB}}V \\
H_{1,\dRB}(A) @>{H_{1\dRB}(\phi)}>> H_{1\dRB}(B).
\end{CD}
$$

\subsubsection{Extensions} For any two objects $M$ and $N$ in $\CdRB,$ we may consider the set $\Ext^1_{\dRB}(M,N)$ of 1-extensions of $M$ by $N$ in $\CdRB,$ namely of diagrams in $\CdRB$ of the form
$$\cE : \;\; 0 \lra N \stackrel{\alpha}{\lra} X \stackrel{\beta}{\lra} M \lra 0$$
such that $\beta \circ \alpha = 0$ and  the diagrams
$$\cE_{\dR} : \;\; 0 \lra N_{\dR} \stackrel{\alpha_{\dR}}{\lra} X_{\dR} \stackrel{\beta_{\dR}}{\lra} M_{\dR} \lra 0$$
and
$$\cE_{\rB} : \;\; 0 \lra N_{\rB} \stackrel{\alpha_{\rB}}{\lra} X_{\rB} \stackrel{\beta_{\rB}}{\lra} M_{\rB} \lra 0$$
are short exact sequences of $\Qb$-vector spaces and of $\Z$-modules respectively.

Equipped with the Baer sum, $\Ext^1_{\dRB}(M,N)$ becomes an abelian group. Actually, for any extension $\cE$ as above, we may choose a $\Qb$-linear splitting $\sigma_{\dR}:  M_{\dR} \rightarrow X_{\dR}$ of $\cE_{\dR}$ and a $\Z$-linear splitting $\sigma_{\rB} : M_{\rB} \rightarrow X_{\rB}$ of $\cE_{\rB}.$ Then $\sigma_{\dR\C}:= \sigma_{\dR} \otimes_{\Qb}1_{\C}$ and $\sigma_{\rB\C}:= \sigma_{\rB}\otimes_{\Z}1_{\C}$ are $\C$-linear splittings of 
$$\cE_{\C} : \;\; 0 \lra N_{\C} \stackrel{\alpha_{\C}}{\lra} X_{\C} \stackrel{\beta_{\C}}{\lra} M_{\C} \lra 0,$$
and consequently $\sigma_{\dR\C}-\sigma_{\rB\C}$ may be written $\alpha_{\C}\circ \phi$ for some uniquely determined $\phi$ in $(M^\vee \otimes N)_{\C}.$ The map
\begin{equation}\label{Extiso}
\begin{array}{rcl}
 \Ext^1_{\dRB}(M,N) & \lrasim    & (M^\vee \otimes N)_{\C}/[(M^\vee \otimes N)_{\dR}+(M^\vee \otimes N)_{\rB}]  \\
 {[\cE]}& \lmt   & [\phi] 
\end{array}
\end{equation}
so defined is easily seen to be an isomorphism of abelian groups.

In particular,we get the usual isomorphisms:
\begin{equation}\label{Extisobis}
\Ext^1_{\dRB}(M,N) \lrasim \Ext^1_{\dRB}(\Z(0), M^\vee\otimes N) \lrasim \Ext^1_{\dRB}(M\otimes N^\vee, \Z(0)).
\end{equation}

\subsection{Abelian varieties over $\Qb$ satisfy $GPC^1$}

We are now in position to complete the proof of Theorem \ref{GPCA}.

As already observed, universal vector extensions of abelian varieties satisfy Condition $\mathbf{LP}$. Corollary \ref{CorSL2} therefore implies that, for any two abelian varieties $A$ and $B$ aver $\Qb$, the map 
$$\begin{array}{rcl}
\LiePer :\Hom_{\gp/\Qb} (E(A), E(B))  & \lra   & \Hom_{\dRB} (\LiePer E(A), \LiePer E(B))   \\
 \psi  & \lmt  &   \LiePer \,\psi := (\Lie \psi, \Lie \psi_{\C \vert \Per E(A)_{\C}}).
\end{array}
$$
is an isomorphism of $\Z$-modules.

Together with the isomorphism (\ref{Eff}), which identifies morphisms between abelian varieties and between their universal vector extensions, this establishes the first assertion in the following theorem; the second assertion follows from the existence of a functorial isomorphism (\ref{JA}) between $\LiePer E(A)$ and $H_{1,\dRB}(A)$:

\begin{theorem}\label{HomAB}
For any two abelian varieties $A$ and $B$ over $\Qb$, the maps 
$$
\begin{array}{rcl}
  \Hom_{\gp/\Qb}(A,B) & \lra   & \Hom_{\dRB} (\LiePer E(A), \LiePer E(B)) \\
\phi  & \lmt   & \LiePer E(\phi)   
\end{array}
$$
and 
$$H_{1,\dRB}:   \Hom_{\gp/\Qb}(A,B)  \lra \Hom_{\dRB}(H_{1,\dRB}(A),H_{1,\dRB}(B))$$
are isomorphisms of $\Z$-modules.
\end{theorem}

In other words, the realization functor $H_{1,\dRB}$ from the category of abelian varieties over $\Qb$ to the category $\CdRB$ is fully faithful. (Compare with  \cite{AndreMotives04}, 7.5.3, where a ``rational'' version of this isomorphism is established, by a reference to  some  advanced transcendence results of W\"ustholz \cite{Wuestholz84}.)  

To complete the proof of Theorem \ref{GPCA}, we consider an abelian variety $A$ over $\Qb$ and we apply Theorem \ref{HomAB} to $A$ and its dual abelian variety $\hat{A}$. In this way, we get an isomorphism
$$H_{1,\dRB}:   \Hom_{\gp/\Qb}(A,\hat{A})  \lrasim \Hom_{\dRB}(H_{1,\dRB}(A),H_{1,\dRB}(\hA)).$$
Composing this isomorphism with the transpose of 
$$\varpi_{A,\dRB}  : H_{1,\dRB}(A) \lrasim H^1_{\dRB}(\hA)\otimes \Z(1),$$
and with the natural identification (\ref{dualhom}), we get an isomorphism
\begin{equation}
\Hom_{\gp/\Qb}(A,\hat{A}) \lrasim \Hom_{dRB}(\Z(0), H^1_{\dRB}(A)\otimes H^1_{\dRB}(A)\otimes \Z(1)).
\end{equation}
The discussion on signs in paragraph \ref{NSsym} (notably the identity (\ref{altphi})) shows that this isomorphism maps the subgroup of \emph{symmetric} morphisms from $A$ to $\hat{A}$ onto the subgroup of skew-symmetric, or \emph{alternating}, elements\footnote{Namely the  elements sent to their opposite by the automorphism of $\Hom_{dRB}(\Z(0), H^1_{\dRB}(A)\otimes H^1_{\dRB}(A)\otimes \Z(1))$ defined by ``switching'' the two copies of $H^1_{\dRB}(A)$.} in $\Hom_{\dRB}(\Z(0), H^1_{\dRB}(A)\otimes H^1_{\dRB}(A)\otimes \Z(1))$:
\begin{equation}\label{isosymalt}
\Hom_{\gp/\Qb}(A,\hat{A})^{\text{sym}} \lrasim \Hom_{dRB}(\Z(0), H^1_{\dRB}(A)\otimes H^1_{\dRB}(A)\otimes \Z(1))^{\text{alt}}.
\end{equation} 

The fact that the morphism of $\Z$-modules in (\ref{isosymalt}) is an isomorphism is nothing but, in a disguised form, the validity of $GPC^1(A)$.  
Indeed, by composition 
with the isomorphism 
$$
\begin{array}{rcl}
 NS(A) := \Pic(A)/\Pic_{0}(A) & \lrasim  & \Hom_{\gp/\Qb}(A,\hA)^{\text{sym}}  \\
{[L]} & \lmt   & \phi_{L},  
\end{array}
$$ 
the isomorphism (\ref{isosymalt})  becomes the  isomorphism
\begin{equation}\label{isosymaltbis}
NS(A) \lrasim \Hom_{\dRB}(\Z(0), H^1_{\dRB}(A)\otimes H^1_{\dRB}(A)\otimes \Z(1))^{\text{alt}}.
\end{equation}
The ``Betti'' component of (\ref{isosymaltbis}) takes its values in $(H^1_{B}(A_{\C})\otimes_{\Z} H^1_{B}(A_{\C}))^{\text{alt}}$ and is well known  to coincide with the classical ``Riemann form'' of elements of the N\'eron-Severi group (see for instance \cite{BirkenhakeLange04}, Chapter 2).  Consequently, after the identification of 
$$\Hom_{\dRB}(\Z(0), H^1_{\dRB}(A)\otimes H^1_{\dRB}(A)\otimes \Z(1))^{\text{alt}}$$ and    $$\Hom_{\dRB}(\Z(0), H^2_{\dRB}(A) \otimes \Z(1)) = H^2_{\text{Gr}}(A),$$ the isomorphism (\ref{isosymalt}) may be read as asserting  that the map
$$c_{1\dRB}: NS(A) \lra H^2_{\text{Gr}}(A)$$
is an isomorphism. This is precisely the content of $GPC^1(A)$. 

\subsection{$\Qb$-points of abelian varieties and extensions in $\CdRB$}\footnote{This section could be skipped at first reading. It has been included since Proposition \ref{KummerTransc} constitutes an application of the theorem of Schneider-Lang close in spirit to the ones in the previous section, and for comparison with Conjecture \ref{Main} \emph{infra}.}\label{kappadRB}

Let $A$ denote an abelian variety over $\Qb$.

Consider some line bundle $L$ over $A$, algebraically equivalent to zero,  equipped with some rigidification $\epsilon : k \simeq L_{0_{A}}.$  
Recall that the $\G_{m}$-torsor $L^\times \xrightarrow{\pi_{L}} A$ over $A$, deduced from the total space of $L$ by deleting its zero section, may be endowed with a unique structure of $\Qb$-algebraic group which makes the diagram
$$0 \lra \G_{m\Qb} \stackrel{\epsilon}{\lra} L^\times \stackrel{\pi_{L}}{\lra} A \lra 0$$
a short exact sequence of commutative $\Qb$-algebraic groups, and that this construction establishes an isomorphism of groups: 
$$\hA(\Qb) \lrasim \Ext^1_{\cgp/\Qb}(A, \G_{m\Qb}).$$

The fiber product
$$E(L^\times) \simeq L^\times \times_{A} E(A)$$
defines a commutative $\Qb$-algebraic group which fits into the following commutative diagram with exact lines:
$$\begin{CD}
0 @>>> \G_{m\Qb} @>{\tilde{\epsilon}}>> E(L^\times) @>{\tilde{\pi}_{L}}>> E(A) @>>> 0 \\
@.           @VV{=}V                                             @VVV                                        @VV{p_{A}}V    @. \\
0 @>>> \G_{m\Qb} @>{\epsilon}>> L^\times  @>{\pi_{L}}>> A @>>> 0.
\end{CD}
$$

By considering the Lie algebra (over $\Qb$) and the periods (over $\C$) of the first line, we get a 1-extension in $\CdRB$:
\begin{equation}\label{ExtdRBL}
0 \lra \Z(1) \xrightarrow{\LiePer\, \tilde{\epsilon}} \LiePer E(L^\times)  \xrightarrow{\LiePer\, \tilde{\pi}_{L}} \LiePer E(A) 
 \lra 0.
\end{equation}
 Thanks to the canonical isomorphisms in $\CdRB$
 $$\LiePer E(A) \lrasim H_{1\dRB}(A) \lrasim H^1_{\dRB}(\hA)\otimes \Z(1) \lrasim H_{1\dRB}(A)^\vee \otimes \Z(1),$$
its class may be seen as an element $\kappa_{\dRB}(L)$ in
 $$\Ext_{\dRB}^1(H_{1\dRB}(A), \Z(1)) \lrasim  \Ext_{\dRB}^1(\Z(0), H_{1\dRB}(\hA))$$
 and defines a morphism of abelian groups
 $$\kappa_{\dRB} : \hA(\Qb) \lra \Ext_{\dRB}^1(\Z(0), H_{1\dRB}(\hA)).$$

The proof of the following proposition is again an application of Corollary \ref{CorSL2}:

\begin{proposition}\label{KummerTransc} The map 
$\kappa_{\dRB}$ 
is injective.
\end{proposition}

We leave the details to the reader, and only emphasize that giving a direct description of the subgroup $\kappa_{\dRB}( \hA(\Qb))$ of  $\Ext_{\dRB}^1(\Z(0), H_{1\dRB}(\hA))$ appears to be an intriguing and difficult issue.

\section{$D$-group schemes}

In this part, we introduce $D$-schemes and $D$-group schemes in a  geometric setting, suitable for the application to Diophantine geometry we want to discuss in the sequel. These definitions are variants of the original definitions by Buium   (\cite{Buium86}, \cite{Buium92}, \cite{Buium94} Chapter 3), which make sense over some fixed differential base field (of characteristic zero). Here we shall consider $D$-schemes and group schemes over some smooth base variety instead : this framework is  the one of Malgrange in \cite{Malgrange10}, with the field of complex numbers  replaced by some arbitrary  field of characteristic zero.  

For simplicity, we shall make  smoothness and quasi-projectivity assumptions which actually could be relaxed in many places. Actually, on a base scheme of finite type over a field of characteristic zero, $D$-schemes are nothing but the ``crystals in relative schemes'' mentioned in  a famous letter of Grothendieck to Tate\footnote{quoted in \cite{Illusie94}, Section 4.1:\emph{``un cristal possède deux propriétés caractéristiques : la rigidité et la faculté de croître dans un voisinage approprié. Il y a des cristaux de toute espèce de substances : des cristaux de soude, de soufre, de modules, d'anneaux,  de sch\'emas relatifs, etc."}}. The approach to $D$-schemes as ``crystals'', defined in terms of infinitesimal sites and stratifications, has much to recommend it (see for instance \cite{Simpson94II}, Section 8), but I have preferred to stick to a more naive approach in the spirit of classical differential geometry,  at the expense of extra regularity assumptions, based on a definition of  $D$-schemes that mimics the one of integrable Ehresmann connections on differentiable fiber bundles (\cite{Ehresmann51}).  

In the following sections we denote $k$ a fixed field  \emph{of characteristic zero}. 

\subsection{Basic definitions}\label{basicD} 

Let $S$ denote a smooth quasi-projective scheme over $k$. 

\subsubsection{$D$-schemes}\label{Ds}

By a \emph{$D$-scheme} over $S$, we shall mean a pair $(X,\cF)$ where $X\stackrel{\pi}{\longrightarrow} S$ is a smooth, quasi-projective scheme over $S$ (hence over $k$), and $\cF$ is an integrable\footnote{In other words, its sheaf of regular sections is closed under Lie bracket.} sub-vector bundle of the ``absolute'' tangent bundle $T_{X/k}$ of $X$ such that
\begin{equation*}
T_{X/k} = T_{X/S} \oplus \cF.
\end{equation*}
This last condition means precisely that $\cF$ determines a splitting of the exact sequence of vector bundles over the $k$-scheme $X$
$$0 \lra T_{X/S} \hlra T_{X/k} \stackrel{D\pi}{\lra} \pi^\ast T_{S/k} \lra 0$$
defined by the differential of $\pi,$ or equivalently that the restriction of $D\pi$ to $\cF$ is an isomorphism:
\begin{equation}
\label{Dpiiso}
D\pi_{\vert \cF} : \cF \lrasim \pi^\ast T_{S/k}.
\end{equation}

A \emph{morphism of $D$-schemes} over $S$
\begin{equation}\label{morD}
\phi : (X_{1},\cF_{1}) \lra (X_{2},\cF_{2})
\end{equation}
is a morphism of $S$-schemes $\phi : X_{1} \rightarrow X_{2}$ whose ``absolute''  differential
$$D\phi : T_{X_{1}/k} \lra \phi^\ast T_{X_{2}/k}$$
maps $\cF_{1}$ to $\phi^\ast\cF_{2}$.

Observe that, if $\phi$ is a morphism of $D$-schemes over $S$ from $(X_{1},\cF_{1})$ to $(X_{2},\cF_{2})$, then Conditions (\ref{Dpiiso}) for $(X_{1},\cF_{1})$ and $(X_{2},\cF_{2})$ imply that  $D\phi $ maps $\cF_{1}$ isomorphically onto $\phi^\ast\cF_{2}$.

Morphisms of $D$-schemes may be obviously composed and define the category of (smooth, quasi-projective) $D$-schemes over $S$. Clearly, this category admits finite products: $(S, T_{S/k})$ is a final object, and the product of two $D$-schemes 
$(X_{1},\cF_{1})$ and $(X_{2},\cF_{2})$ over $S$ may be constructed  as the $D$-scheme $(X,\cF)$ consisting of their product as schemes over $S$,
$$X:=X_{1} \times_{S} X_{2},$$
equipped with the sub-vector bundle $\cF$ of $T_{X/k}$ which is the ``direct sum of $\cF_{1}$ and $\cF_{2}$ over $T_{S/k}$,'' formally defined as the kernel of the surjective morphism of vector bundles over $X$:
$$(D_{\pi_{1}}, -D_{\pi_{2}}) : (\cF_{1}\boxplus \cF_{2})_{\vert X} \lra \pi^\ast T_{S/k}.$$
(It lies inside the kernel of 
$$(D_{\pi_{1}}, -D_{\pi_{2}}) : (T_{X_{1}/k}\boxplus T_{X_{2}/k})_{\vert X} \lra \pi^\ast T_{S/k},$$
which may be identified with $T_{X/k}$.) 

A \emph{closed $D$-subscheme} of a $D$-scheme $(X,\cF)$ over $S$ is the image of a morphism of $D$-schemes with range $(X,\cF)$ that is also a closed immersion. Equivalently it is a closed, smooth subscheme $Y$ of $X$ such that its tangent bundle $T_{Y/k},$ which may be identified to a sub-vector bundle of $T_{X/k\vert Y},$ contains $\cF_{\vert Y}.$

A \emph{horizontal section}  of some $D$-scheme $(X,\cF)$ over $S$ is a right inverse of the structural morphism $X \lra S$ in the category of $D$-schemes over $S$. In other words, it is a section $\cP$ of this morphism over $S$, the differential of which $D\cP:T_{S/k}\lra \cP^\ast T_{X/k}$ takes its values in $\cP^\ast \cF,$ or equivalently, the image of which is a $D$-subscheme of $(X,\cF).$  

From the integrable sub-vector bundle $\cF$ of $T_{X/k}$, the normal bundle $\cP^\ast T_{X/S}$ of any horizontal section $\cP$ inherits an integrable connection.

\subsubsection{$D$-group schemes}\label{Dgs} A (smooth, quasi-projective) \emph{$D$-group scheme} over $S$ is defined as a group object in the category of $D$-schemes over $S$. 

A $D$-group scheme $\bG$ over $S$ may be identified with a pair $(G,\cF)$ where $G$ is a smooth, quasi-projective group scheme over $S$ and $\cF$ a sub-vector bundle of $T_{G/k}$ which makes $(G,\cF)$ a $D$-scheme over $S$, in such a way that the graphs of the unit section $e_{G}$, of the inverse map, and of the composition map of the group scheme $G$ become $D$-subschemes of the $D$-schemes $G$, $G^2$, and $G^3$ over $S$.  

In intuitive terms, a $D$-group scheme may be   thought as a smooth group scheme over $S$ equipped with some ``algebraic connection'' compatible with its group structure.

Since its unit section $e_{G}$ is horizontal, the relative Lie algebra $\Lieb_{S}G := e_{G}^\ast T_{G/S}$ of the group scheme $G$ over $S$ underlying some $D$-group scheme $\bG$ over $S$ becomes endowed with a natural integrable connection. The so-defined module with integrable connection shall be denoted $\Lieb_{S}\bG$.  


Assume that $S$ is integral (or equivalently, connected), of dimension $s$, and consider its field of rational functions $k(S).$ Let us choose some $k(S)$-basis $(v_{1},\ldots,v_{s})$ of the $k(S)$-vector space of rational sections of $T_{S/k}$ such that  the Lie brackets $[v_{i},v_{j}]$ all vanish\footnote{Such bases exist: simply write $k(S)$ as a finite degree extension of $k(X_{1},\ldots,X_{s})$, and lift the standard basis $(\partial/\partial X_{1},\ldots,\partial/\partial X_{s})$.}. Then the field $k(S)$ equipped with the derivations $(\delta_{1},\ldots,\delta_{s})$ becomes a differential field in the classical  sense of Ritt and Kolchin. Let us finally choose a differential closure $(K; \delta_{1},\ldots,\delta_{s})$ of $(k(S); \delta_{1},\ldots,\delta_{s}).$
Through the base changes 
$$\Spec K \lra \Spec k(S) \hlra S,$$
any $D$-group scheme $(G,\cF)$ over $S$ in our sense defines  $D$-group schemes in the sense of Buium over the differential fields 
$(k(S); \delta_{1},\ldots,\delta_{s})$ and $(K; \delta_{1},\ldots,\delta_{s})$, and a  $\Delta_{0}$-group, that is, a differential algebraic group of finite dimension in the sense of Kolchin,  by considering the subgroup of the group $G(K)$ of $K$-points of $G$ consisting of its ``horizontal points". (We refer the reader to \cite{Buium92}, Chapter 5, \cite{Pillay97},  \cite{Pillay04}, and \cite{BertrandPillay10}  for discussions of the relations between Buium's $D$-groups and differential algebraic groups.)

\subsubsection{Extensions}\label{subsubsec:Ext}

Let $\bG_{1} =(G_{1}, \cF_{1})$ and $\bG_{2}=(G_{2},\cF_{2})$ be two commutative $D$-group schemes over $S$. An \emph{extension} of $\bG_{1}$ by $\bG_{2}$ in the category of commutative $D$-group schemes over $S$ is a diagram
\begin{equation}\label{extDtyp}
0 \lra \bG_{2} \stackrel{i}{\lra} \bG \stackrel{p}{\lra} \bG_{1} \lra 0
\end{equation}
in this category such that the underlying diagram of commutative group schemes over $S$ 
$$0 \lra G_{2} \stackrel{i}{\lra} G \stackrel{p}{\lra} G_{1} \lra 0$$
is a short exact sequence\footnote{As usual, by this we mean a short exact sequence of fppf sheaves over $S$. Since we work over  a base field $k$ of characteristic zero, this is equivalent to the following ``geometric'' condition, expressed in terms of some algebraic closure $\oli{k}$ of $k$: for any point $s \in S(\oli{k}),$ the diagram
$$0 \lra G_{2s}(\oli{k}) \stackrel{i_{s}}{\lra} G_{s}(\oli{k}) \stackrel{p_{s}}{\lra} G_{1s}(\oli{k}) \lra 0$$
is a short exact sequence of abelian groups.} (compare \cite{KowalskiPillay06}).

The Baer sum of two extensions of $\bG_{1}$ by $\bG_{2}$ may be defined in an obvious way. Equipped with this operation, the set $\Ext^1_{\cDgp/S}(\bG_{1},\bG_{2})$ of isomorphism classes of these extensions defines an abelian group, which satisfies the usual functorialities in $S$, $\bG_{1}$, and $\bG_{2}$.

We may apply the functor $\Lieb_{S}$ to the extension (\ref{extDtyp}). We obtain a short exact sequence of modules with integrable connections over $S$:
$$0 \lra \Lieb_{S}\bG_{2} \xrightarrow{\Lieb_{S} i}
\Lieb_{S}\bG \xrightarrow{\Lieb_{S}p}
\Lieb_{S}\bG_{1} \lra 0.$$
This construction defines an additive map, say when $S$ is projective:
$$\Lieb_{S}^1 : \Ext^1_{\cDgp/S}(\bG_{1},\bG_{2}) \lra \Ext^1_{\mic/S}(\Lieb_S{\bG_{1}}, \Lieb_{S}\bG_{2}) \simeq
H^1_{\dR}(S, (\Lieb_S{\bG_{1}})^\vee\otimes\Lieb_{S}\bG_{2}),$$
where we use the notation introduced in paragraph \ref{deRhamcoeff}, formula (\ref{coeffalg}).

\subsubsection{Functoriality in $S$}
If $\phi: S' \lra S$ is a morphism of projective schemes over $k$, then, from any $D$-scheme $(X,\cF)$ over $S$, we may deduce a $D$-scheme $(X',\cF')$ over $S'$ by ``pulling it back'' by $\phi$ as follows : $X'$ is the smooth, quasi-projective $S'$-scheme defined as the fiber product $X\times_{S}S'$;  if $\tilde{\phi} : X' \lra X$ denotes the canonical ``first projection'' morphism and $D\tilde{\phi}: T_{X'/k} \lra \tilde{\phi}^\ast T_{X/k}$ its differential,  the $D$-structure on $X'$ over  $S'$ is defined by the integrable sub-vector bundle of $T_{X'/k}$
$$\cF' :=  D\tilde{\phi}^{-1} (\tilde{\phi}^\ast \cF).$$

This construction of ``base change'' is functorial, and transforms $D$-group schemes over $S$ into $D$-group schemes over $S'$. It satisfies an obvious compatibility with the Lie algebra functor (from $D$-group schemes to modules with integrable connections) and the pullback of modules with integrable connections.  

The $D$-schemes over $\Spec k$ are nothing but the smooth, quasi-projective schemes over $k$. A \emph{constant} $D$-scheme over $S$ is a $D$-scheme isomorphic to the pullback by the $k$-morphism $S\lra \Spec k$ of some  smooth, quasi-projective schemes over $k$. In the sequel, we shall denote $\bG_{m,S}$ the constant multiplicative group scheme over $S$, defined as the pullback of the algebraic group $\G_{m,k}$. After the change of base $S\lra \Spec k$, the isomorphism
$$
\begin{array}{rcl}
 \Lie \G_{m,k} & \lrasim   & k   \\
X.\partial /\partial X  & \longmapsto   & 1
\end{array}  
$$
becomes an isomorphism of modules with integrable connections:
$$\Lieb_{S} \bG_{m,S} \lrasim (\cO_{S}, d).$$

\subsubsection{Change of base fields}\label{kk'} If $k'$ is a field extension of $k$, the extension of scalars from $k$ to $k'$ associates a $D$-scheme $(X_{k'}, \cF_{k'})$ over $S_{k'}$, defined over the base field $k'$, to any $D$-scheme $(X, \cF)$ over $S$. This operation satisfies obvious functoriality properties that we shall use freely in the sequel. In particular, it attaches $D$-group schemes over $S_{k'}$ to $D$-group schemes over $S$, and defines morphisms of extension groups:
$$\Ext^1_{\cDgp/S}(\bG_{1},\bG_{2}) \lra \Ext^1_{\cDgp/S_{k'}}(\bG_{1k'},\bG_{2k'}).$$

\subsection{$D$-schemes and analytification}\label{Dan} 
When the base field $k$ is $\C,$ a $D$-scheme $(X,\cF)$ (resp. a $D$-group scheme $(G,\cF)$) over $S$ determines, through analytification, a ``$D$-analytic space''  $(X^\an, \cF^\an)$ (resp. a ``$D$-complex Lie group'' $(G^\an, \cF^\an)$) over the complex manifold $S^\an$. We shall omit the formal definitions of these notions --- just ``copy''  the above ones in the analytic context --- and content ourselves with a few observations.

First, after analytification, a $D$-scheme $(X,\cF)$ \emph{projective} over $S$ becomes locally constant in the analytic category. Namely, for any point $s_{0}$ of $S^\an$, there exists an open neighbourhood $\Omega$ of $s_{0}$ in $S^\an$ and an isomorphism of $\C$-analytic spaces over $\Omega$
\begin{equation}\label{psiun}
\Psi_{s_{0}}: \Omega \times X^\an_{s_{0}} \lrasim X^\an_{\Omega}
\end{equation}
such that
\begin{equation}\label{init}
\Psi_{s_{0}}(s_{0}, . ) = Id_{X^\an_{s_{0}}}
\end{equation}
and, for any $(s,x)$ in $\Omega \times X^\an_{s_{0}},$
\begin{equation}\label{horiz}
\cF_{\Psi_{s_{0}}(s,x)} = D\Psi_{s_{0}}(s,x)(T_{s}\Omega \oplus 0).
\end{equation}
This follows from the analytic integrability of $\cF^\an$, together with the properness of the structural morphism $X^\an \lra S^\an$ in the analytic topology. (Observe that Conditions (\ref{init}) and  (\ref{horiz}) uniquely determine $\Psi_{s_{0}}$ for $\Omega$ connected.)      
  
  Second, as pointed out by Hamm (\cf \cite{Buium92}, Chapter 2, 1.3), a similar statement holds for any $D$-group scheme $(G, \cF)$ over $S$. Thus we get a (unique) isomorphism of complex Lie groups over\footnote{By a ``complex Lie group  over a complex analytic manifold $M$'', we mean a group object in the category of complex analytic manifolds ``smooth'' (in the ``algebrogeometric'' sense, that is ``submersive'') over $M$.} $\Omega$ (assumed to be small enough and connected)
  \begin{equation}\label{psideux}
\Psi_{s_{0}}: \Omega \times G^\an_{s_{0}} \lrasim G^\an_{\Omega}
\end{equation}
which  satisfy the initial condition (\ref{init}) and the horizontality condition (\ref{horiz}).

Consider in particular the case of a commutative $D$-group scheme $\bG = (G,\cF)$ over $S$, with connected fibers. Then the ``relative'' exponential map
$$\exp_{G/S} : \Lieb_{S} G \lra G^\an$$
defines a surjective morphism of complex Lie groups over $S^\an$. It is compatible with the ``horizontal'' structures defined by the integrable connection  on $\Lieb_{S} \bG$ and by (\ref{horiz}), and consequently its kernel 
$$\Perb_{S} G := \ker \exp_{G/S}$$
is a local system (that is, a locally free sheaf) of $\Z$-modules of finite rank over $S^\an$, which fits into a short exact sequence in the category of commutative complex Lie groups over $S^\an$:
$$0 \lra \Perb_{S} G \hlra \Lieb_{S} G \xrightarrow{\exp_{G/S}}
 G^\an \lra 0.$$
 
This is even a short exact sequence of commutative $D$-complex Lie groups, which should be denoted 
$$0 \lra \Perb_{S} G \hlra \Lieb_{S} \bG \xrightarrow{\exp_{G/S}}
 \bG^\an \lra 0.$$
 This shows, in particular, that when $s$ varies in $S^\an$, the dimension of the complex sub-vector space of $\Lie G_{s}$ generated by its period lattice $\Per G_{s}$ is locally constant (in the analytic topology). Consequently, if $S$ (hence $S^\an$) is connected and if, for some $s_{0} \in S$, $G_{s_{0}}$ satisfies condition $\mathbf{LP}$ (\cf Section \ref{Morag}),  then $G_{s}$ satisfies $\mathbf{LP}$ for every $s$ in $S$, and the structure of $\bG$ as a $D$-group scheme over $S$ is uniquely determined by its structure of a group scheme.  Similarly, if $\bG_{1}$ and $\bG_{2}$ are two commutative $D$-groups schemes  over $S$, and if $\bG_{1}$ has connected  fibers satisfying $\mathbf{LP},$ then any morphism of group schemes from $G_{1}$ to $G_{2}$ is a morphism of $D$-group schemes from $\bG_{1}$ to $\bG_{2}$.
 
 These remarks will apply to the $D$-group schemes associated to abelian schemes and to  their extension by multiplicative groups considered in Sections \ref{ED} and \ref{ExtD} \emph{infra}. (See also \cite{BertrandPillay10}, Lemma 3.4, for similar unicity statements in a more ``differential algebraic" formulation.) 
 
 Associating its   local system of periods  $\Perb_{S}G$ to a $D$-group scheme $\bG$ is a functorial construction (in $S$ and $\bG$). Applied to extensions, it defines a morphism of $\Z$-modules, for any two commutative $D$-groups schemes $\bG_{1}$ and $\bG_{2}$ with connected fibers over $S$:
\begin{equation}\label{Perext}
\Perb^1_{S}: \Ext^1_{\cDgp/S}(\bG_{1},\bG_{2}) \lra \Ext^1_{\text{Ab-Sheaves}/S^\an}(\Perb_{S}G_{1},\Perb_{S}G_{2})
\simeq H^1(S^\an,(\Perb_{S}G_{1})^\vee \otimes
\Perb_{S}G_{2}).
\end{equation}

\subsection{Moduli spaces of vector bundles with connections as $D$-schemes}\label{ModuliD} If the $S$-scheme $X$ underlying some $D$-scheme $(X,\cF)$ as above is projective over $S,$ then, locally in the \'etale topology of $S$, $X$ is ``constant"  over $S$ (namely, when $k$ is algebraically closed, of the form $X_{0}\times_{k} S$, after replacing $S$ by some \'etale neighborhood of any given point of $S$).  This follows from the representability of the Isom-functors in the projective case, together with the formal integrability of $\cF$ and Artin's algebraization theorem (compare with \cite{Buium86},  II.1, and \cite{Gillet02}, Section 3).

This property is a refinement, which makes sense in pure algebraic geometry, of the local analytic triviality of projective $D$-schemes when $k= \C$. It strongly limits the possible constructions of smooth projective $D$-schemes.

It is remarkable that, in contrast, highly ``nonconstant'' smooth \emph{quasi-projective} $D$-schemes arise naturally. Indeed the construction of the moduli spaces  $\MIC_{N}(M,o)$ of vector bundles with connection recalled in  paragraph \ref{MIC} above, applied to smooth projective families of pointed projective varieties parameterized by $S$, provides quasi-projective $D$-schemes over $S$. 

Namely, if $M$ is a smooth, projective $S$-scheme with geometrically connected fibers, and if $o$ denotes a section of $M$ over $S$, then Simpson's techniques apply to this relative situation. They lead to  the construction of  a flat, quasi-projective $S$-scheme\footnote{In \cite{Simpson94II}, this $S$-scheme is denoted ${\mathbf R}_{\rm DR}(M/S,o,N)$.} $\MIC_{N}(M/S, o)$, the fiber of which 
 over some point $s\in S(\oli{k})$ may be identified with the moduli space $\MIC_{N}(M_{s}, o(s))$. Formally, for any $S$-scheme $\Sigma,$ $\MIC_{N}(M/S, o)(\Sigma)$ classifies vector bundles of rank $N$ over $X_{\Sigma}:=X\times_{S}\Sigma,$ rigidified over $o_{\Sigma},$ and equipped with an integrable connection relative to $\Sigma.$

The $S$-scheme $\MIC_{N}(M/S, o)$ admits a canonical structure of $D$-scheme over $S$, which reflects its so-called crystalline nature. For general $M$ and $N$, this scheme may actually not be smooth over $S$, and properly speaking  it is not covered by the above definition of $D$-schemes (which should be replaced by a suitable definition in terms of the infinitesimal site and stratifications associated to $X/k$). However, in the sequel, we shall be mainly concerned by the situation where $N=1,$ in which case $\MIC_{1}(M/S, o)$ is a \emph{smooth}, quasi-projective, group scheme over $S$, and we allow ourselves to neglect this issue of regularity.

When $k=\C,$ the $D$-scheme structure of  $\MIC_{N}(M/S, o)$ may be described as follows. When $s$ varies in the complex manifold $S^\an,$ the family of fundamental groups
$$\Gamma_{s} := \pi_{1}(M^\an_{s}, o(s))$$
define a local system (that is, a locally constant sheaf) of groups on $S^\an.$ Over any simply connected open subset $\Omega$ in $S^\an$, it may be trivialized :  for any pair of points $(s_{0},s_{1})$ in $\Omega,$ we get a canonical isomorphism
$$\gamma_{s_{1},s_{0}}: \Gamma_{s_{0}} \lrasim \Gamma_{s_{1}},$$
which clearly induces an isomorphism of representations spaces:
\begin{equation*}
\begin{array}{crcl}
\Phi^{\text{Rep}}_{s_{1},s_{0}}:& \Rep_{N}(\Gamma_{s_{0}}) & \lrasim  & \Rep_{N}(\Gamma_{s_{1}})   \\
 & \rho& \lmt   & \rho \circ \gamma_{s_{0},s_{1}}^{-1}\;\; . 
\end{array} 
\end{equation*}
Moreover the monodromy isomorphisms (\ref{aniso}) 
$$\mon_{o(s)}: \MIC_{N}(M/S,o)_{s} =  \MIC_{N}(M_{s},o(s)) \lrasim \Rep_{N}(\Gamma_{s})$$
and their inverses
depend analytically  on $s$, in the sense that, if $s_{0}$ denotes a base point in $\Omega,$
the bijection of sets
\begin{equation}\label{monOm}
\begin{array}{crcl}
\Psi_{s_{0}}:& \Omega \times \Rep_{N}(\Gamma_{s_{0}}) & \lrasim  &  \MIC_{N}(M/S,o)_{\Omega}  \\
 & (s, \rho) & \lmt   & \mon^{-1}_{s}(\Phi^{\text{Rep}}_{s, s_{0}}(\rho))
\end{array} 
\end{equation}
is an isomorphism of $\C$-analytic spaces over $\Omega.$

The $D$-scheme structure over $s$ of $X:=\MIC_{N}(M/S,o)$ is compatible with the ``analytic trivialization'' (\ref{monOm}). Assume indeed that $\MIC_{N}(M/S,o)$ is smooth over $S$ (for instance, suppose that $N=1$); then the subvector bundle $\cF$ of $T_{X/\C}$ which defines this structure becomes ``horizontal'' via the above isomorphism :
$$\mbox{for any $(s,\rho)  \in  \Omega \times \Rep_{N}(\Gamma_{s_{0}})$, } \cF_{\Phi(s,\rho)} = D\Psi_{s_{0}}(s,\rho)(T_{s}\Omega \oplus 0).$$

It is quite remarkable that the \emph{analytic} sub-vector bundle $\cF$ of $T_{X/\C}$ defined through this formula in terms of the local analytic trivializations  (\ref{monOm}) of $X$ over $S$ is an \emph{algebraic} subvector bundle of $T_{X/\C}.$

This is due to Grothendieck and to Mazur and Messing (\cite{MazurMessing74}) when $N=1$ (see also \cite{Buium92}), and to Simpson (\cite{Simpson94II}, Section 8) in general. 
Basically their proof consists in considering the avatar in formal geometry (over the formal completion $\widehat{S}_{s_{0}}$ of $S$ at $s_{0}$) of the local analytic trivialization of $\MIC_{N}(M/S,o)$ over $\Omega$ induced by (\ref{monOm}):
\begin{equation}\label{monOmbis}
 \Psi^{\text{MIC}}_{s_{0}}:= \Psi_{s_{0}}\circ (Id_{\Omega}\times \mon_{o(s_{0})}) : \Omega \times  \MIC_{N}(M/S,o)_{0}  \lrasim  \MIC_{N}(M/S,o)_{\Omega}. 
\end{equation}
It turns out that the formal analogue of (\ref{monOmbis}) over $\widehat{S}_{s_{0}}$ may be directly constructed in (formal) algebraic geometry, with no recourse to analytic techniques, over any base field $k$ of characteristic zero. 

 The existence of the local analytic trivializations $\Psi^{\text{MIC}}_{s_{0}}$ is indeed a direct consequence of the following basic observation : if $(E,\nabla)$ is an analytic vector bundle with integrable connection over some connected analytic submanifold $Y$ of some analytic manifold $X$, then $(E,\nabla)$ uniquely extends, as a vector bundle with integrable connection, over any sufficiently small open  connected neighbourhood of $Y$ in $X$. This property admits a natural avatar in formal geometry, valid over any base field of characteristic zero, which implies the existence of a formal analogue of $\Psi^{\text{MIC}}_{s_{0}}$. This construction, with $s_{0}$ varying in $S,$ endows $\MIC_{N}(M/S,o)$ with a structure of $D$-scheme over $S.$  
 
\subsection{Universal vector extensions as $D$-group schemes}\label{ED}
\footnote{The content of Sections \ref{ED} and \ref{ExtD} is thoroughly discussed, with a slightly different perspective, in \cite{BertrandPillay10}, Part 3 and Appendix, which constitutes the main reference  for these two sections.} The above discussion may be specialized to the case $N=1$. Then $\MIC_{1}(M/S,o)$ is a smooth, quasi-projective group scheme over $S$ --- its group structure is induced by the tensor product of rigidified line bundles with connections --- and its neutral component $\MIC_{1}(M/S,o)^0$ may be identified with the universal vector extension $E(\Pic_{0}(M/S))$ of the connected relative Picard variety $\Pic_{0}(M/S)$ of $M$ over $S$. Moreover, the structure of a $D$-scheme over $S$ on $\MIC_{1}(M/S,o)$ is compatible with its structure of a group scheme.

Let us introduce the relative Albanese variety of $M$ over $S$, namely the abelian scheme over $S$ defined as
$$\cA := \widehat{\Pic_{0}(M/S)},$$
and the relative Albanese morphism 
$$\alpha_{o}: M \lra \cA$$
attached to the section $o$. It induces an isomorphism of group schemes over $S$ (see for instance \cite{BK09}, Appendix B):
$$\alpha_{o}^\ast : \MIC_{1}(M/S, o)^0 \lrasim \MIC_{1}(\cA/S,0_{\cA, 0})^0,$$
compatible with their structure of $D$-schemes. Together with the identification of group schemes over $S$
$$\MIC_{1}(\cA/S,0_{\cA})^0 = \MIC_{1}(\cA/S,0_{\cA}) \lrasim E(\widehat{\cA}),$$
this shows that (i) to study $\MIC_{1}(M/S, o)^0$, we may consider the case where $M$ is some abelian scheme over $S;$ and (ii) that the universal vector extension $E(\widehat{\cA})$ --- hence by duality the universal vector extension of any abelian scheme over $S$ --- is endowed with a natural structure of $D$-group schemes, that we shall denote $\bfE(\widehat{\cA})$.  

The analytic description of the $D$-structure on the moduli spaces $\MIC_{N}(M/S, o)$ boils down in the present situation to the following description of the $D$-group scheme $\bfE(\cB)$  defined by the universal vector extension $E(\cB)$ attached to some abelian scheme $\cB$ (see also \cite{MazurMessing74}, 4.4).

Assume that $k=\C,$ and consider an abelian scheme  over $S$, of relative dimension $g$, 
$$\pi: \cB \lra S.$$
 As in Section \ref{Dan}, we may consider 
 the analytic description of the complex Lie group $\cB^\an$ over $S^\an$ as a quotient of $\Lieb_{S}\cB$ by its local system of periods :
$$0 \lra \Perb_{S} \cB \hlra \Lieb_{S} \cB \xrightarrow{\exp_{\cB/S}}
 \cB^\an \lra 0.$$

 This local system $\Perb_{S}\cB$ is locally free of rank $2g,$ and may be identified with the local systems of fundamental groups, of fiber at $s\in S$:
 $$\pi_{1}(\cB_{s}, 0_{\cB_{s}}) \simeq H_{1}(\cB^\an_{s}, \Z).$$
 In the sequel, we shall denote it $\cH_{1\rB}(\cB^\an/S^\an)$. In turn, the dual local system
 $$\cH^1_{\rB}(\cB^\an/S^\an) := \cH_{1\rB}(\cB^\an/S^\an)^\vee$$
 may be identified with $R^1\pi^\an_{\ast} \Z_{\cB^\an}.$ 
 
 As discussed in paragraph \ref{ExamHom}, for any $s\in S^\an,$ we have a canonical isomorphism
 \begin{equation}
J_{\cB_{s}}: \Lie E(\cB_{s}) \lrasim H_{1\dR}(\cB_{s}/\C) \simeq H_{1}(\cB^\an_{s}, \C) \simeq H_{1}(\cB^\an_{s}, \Z) \otimes_{\Z}\C,
\end{equation}
which sends $\Per \cB_{s}$ isomorphically onto $H_{1}(\cB^\an_{s}, \Z).$ These isomorphisms depend analytically on $s\in S^\an,$ and define  isomorphisms $J_{\cB}$ of analytic vector bundles and local systems over $S^\an$, which fit into a commutative diagram:
$$\begin{CD}
\Perb_{S} E(\cB)   @>{J_{\cB}}>\sim> \cH_{1\rB}(\cB^\an/S^\an) \\
@VVV                             @VVV \\
\Lieb_{S}E(\cB) @>{J_{\cB}}>\sim> \cH_{1\rB}(\cB^\an/S^\an)\otimes_{\Z}\C,
\end{CD}
$$
where the vertical maps are the obvious injections. They induce an isomorphism of complex Lie groups over $S^\an$:
$$J_{\cB}^\times : E(\cB)^\an \lrasim \cH_{1\rB}(\cB^\an/S^\an) \otimes_{\Z}\G^\an_{m\C}$$
which makes the following diagram commutative:
\begin{equation}\label{extaniso}
\begin{CD}
0 @>>> \cH_{1\rB}(\cB^\an/S^\an) @>{J_{\cB}^{-1}}>> \Lieb_{S}E(\cB) @>{\exp_{E(\cB)/S}}>>                       E(\cB)^\an @>>>0 \\
@.                   @VV{=}V                                                         @VV{J_{\cB}}V                                                                    @VV{J_{\cB}^\times}V     @. \\
0 @>>> \cH_{1\rB}(\cB^\an/S^\an) @>{. \otimes_{\Z}1_{\C}}>>\cH_{1\rB}(\cB^\an/S^\an)\otimes_{\Z} \C @>{Id_{\cH_{1\rB}}\otimes_{\Z}\bfe}>> \cH_{1\rB}(\cB^\an/S^\an) \otimes_{\Z}\G^\an_{m\C} @>>>0.
\end{CD}
\end{equation}
(Recall that $\bfe := \exp (2\pi i .)$.)

 In (\ref{extaniso}),  both lines are short exact sequences of commutative complex Lie groups over $S^\an$, and
 the vertical arrows are isomorphisms. These isomorphisms are actually compatible with the $D$-structures in the analytic category: the connection on $\Lieb_{S}\bfE(\cB)$ is the dual of the Gauss-Manin connection on $\cH^1_{\dR}(\cB/S)$, and is mapped by $J_{\cB}$ to the connection on $\cH_{1\rB}(\cB^\an/S^\an)\otimes_{\Z}\C$ which makes horizontal the sections of the local system $\cH_{1\rB}(\cB^\an/S^\an)$; the local analytic trivializations of $E(\cB)^\an$ induced by the $D$-structure of $\bfE(\cB)$, become, under the isomorphism $J_{\cB}^\times$, the local  trivializations of $\cH_{1\rB}(\cB^\an/S^\an) \otimes_{\Z}\G^\an_{m\C}$ induced by local trivializations of $\cH_{1\rB}(\cB^\an/S^\an)$.
 
\subsection{Extensions of abelian schemes by $\G_{m}$ and $D$-group schemes}\label{ExtD}

The construction of the algebraic groups $L^\times$ and $E(L^\times)$ attached to some line bundle $L$ algebraically equivalent to zero on some abelian variety $A$ discussed  in Section \ref{kappadRB} extends to a relative situation. 

Consider for instance an abelian scheme $\cB$ over $S$ as in the previous section. If $\cL$ is a line bundle over $\cB,$ equipped with a rigidification along the zero section of $\cB$
$$\epsilon : \cO_{S} \lrasim 0_{\cB}^\ast \cL,$$
and algebraically equivalent to zero on the fibers of $\cB$ --- in other words, if $(\cL, \epsilon)$ defines a section $\cP$ over $S$ over the dual abelian scheme $\widehat{\cB}$ ---  then the $\G_{m}$-torsor $\pi_{\cL} : \cL^\times \lra \cB$, deduced from the total space of $\cL$ by deleting its zero section, admits a unique structure of a commutative group scheme over $S$ which makes the diagram
\begin{equation}\label{extcL}
0 \lra \G_{m S} \stackrel{\epsilon}{\lra} \cL^\times \stackrel{\pi_{\cL}}{\lra} \cB \lra 0
\end{equation}
an extension of smooth commutative group schemes over $S$. By pulling back this extension along the morphism 
$$p_{\cB}: E(\cB) \lra \cB,$$
 we define a smooth commutative group scheme 
 $$E(\cL^\times) := \cL^\times \times_{\cB} E(\cB)$$
 which fits into an short exact sequence of group schemes over $S$:
  \begin{equation}\label{extEcL}
0 \lra \G_{m S} \stackrel{\epsilon'}{\lra} E(\cL^\times) \stackrel{\tilde{\pi}_{\cL}}{\lra} E(\cB) \lra 0
\end{equation}

In the sequel we shall use that $E(\cL^\times)$ may be canonically equipped with a $D$-structure, so that it becomes a commutative $D$-group scheme $\bfE(\cL^\times)$ over $S$ and the extension of commutative group schemes (\ref{extEcL}) becomes an extension of commutative $D$-group schemes:   
  \begin{equation}\label{DextEcL}
0 \lra \bG_{m S} \stackrel{\epsilon'}{\lra} \bfE(\cL^\times) \stackrel{\tilde{\pi}_{\cL}}{\lra} \bfE(\cB) \lra 0
\end{equation}

This construction is alluded to in \cite{Brylinski83} (2.2.2.1), appears  in a ``differential algebraic context''  in \cite{BertrandPillay10} Lemma 3.4 (i-ii),  and  in a ``geometric context'' in \cite{AndreattaBarbieri-Viale05} (see also \cite{AndreattaBertapelle11}). The construction of the $D$-structure on $\bfE(\cL^\times)$ and of the extension (\ref{DextEcL}) may be understood as follows in terms of moduli spaces of vector bundles with integrable connections. 

The construction of the relative moduli spaces $\MIC_{N}(M/S, o)$ and of their $D$-structure discussed   in Section \ref{ModuliD} directly extends to the moduli spaces $\MIC_{N}(M/S, o,o')$ of vector bundles equipped with a relative integrable connection rigidified along two sections $o$ and $o'$ of $M$ over $S$. Besides, as explained in Section \ref{ED}, $\bfE(\cB)$ may be identified with the $D$-group scheme $\MIC_{1}(\widehat{\cB}/S, 0_{\widehat{\cB}})$. The discussion of Section \ref{VarMic} may be extended to the relative case, and allows one to identify $E(\cL^\times)$ with $\MIC_{1}(\widehat{\cB}/S, 0_{\widehat{\cB}}, \cP)$, in a way compatible with their respective structure of $\G_{m}$-torsors over $E(\cB)$ and $\MIC_{1}(\widehat{\cB}/S, 0_{\widehat{\cB}})$. The canonical $D$-structure on  $E(\cL^\times)$ is the $D$-structure   deduced from the one  on $\MIC_{1}(\widehat{\cB}/S, 0_{\widehat{\cB}}, \cP)$ through this identification.


\section{A conjecture}\label{sec:Conj}

 In this final part, we consider the following geometric data: a smooth projective connected curve $C$ over $\Qb$, and an abelian scheme over $C$,
$\pi: \cA \lra C.$

As before we denote $E(\cA)$ the universal vector extension of this abelian scheme. It is a smooth connected commutative group scheme over $C,$ endowed with a canonical structure of a $D$-group scheme. If necessary, we shall use the notation $\bfE(\cA)$ to denote $E(\cA)$ considered as a $D$-group scheme over $C$, to distinguish it from the ``plain'' group scheme $E(\cA)$ over $C$. 

As usual, we  denote $\widehat{\cA}$ the  abelian scheme over $C$ dual to $\cA$. 
 
 We shall make the following simplifying assumption :
 \begin{equation}\label{Hodgepositive}
\mbox{\emph{the vector bundle $\E_{\cA} := (\Lieb_{C} \cA)^\vee$ is ample.}}
\end{equation}
Recall that, in general, $\E_{\cA}$ is only semipositive.
 Condition (\ref{Hodgepositive}) implies the vanishing of the $\overline{\Qb(C)}/\Qb$-trace of the geometric generic fiber $\cA_{\overline{\Qb(C)}}$ of $\cA$, and shall 
 ensure
 that the extensions of formal $D$-groups (\ref{forC}) and local systems (\ref{extGamma}) considered below have no nontrivial automorphisms (hence have their middle term defined, up to unique isomorphism, by their extension class). 
 
 \subsection{A construction}\label{Cons}
 
 Suppose that we are given the following datum:
 
 \noindent (i) \emph{a section $\cP$ over $C$ of the dual abelian scheme} $\cA.$
 
 By the very definition of $\cA$, it defines 
 
 \noindent (ii) \emph{a line bundle $\cL$ over $\cA$, equipped with a rigidification $\epsilon : \cO_{C} \lrasim 0_{\cA}^\ast \cL$ along the zero section, algebraically equivalent to zero on the fiber of $\pi: \cA \lra C.$}
 
 As recalled above, the $\G_{m}$-torsor $\cL^\times$ over $\cA$ defines in a unique way 
 
 \noindent (iii) \emph{an extension of smooth commutative group schemes over $C$,}
 $$0 \lra \G_{m, S} \stackrel{\epsilon}{\lra} \cL^\times \lra \cA \lra 0.$$
 
 Finally, through the construction descibed in Section \ref{ExtD}, we obtain:
 
 \noindent (iv) \emph{an extension of commutative $D$-group scheme over $C$,}
 $$0 \lra  \bfGm_{S}   \lra \bfE(\cL^\times)\lra \bfE(\cA) \lra 0.$$
 
These successive operations are easily seen to establish a bijective correspondence between the four kinds of data (i)--(iv) above, and to be additive:

 \begin{lemma} The above construction defines isomorphisms of $\Z$-modules:
 \emph{
 $${\widehat{\cA}(C) \lrasim \Ext^1_{{\cgp}/C}(\cA, \G_{m S}) \lrasim}\Ext_{\cDgp}^1(\bfE(\cA),\bfGm_{S}).$$}
\end{lemma}

This would actually hold in the general situation considered in Section \ref{ExtD}, without any further assumption on the base scheme $S$.

\subsection{$\Lieb_{C}^1$ and $\Perb^1_{C_{\C}^\an}$}\label{subsec:LiePer1}

Recall that the dual of the module with integrable connection $\Lieb_{C} \bfE(\cA)$ over $C$ may be identified with the relative de Rham cohomology of $\cA$ over $C$ equipped with the Gauss-Manin connection
$(\cH^1_{\dR}(\cA/C), \nabla_{GM})$, and the local system of periods $\Perb_{C_{\C}}E(\cA)_{\C}$ over $C^\an_{\C}$ with the local system defined by the relative Betti first homology of $\cA^\an_{\C}$ over $\C_{\C}^\an$, which we denote $\cH^1_{\rB}(\cA_{\C}^\an/C^\an_{\C}).$

Besides, the module with integrable connection $\Lieb_{C} \bG_{m,C}$ over $C$ may be identified with the trivial module with integrable connection $(\cO_{C}, d)$, and the local system of periods $\Perb_{C_{\C}} \G_{m,C_{\C}}$ over $C^\an_{\C}$ with the constant local system $\Z_{C^\an_{\C}}.$

Consequently the maps $\Lieb^1_{S}$ and $\Perb^1_{S}$ defined on extension classes of commutative $D$-group schemes in paragraph \ref{subsubsec:Ext} and Section \ref{Dan} take here the following form: 
 $$\Lieb_{C}^1 : \Ext^1_{\cDgp/C}(\bfE(\cA),\bfGm_{S}) \lra H^1_{\dR}(C, (\cH^1_{\dR}(\cA/C), \nabla_{GM}))$$
and
$$\Perb^1_{C_{\C}^\an} : \Ext^1_{\cDgp/C_{\C}}(\bfE(\cA)_{\C},\bfGm_{S_{\C}}) \lra H^1(C^\an_{\C}, \cH^1_{\rB}(\cA_{\C}^\an/C^\an_{\C})).$$

Observe that, after tensoring with $\C$, the range spaces of these two maps become canonically isomorphic. Indeed we have ``elementary'' isomorphisms defined by the base change from $\Qb$ to $\C$
\begin{equation}\label{elem1}
H^1_{\dR}(C, (\cH^1_{\dR}(\cA/C), \nabla_{GM})) \otimes_{\Qb} \C \lrasim H^1_{\dR}(C_{\C}, (\cH^1_{\dR}(\cA_{\C}/C_{\C}), \nabla_{GM}))
\end{equation}
and by extension of coefficients from $\Z$ to $\C$
\begin{equation}\label{elem2}
H^1(C^\an_{\C}, \cH^1_{\rB}(\cA_{\C}^\an/C^\an_{\C}))\otimes_{\Z}\C \lrasim
H^1(C^\an_{\C}, \cH^1_{\rB}(\cA_{\C}^\an/C^\an_{\C})_{\C}),
\end{equation}
and the complex vector spaces in the right-hand sides of (\ref{elem1}) and (\ref{elem2}) may be identified by means of the comparison isomorphisms between Betti and algebraic de Rham cohomology (with coefficients) discussed in paragraph \ref{deRhamcoeff}.

If $\cE$ is an element of $\Ext^1_{\cDgp}(\bfE(\cA),\bfGm_{S})$, we shall denote  $\cE_{\C}$ its ``complexification'' in the group $ \Ext^1_{\cDgp}(\bfE(\cA)_{\C},\bfGm_{S_{\C}})$ (in the sense of paragraph \ref{kk'}).

The following lemma is  proved in the same way as Lemma \ref{compChern}, which compared the first Chern classes in de Rham and Betti cohomology (see also the discussion in paragraph \ref{final} \emph{infra}).

\begin{lemma}\label{compdRB1}
For any extension class $\cE$ in \emph{$\Ext^1_{\cDgp}(\bfE(\cA),\bfGm_{C})$},  the  equality 
\begin{equation}\label{eqcompdRB1}
(\Lieb_{C}^1 \cE)\otimes_{\Qb} 1_{\C} = 2 \pi i (\Perb^1_{C_{\C}^\an} \cE_{\C}) \otimes_{\Z} 1_{\C}
\end{equation}
holds in 
$$H^1_{\dR}(C, (\cH^1_{\dR}(\cA/C), \nabla_{GM})) \otimes_{\Qb}\C \simeq H^1(C^\an_{\C}, \cH^1_{\rB}(\cA_{\C}^\an/C^\an_{\C}))\otimes_{\Z}\C.$$

\end{lemma}

\subsection{A conjecture}\label{subsec:Conj}

\subsubsection{} We  finally arrive at the formulation of the conjecture which constitutes the aim of this article.

\begin{conjecture}\label{Main} Any pair of classes of extensions  $(\alpha, \beta)$  with $\alpha$ in 
$H^1_{\dR}(C, (\cH^1_{\dR}(\cA/C), \nabla_{GM}))$ and $\beta$ in $H^1(C^\an_{\C}, \cH^1_{\rB}(\cA_{\C}^\an/C^\an_{\C}))$ which satisfies the compatibility relation
\begin{equation}\label{eqcompdRB1bis}
\alpha\otimes_{\Qb} 1_{\C} = 2 \pi i \, \beta \otimes_{\Z} 1_{\C}
\end{equation}
in 
$$H^1_{\dR}(C, (\cH^1_{\dR}(\cA/C), \nabla_{GM})) \otimes_{\Qb}\C \simeq H^1(C^\an_{\C}, \cH^1_{\rB}(\cA_{\C}^\an/C^\an_{\C}))\otimes_{\Z}\C$$
is of the form $(\Lieb_{S}\cE, \Perb^1_{C_{\C}^\an} \cE_{\C})$ for some  class $\cE$ in \emph{$\Ext^1_{\cDgp}(\bfE(\cA),\bfGm_{C})$} and hence is obtained from some  section $\cP$ of the dual abelian scheme $\widehat{\cA}$ over $C$.
\end{conjecture}

The class $\cE$ and the section $\cP$, if they exist, are uniquely determined by these conditions.

By using the Leray-Serre spectral sequence to analyze the group $H^2_{\rm Gr}(\cA)$ attached to $\cA$ (seen as a smooth projective variety over $\Qb$) by means of the fibering 
$\pi : \cA \lra C,$ and by using a relative generalization (over $C$) of Theorem \ref{GPCA}, we may prove the following:

\begin{proposition}\label{MainGPCA} With the above notation, Conjecture \ref{Main} holds if and only if  the smooth projective variety $\cA$ over $\Qb$ satisfies $GPC^1(\cA).$ 
\end{proposition}

\subsubsection{} Consider $f : S \lra C$ a smooth projective connected surface $S$ over $\Qb$ fibered over $C$. Assume for simplicity that $f$ is a smooth morphism (all fibers of $f$ are therefore smooth projective curve) and admits a section $o$. Then we may introduce the relative Jacobian 
$$\cJ := \Jac (S/C)$$
of $S$ over $C$. It is an abelian scheme over $C$. Using the section $o$, we may define a relative Jacobian embedding
$$j_{o}: S \hlra \cJ.$$
(It is a closed embedding, over $S$, which maps $o$ to the zero section $0_{\cJ}$ of $\cJ$ over $C$.) Pulling back by $j_{o}$ establishes a bijection between line bundles $\cL$ over $\cJ$, defining as above sections over $C$ of the dual abelian schemes $\widehat{\cJ}$\footnote{that is, line bundles rigidified along $\cJ$, and algebraically equivalent to zero in the fibers of $\cJ$ over $C$.}, and line bundles $\cM$ over $S$, rigidified along $o$ and of degree zero on the fibers of $f$. 

With this notation, we have the following variant\footnote{This variant is actually simpler than Proposition \ref{MainGPCA}: its proof does not require Theorem \ref{GPCA} and its relative generalization.}  of Proposition \ref{MainGPCA}:
\begin{proposition}\label{MainGPCS} The validity of $GPC^1(S)$ is equivalent to the validity of Conjecture \ref{Main} for $\cA = \cJ$. 
\end{proposition} 

Conjecture \ref{Main} may be extended to possibly degenerating families of abelian varieties over $C$ (say, with semi-abelian bad fibers). This generalized version may be applied to the relative Jacobian of any smooth projective surface fibered over $C$ (say, with semi-stable fibers) and would imply the validity of $GPC^1$ for any smooth projective surface and, actually, for any smooth projective variety over $\Qb$. This approach to $GPC^1$ through fibrations of surfaces over curves and associated families of Jacobian varieties is very much in the spirit of the classical works of Picard, Poincar\'e, and Lefschetz which constituted our starting point in Section \ref{subsec:algline}.

\subsubsection{}\label{final} To avoid technicalities, I prefer not to discuss this in detail, and would instead stress the fact that Conjecture \ref{Main}  may be rephrased as an algebraization criterion concerning formal line bundles, satisfying suitable ``differential algebraic'' and ``analytic'' conditions, in the spirit of Theorems \ref{SL1} and \ref{SL2} \emph{\`a la} Schneider-Lang, as expected in Section  \ref{subsec:analogy}.

Indeed, consider a pair of classes $(\alpha, \beta)$ as in  Conjecture $\ref{Main}$. 

The class $\alpha$ lies in 
$$H^1_{\dR}(C, (\cH^1_{\dR}(\cA/C), \nabla_{GM})) \simeq \Ext^1_{\mic/C}(\Lieb_{C}\bfE(\cA),\Lieb_{C}\bfGm_{S}),$$
and defines an extension of vector bundles with (integrable) connections over $C$, defined over $\Qb$:
$$ 0 \lra \Lieb_{C}\bfGm_{C} \lra (M,\nabla) \lra \Lieb_{C}\bfE(\cA) \lra 0.$$   
It may be interpreted as an extension of ``formal commutative $D$-group schemes over $C$'':
\begin{equation}\label{forC}
0 \lra \widehat{\bfGm}_{C} \lra \mathbf{G}_{\text{for}} \lra \widehat{\bfE(\cA)} \lra 0,
\end{equation}
where $\widehat{\bfGm}_{C}$ (resp. $\widehat{\bfE(\cA)}$) denotes the completion of the $D$-group scheme
$\bfGm_{C}$ (resp. $\bfE(\cA)$) over $C$ along its unit (resp. zero) section. (Here we use that the base field $\Qb$ has characteristic zero, so that we have formal exponential maps at our disposal.)

Observe that, by forgetting the $D$-structure, from (\ref{forC}) we deduce an extension of formal groups over $C$,
$$0 \lra \widehat{\Gm}_{C} \lra {G}_{\text{for}} \lra \widehat{E(\cA)} \lra 0,$$
which in turn defines a $\G_{m}$-torsor or, equivalently, a line bundle $\cN_{\text{for}}$, on the formal completion  $\widehat{E(\cA)}$.

The class $\beta$ lies in 
$$H^1(C^\an_{\C}, \cH^1_{\rB}(\cA_{\C}^\an/C^\an_{\C})) 
\simeq
\Ext^1_{\text{Ab-Sheaves}}(\Perb_{C_{\C}}\cA_{\C}, \Z_{C^\an_{\C}})$$
and defines an extension of local systems over free $\Z$-modules of finite rank over $C^\an_{\C}$:
\begin{equation}\label{extGamma}
0 \lra \Z_{C^\an_{\C}} \lra \Gamma \lra \Perb_{C_{\C}}\cA_{\C} \lra 0.
\end{equation}
After tensoring with the multiplicative group $\G_{m\C}^\an$, we deduce from (\ref{extGamma}) an extension of ``commutative $D$-complex Lie groups'' over $C_{\C}^\an$:
\begin{equation}\label{extGammaLie}
0 \lra {\bG}^\an_{m,C_{\C}} \lra \Gamma \otimes \bG^\an_{m,C_{\C}} \lra \bfE(\cA)^\an_{\C} \lra 0.
\end{equation}
This construction is easily seen to establish a one-to-one correspondence between extensions of local systems (\ref{extGamma}) and extension in the analytic category of  $\bfE(\cA)^\an_{\C}$ by ${\bG}^\an_{m,C_{\C}}$. When $\beta$ is the image by $\Perb^1_{C_{\C}}$ of some extension class $[\cE_{\C}]$, the extension (\ref{extGammaLie}) is nothing but the analytification $\cE^\an_{\C}$ of $\cE_{\C}$.

Here again the extension (\ref{extGammaLie}) defines some analytic line bundle $\cN^\an$ over $E(\cA)_{\C}^\an$, by forgetting the $D$-structure and part of the group structure on
$\Gamma \otimes \G^\an_{m,C_{\C}}$.

The equality (\ref{eqcompdRB1bis})
$$\alpha\otimes_{\Qb} 1_{\C} = 2 \pi i \, \beta \otimes_{\Z} 1_{\C}$$
expresses the fact that the extension of ``commutative formal analytic $D$-groups'' over $C^\an_{\C}$ deduced from  (\ref{extGammaLie}) by completion along the zero sections coincides with the analytification of  the ``commutative formal $D$-groups'' over $C_{\C}$ deduced from (\ref{forC}) by extending the base field from $\Qb$ to $\C.$

Finally Conjecture \ref{Main} may be rephrased as asserting the algebraicity of any pair $(\cN^{\text{for}}, \cN^\an)$, consisting of a formal line bundle $\cN^{\text{for}}$ on the formal completion $\widehat{E(\cA)}$ of $E(\cA)$ along its zero section and of some analytic line bundle $\cN^\an$ over $E(\cA)^\an_{\C}$ such that the associated $\G_{m}$-torsors $\cN^{\text{for} \times}$ and
 $\cN^{\an \times}$ may be endowed with suitably compatible structures of $D$-group schemes over $C$ and $C^\an_{\C}$ (in the respective formal and analytic categories).

\newcommand{\etalchar}[1]{$^{#1}$}



\begin{thebibliography}{KGG{\etalchar{+}}04}

\bibitem[And04]{AndreMotives04}
Y.~Andr\'e.
\newblock {\em { Une introduction aux motifs. Motifs purs, motifs mixtes,
  p\'eriodes}}.
\newblock {Panoramas et Synth\`eses 17. Soci\'et\'e Math\'ematique de France},
  2004.

\bibitem[ABV05]{AndreattaBarbieri-Viale05}
F.~Andreatta and L.~Barbieri-Viale.
\newblock Crystalline realizations of 1-motives.
\newblock {\em Math. Ann.}, 331(1):111--172, 2005.

\bibitem[AB11]{AndreattaBertapelle11}
F.~Andreatta and A.~Bertapelle.
\newblock Universal extension crystals of 1-motives and applications.
\newblock {\em J. Pure Appl. Algebra}, 215(8):1919--1944, 2011.



\bibitem[And63]{Andreotti63}
A.~Andreotti.
\newblock Th\'eor\`emes de d\'ependance alg\'ebrique sur les espaces complexes
  pseudo-concaves.
\newblock {\em Bull. Soc. Math. France}, 91:1--38, 1963.



\bibitem[Ano56]{Anonymous56}
Anonymous.
\newblock Correspondence.
\newblock {\em Amer. J. Math.}, 78:898, 1956.

\bibitem[BW07]{BakerWuestholz07}
A.~Baker and G.~W{\"u}stholz.
\newblock {\em Logarithmic forms and {D}iophantine geometry}, volume~9 of {\em
  New Mathematical Monographs}.
\newblock Cambridge University Press, Cambridge, 2007.

\bibitem[BBM82]{BerthelotBreenMessing82}
P.~Berthelot, L.~Breen, and W.~Messing.
\newblock {\em Th\'eorie de {D}ieudonn\'e cristalline. {II}}, volume 930 of
  {\em Lecture Notes in Mathematics}.
\newblock Springer-Verlag, Berlin, 1982.

\bibitem[Ber83]{Bertrand83}
D.~Bertrand.
\newblock Endomorphismes de groupes alg\'ebriques; applications
  arithm\'etiques.
\newblock In {\em Diophantine approximations and transcendental numbers
  ({L}uminy, 1982)}, volume~31 of {\em Progr. Math.}, pages 1--45. Birkh\"auser
  Boston, Boston, MA, 1983.

\bibitem[BP10]{BertrandPillay10}
D.~Bertrand and A.~Pillay.
\newblock A {L}indemann-{W}eierstrass theorem for semi-abelian varieties over
  function fields.
\newblock {\em J. Amer. Math. Soc.}, 23(2):491--533, 2010.


\bibitem[BL04]{BirkenhakeLange04}
C.~Birkenhake and H.~Lange.
\newblock {\em Complex abelian varieties}, volume 302 of {\em Grundlehren der
  Mathematischen Wissenschaften}.
\newblock Springer-Verlag, Berlin, second edition, 2004.

\bibitem[BM01]{BogomolovMcQuillan01}
F.~Bogomolov and M.~L. McQuillan.
\newblock Rational curves on foliated varieties.
\newblock Preprint I.H.E.S., 2001.

\bibitem[Bom70]{Bombieri70}
E.~Bombieri.
\newblock {Algebraic values of meromorphic maps.}
\newblock {\em Invent. Math.}, 10:267--287, 1970.



\bibitem[Bos01]{Bost01}
J.-B. Bost.
\newblock Algebraic leaves of algebraic foliations over number fields.
\newblock {\em Publ. Math. I.H.E.S.}, 93:161--221, 2001.

\bibitem[Bos04]{Bost04}
J.-B. Bost.
\newblock {Germs of analytic varieties in algebraic varieties: canonical
  metrics and arithmetic algebraization theorems.}
\newblock In {\em {A. Adolphson et al. (ed.), Geometric aspects of Dwork
  theory}}, volume~II, pages 371--418. Walter de Gruyter, Berlin, 2004.

\bibitem[Bos06]{Bost06}
J.-B. Bost.
\newblock Evaluation maps, slopes, and algebraicity criteria.
\newblock In {\em Proceedings of the {I}nternational {C}ongress of
  {M}athematicians, {M}adrid 2006}, volume~II, pages 537--562. European
  Mathematical Society, 2006.

\bibitem[BCL09]{BostChambert-Loir07}
J.-B. Bost and A.~Chambert-Loir.
\newblock Analytic curves in algebraic varieties over number fields.
\newblock In {\em Algebra, arithmetic, and geometry: in honor of {Y}u. {I}.
  {M}anin. {V}ol. {I}}, volume 269 of {\em Progr. Math.}, pages 69--124.
  Birkh\"auser, Boston, MA, 2009.

\bibitem[BK09]{BK09}
J.-B. Bost and K.~K{\"u}nnemann.
\newblock Hermitian vector bundles and extension groups on arithmetic schemes
  {II}. {T}he arithmetic {A}tiyah extension.
\newblock {\em Ast\'erisque}, (327):361--424 (2010), 2009.





\bibitem[Bou98]{Bouscaren98}
E.~Bouscaren, editor.
\newblock {\em Model theory and algebraic geometry. {A}n introduction to {E}.
  {H}rushovski's proof of the geometric {M}ordell-{L}ang conjecture}, volume
  1696 of {\em Lecture Notes in Mathematics}.
\newblock Springer-Verlag, Berlin, 1998.


\bibitem[Bry83]{Brylinski83}
J.-L. Brylinski.
\newblock {$``1$}-motifs'' et formes automorphes (th\'eorie arithm\'etique des
  domaines de {S}iegel).
\newblock In {\em Conference on automorphic theory ({D}ijon, 1981)}, volume~15
  of {\em Publ. Math. Univ. Paris VII}, pages 43--106. Univ. Paris VII, Paris,
  1983.

\bibitem[Bui86]{Buium86}
A.~Buium.
\newblock {\em Differential function fields and moduli of algebraic varieties},
  volume 1226 of {\em Lecture Notes in Mathematics}.
\newblock Springer-Verlag, Berlin, 1986.

\bibitem[Bui92a]{Buium92}
A.~Buium.
\newblock {\em Differential algebraic groups of finite dimension}, volume 1506
  of {\em Lecture Notes in Mathematics}.
\newblock Springer-Verlag, Berlin, 1992.

\bibitem[Bui92b]{Buium92Annals}
A.~Buium.
\newblock Intersections in jet spaces and a conjecture of {S}. {L}ang.
\newblock {\em Ann. of Math. (2)}, 136(3):557--567, 1992.

\bibitem[Bui93a]{Buium93eff}
A.~Buium.
\newblock Effective bound for the geometric {L}ang conjecture.
\newblock {\em Duke Math. J.}, 71(2):475--499, 1993.

\bibitem[Bui93b]{Buium93}
A.~Buium.
\newblock Geometry of differential polynomial functions. {I}. {A}lgebraic
  groups.
\newblock {\em Amer. J. Math.}, 115(6):1385--1444, 1993.

\bibitem[Bui94]{Buium94}
A.~Buium.
\newblock {\em Differential algebra and {D}iophantine geometry}.
\newblock Actualit\'es Math\'ematiques. Hermann, Paris, 1994.

\bibitem[BV93]{BuiumVoloch93}
A.~Buium and J.~F. Voloch.
\newblock Integral points of abelian varieties over function fields of
  characteristic zero.
\newblock {\em Math. Ann.}, 297(2):303--307, 1993.

\bibitem[CS53]{CartanSerre53}
H.~Cartan and J.-P. Serre.
\newblock Un th\'eor\`eme de finitude concernant les vari\'et\'es analytiques
  compactes.
\newblock {\em C. R. Acad. Sci. Paris}, 237:128--130, 1953.

\bibitem[CL02]{Chambert01}
A.~Chambert-Loir.
\newblock Th\' eor\`emes d'alg\'ebricit\' e en g\'eom\'etrie diophantienne.
  {S}\'eminaire {B}ourbaki, {V}ol.\ 2000/2001, {E}xpos\'e 886.
\newblock {\em Ast\'erisque}, 282:175--209, 2002.

\bibitem[Cho49]{Chow49}
W.-L. Chow.
\newblock On compact complex analytic varieties.
\newblock {\em Amer. J. Math.}, 71:893--914, 1949.

\bibitem[Cho86]{Chow86}
W.-L. Chow.
\newblock Formal functions on homogeneous spaces.
\newblock {\em Invent. Math.}, 86(1):115--130, 1986.

\bibitem[CC85a]{ChudnovskysGroth85}
D.~V. Chudnovsky and G.~V. Chudnovsky.
\newblock Applications of {P}ad\'e approximations to the {G}rothendieck
  conjecture on linear differential equations.
\newblock In {\em Number theory, Semin. New York 1983-84}, volume 1135 of {\em
  Lectures Notes in Mathematics}, pages 52--100. Springer, Berlin, 1985.

\bibitem[CC85b]{ChudnovskysAcad85}
D.~V. Chudnovsky and G.~V. Chudnovsky.
\newblock Pad\'e approximations and {D}iophantine geometry.
\newblock {\em Proc. Nat. Acad. Sci. U.S.A.}, 82(8):2212--2216, 1985.





\bibitem[Col91]{Coleman91}
R.~F. Coleman.
\newblock The universal vectorial bi-extension and {$p$}-adic heights.
\newblock {\em Invent. Math.}, 103(3):631--650, 1991.

\bibitem[Con48]{Conforto48}
F.~Conforto.
\newblock Sopra le trasformazioni in s\`e della variet\`a di {J}acobi relativa
  ad una curva di genere effettivo diverso dal genere virtuale, in ispecie nel
  caso di genere effettivo nullo.
\newblock {\em Ann. Mat. Pura Appl. (4)}, 27:273--291, 1948.

\bibitem[Con49]{Conforto49II}
F.~Conforto.
\newblock Sulla nozione di corpi equivalenti e di corpi coincidenti nella
  teoria delle funzioni quasi abeliane.
\newblock {\em Rend. Sem. Mat. Univ. Padova}, 18:292--310, 1949.


\bibitem[Del71]{Deligne71}
P.~Deligne.
\newblock Th\'eorie de {H}odge. {II}.
\newblock {\em Inst. Hautes \'Etudes Sci. Publ. Math.}, (40):5--57, 1971.

\bibitem[Dem82]{Demailly82}
J.-P. Demailly.
\newblock Formules de {J}ensen en plusieurs variables et applications
  arithm\'etiques.
\newblock {\em Bull. Soc. Math. France}, 110(1):75--102, 1982.

\bibitem[Ehr51]{Ehresmann51}
C.~Ehresmann.
\newblock Les connexions infinit\'esimales dans un espace fibr\'e
  diff\'erentiable.
\newblock In {\em Centre Belge Rech. Math., Colloque de Topologie, Bruxelles,
  Juin 1950}, pages 29--55. Masson, Paris, 1951.

\bibitem[Fal79]{Faltings79}
G.~Faltings.
\newblock Algebraisation of some formal vector bundles.
\newblock {\em Ann. of Math. (2)}, 110(3):501--514, 1979.

\bibitem[Fal80]{Faltings80}
G.~Faltings.
\newblock Some theorems about formal functions.
\newblock {\em Publ. Res. Inst. Math. Sci.}, 16(3):721--737, 1980.

\bibitem[Fal81]{Faltings81}
G.~Faltings.
\newblock Formale {G}eometrie und homogene {R}\"aume.
\newblock {\em Invent. Math.}, 64(1):123--165, 1981.

\bibitem[Gas10]{Gasbarri10}
C.~Gasbarri.
\newblock {Analytic subvarieties with many rational points.}
\newblock {\em Math. Ann.}, 346(1):199--243, 2010.

\bibitem[Gil02]{Gillet02}
H.~Gillet.
\newblock Differential algebra---a scheme theory approach.
\newblock In {\em Differential algebra and related topics ({N}ewark, {NJ},
  2000)}, pages 95--123. World Sci. Publ., River Edge, NJ, 2002.

\bibitem[Gro61]{EGAIII1}
A.~Grothendieck.
\newblock \'{E}l\'ements de g\'eom\'etrie alg\'ebrique. {III}. \'{E}tude
  cohomologique des faisceaux coh\'erents. {I}.
\newblock {\em Inst. Hautes \'Etudes Sci. Publ. Math.}, (11):167, 1961.

\bibitem[Gro62]{GrothendieckFGA}
A.~Grothendieck.
\newblock {\em Fondements de la g\'eom\'etrie alg\'ebrique. [{E}xtraits du
  {S}\'eminaire {B}ourbaki, 1957--1962.]}.
\newblock Secr\'etariat math\'ematique, Paris, 1962.

\bibitem[Gro66]{Grothendieck66}
A.~Grothendieck.
\newblock On the de {R}ham cohomology of algebraic varieties.
\newblock {\em Inst. Hautes \'Etudes Sci. Publ. Math.}, (29):95--103, 1966.

\bibitem[Gro68]{GrothendieckSGA2}
A.~Grothendieck.
\newblock {\em Cohomologie locale des faisceaux coh\'erents et th\'eor\`emes de
  {L}efschetz locaux et globaux}.
\newblock S{\'e}minaire de G{\'e}om{\'e}trie Alg{\'e}brique du Bois-Marie, SGA
  2, 1962, Advanced Studies in Pure Mathematics, Vol. 2. North-Holland
  Publishing Co., Amsterdam, 1968.

\bibitem[Gun90]{Gunning90}
R.~C. Gunning.
\newblock {\em Introduction to holomorphic functions of several variables.
  {V}ol. {II}. Local theory}.
\newblock The Wadsworth \& Brooks/Cole Mathematics Series. Wadsworth \&
  Brooks/Cole Advanced Books \& Software, Monterey, CA, 1990.

\bibitem[Har68]{Hartshorne68}
R.~Hartshorne.
\newblock Cohomological dimension of algebraic varieties.
\newblock {\em Ann. of Math.}, 88:403--450, 1968.

\bibitem[Har70]{Hartshorne70}
R.~Hartshorne.
\newblock {\em Ample subvarieties of algebraic varieties}, volume 156 of {\em
  Lecture Notes in Mathematics}.
\newblock Springer-Verlag, Berlin, 1970.

\bibitem[Har75]{Hartshorne75}
R.~Hartshorne.
\newblock On the {D}e {R}ham cohomology of algebraic varieties.
\newblock {\em Inst. Hautes \'Etudes Sci. Publ. Math.}, (45):5--99, 1975.

\bibitem[Her12]{Herblot12}
M.~Herblot.
\newblock Algebraic points on meromorphic curves.
\newblock arXiv:1204.6336 [math. NT], 2012.

\bibitem[HM68]{HironakaMatsumura68}
H.~Hironaka and H.~Matsumura.
\newblock Formal functions and formal embeddings.
\newblock {\em J. Math. Soc. Japan}, 20:52--82, 1968.

\bibitem[Hod41]{Hodge41}
W.~V.~D. Hodge.
\newblock {\em The {T}heory and {A}pplications of {H}armonic {I}ntegrals}.
\newblock Cambridge University Press, Cambridge, England, 1941.

\bibitem[Hru96]{Hrushovski96}
E.~Hrushovski.
\newblock The {M}ordell-{L}ang conjecture for function fields.
\newblock {\em J. Amer. Math. Soc.}, 9(3):667--690, 1996.

\bibitem[HZ96]{HrushovskiZilber96}
E.~Hrushovski and B.~Zilber.
\newblock {Zariski geometries.}
\newblock {\em J. Am. Math. Soc.}, 9(1):1--56, 1996.

\bibitem[Ill94]{Illusie94}
L.~Illusie.
\newblock Crystalline cohomology.
\newblock In {\em Motives ({S}eattle, {WA}, 1991)}, volume~55 of {\em Proc.
  Sympos. Pure Math.}, pages 43--70. Amer. Math. Soc., Providence, RI, 1994.

\bibitem[Ill05]{Illusie05}
L.~Illusie.
\newblock Grothendieck's existence theorem in formal geometry.
\newblock In {\em Fundamental algebraic geometry}, volume 123 of {\em
  Mathematical Surveys and Monographs}, pages 179--233. Amer. Math. Soc.,
  Providence, RI, 2005.

\bibitem[Jan90]{Jannsen90}
U.~Jannsen.
\newblock {\em Mixed motives and algebraic {$K$}-theory}, volume 1400 of {\em
  Lecture Notes in Mathematics}.
\newblock Springer-Verlag, Berlin, 1990.

\bibitem[KSCT07]{KebekusSolaToma07}
S.~Kebekus, L.~Sol\'a~Conde, and M.~Toma.
\newblock {Rationally connected foliations after Bogomolov and McQuillan.}
\newblock {\em J. Algebr. Geom.}, 16(1):65--81, 2007.


\bibitem[KS53a]{KodairaSpencer53II}
K.~Kodaira and D.~C. Spencer.
\newblock Divisor class groups on algebraic varieties.
\newblock {\em Proc. Nat. Acad. Sci. U. S. A.}, 39:872--877, 1953.

\bibitem[KS53b]{KodairaSpencer53I}
K.~Kodaira and D.~C. Spencer.
\newblock Groups of complex line bundles over compact {K}\"ahler varieties.
\newblock {\em Proc. Nat. Acad. Sci. U. S. A.}, 39:868--872, 1953.

\bibitem[KGG{\etalchar{+}}04]{SpencerVita04}
J.~J. Kohn, P.~A. Griffiths, H.~Goldschmidt, E.~Bombieri, B.~Cenkl,
  P.~Garabedian, and L.~Nirenberg.
\newblock Donald {C}. {S}pencer (1912--2001).
\newblock {\em Notices Amer. Math. Soc.}, 51(1):17--29, 2004.

\bibitem[KP06]{KowalskiPillay06}
P.~Kowalski and A.~Pillay.
\newblock Quantifier elimination for algebraic {$D$}-groups.
\newblock {\em Trans. Amer. Math. Soc.}, 358(1):167--181 (electronic), 2006.


\bibitem[Lan62]{Lang62}
S.~Lang.
\newblock {Transcendental points on group varieties.}
\newblock {\em Topology}, 1:313--318, 1962.

\bibitem[Lan65]{Lang65}
S.~Lang.
\newblock {Algebraic values of meromorphic functions.}
\newblock {\em Topology}, 3:183--191, 1965.

\bibitem[Lan66a]{Lang66}
S.~Lang.
\newblock {Algebraic values of meromorphic functions. II.}
\newblock {\em Topology}, 5:363--370, 1966.

\bibitem[Lan66b]{Lang66b}
S.~Lang.
\newblock {\em {Introduction to transcendental numbers.}}
\newblock Addison-Wesley Series in Mathematics. {Addison-Wesley Publishing
  Company}, 1966.

\bibitem[LP91]{LePotier91}
J.~Le~Potier.
\newblock Fibr\'es de {H}iggs et syst\`emes locaux. {S}{\'e}minaire {B}ourbaki,
  {V}ol. 1990/91, {E}xpos\'e 737.
\newblock {\em Ast\'erisque}, (201-203):221--268, 1991.

\bibitem[Mal10]{Malgrange10}
B.~Malgrange.
\newblock Differential algebraic groups.
\newblock In {\em Algebraic approach to differential equations}, pages
  292--312. World Sci. Publ., Hackensack, NJ, 2010.

\bibitem[Man58]{Manin58}
Ju.~I. Manin.
\newblock Algebraic curves over fields with differentiation.
\newblock {\em Izv. Akad. Nauk SSSR. Ser. Mat.}, 22:737--756, 1958.

\bibitem[Man61]{Manin61}
Ju.~I. Manin.
\newblock The {H}asse-{W}itt matrix of an algebraic curve.
\newblock {\em Izv. Akad. Nauk SSSR Ser. Mat.}, 25:153--172, 1961.

\bibitem[Man63]{Manin63}
Ju.~I. Manin.
\newblock Rational points on algebraic curves over function fields.
\newblock {\em Izv. Akad. Nauk SSSR Ser. Mat.}, 27:1395--1440, 1963.

\bibitem[Mar00]{Marker2000}
D.~Marker.
\newblock Manin kernels.
\newblock In {\em Connections between model theory and algebraic and analytic
  geometry}, volume~6 of {\em Quad. Mat.}, pages 1--21. Dept. Math., Seconda
  Univ. Napoli, Caserta, 2000.

\bibitem[MM74]{MazurMessing74}
B.~Mazur and W.~Messing.
\newblock {\em Universal extensions and one dimensional crystalline
  cohomology}.
\newblock Lecture Notes in Mathematics, Vol. 370. Springer-Verlag, Berlin,
  1974.

\bibitem[Mes73]{Messing73}
W.~Messing.
\newblock {The universal extension of an abelian variety by a vector group.}
\newblock In {\em Symposia Mathematica, Vol. XI (Algebra commut., Geometria,
  Convegni 1971/1972, Roma INDAM)}, pages 359--372. Academic Press, London--New
  York, 1973.



\bibitem[Mum76]{Mumford76}
D.~Mumford.
\newblock {\em Algebraic geometry. {I}. {C}omplex projective varieties}, volume
  221 of {\em Grundlehren der Mathematischen Wissenschaften}.
\newblock Springer-Verlag, Berlin, 1976.

\bibitem[Oda69]{Oda69}
T.~Oda.
\newblock The first de {R}ham cohomology group and {D}ieudonn\'e modules.
\newblock {\em Ann. Sci. \'Ecole Norm. Sup. (4)}, 2:63--135, 1969.

\bibitem[Pil97a]{Pillay97BAMS}
A.~Pillay.
\newblock Model theory and {D}iophantine geometry.
\newblock {\em Bull. Amer. Math. Soc. (N.S.)}, 34(4):405--422, 1997.

\bibitem[Pil97b]{Pillay97}
A.~Pillay.
\newblock Some foundational questions concerning differential algebraic groups.
\newblock {\em Pacific J. Math.}, 179(1):179--200, 1997.

\bibitem[Pil04]{Pillay04}
A.~Pillay.
\newblock Algebraic {$D$}-groups and differential {G}alois theory.
\newblock {\em Pacific J. Math.}, 216(2):343--360, 2004.

\bibitem[Poi02]{Poincare02}
H.~Poincar\'e.
\newblock Sur les fonctions ab\'eliennes.
\newblock {\em Acta Math.}, 26:43--98, 1902.

\bibitem[Pui50]{Puiseux50}
V.~Puiseux.
\newblock Recherches sur les fonctions alg{\'e}briques.
\newblock {\em Journal de Math{\'e}matiques pures et appliqu{\'e}es},
  15:365--480, 1850.

\bibitem[Pui51]{Puiseux51}
V.~Puiseux.
\newblock Nouvelles recherches sur les fonctions alg{\'e}briques.
\newblock {\em Journal de Math{\'e}matiques pures et appliqu{\'e}es},
  16:228--240, 1851.

\bibitem[Ray75]{RaynaudMe75}
Mich{\`e}le Raynaud.
\newblock {\em Th\'eor\`emes de {L}efschetz en cohomologie coh\'erente et en
  cohomologie \'etale}.
\newblock Soci\'et\'e Math\'ematique de France, Paris, 1975.
\newblock Bull. Soc. Math. France, M{\'e}m. No. 41, Suppl{\'e}ment au Bull.
  Soc. Math. France, Tome 103.

\bibitem[Rem56]{Remmert56}
R.~Remmert.
\newblock Meromorphe {F}unktionen in kompakten komplexen {R}\"aumen.
\newblock {\em Math. Ann.}, 132:277--288, 1956.

\bibitem[Rie57]{Riemann57IV}
B.~Riemann.
\newblock Theorie der {A}bel'schen {F}unctionen.
\newblock {\em J. reine angew. Math.}, 54:115--155, 1857.

\bibitem[Ros58]{Rosenlicht58}
M.~Rosenlicht.
\newblock Extensions of vector groups by abelian varieties.
\newblock {\em Amer. J. Math.}, 80:685--714, 1958.

\bibitem[Sch41]{Schneider41}
Th. Schneider.
\newblock Zur {T}heorie der {A}belschen {F}unktionen und {I}ntegrale.
\newblock {\em J. Reine Angew. Math.}, 183:110--128, 1941.

\bibitem[Sch57]{Schneider57}
Th. Schneider.
\newblock {\em Einf\"uhrung in die transzendenten {Z}ahlen}.
\newblock Springer-Verlag, Berlin, 1957.

\bibitem[Ser54]{Serre53}
J.~P. Serre.
\newblock Fonctions automorphes: quelques majorations dans le cas o\`u
  ${X}/{G}$ est compact.
\newblock {\em S\'eminaire H. Cartan}, 6, 1953--1954.
\newblock Expos\'e 2.

\bibitem[Ser56]{Serre56}
J.~P. Serre.
\newblock G\'eom\'etrie alg\'ebrique et g\'eom\'etrie analytique.
\newblock {\em Ann. Inst. Fourier, Grenoble}, 6:1--42, 1955--1956.

\bibitem[Ser59]{Serre59}
J.-P. Serre.
\newblock {\em Groupes alg\'ebriques et corps de classes}.
\newblock Publications de l'institut de math\'ematique de l'universit\'e de
  Nancago, VII. Hermann, Paris, 1959.

\bibitem[Sev61]{Severi61}
F.~Severi.
\newblock {\em Funzioni quasi abeliane}, volume~20 of {\em Pontificiae
  Academiae Scientiarum Scripta Varia}.
\newblock Pontificia Academia Scientiarum, Vatican City, second augmented
  edition, 1961.

\bibitem[Sha77]{Shafarevich77}
I.R. Shafarevich.
\newblock {\em {Basic algebraic geometry. 2nd ed.}}
\newblock {Springer Study Edition. Berlin-Heidelberg-New York:
  Springer-Verlag}, 1977.

\bibitem[Sie55]{Siegel55}
C.~L. Siegel.
\newblock Meromorphe {F}unktionen auf kompakten analytischen
  {M}annigfaltigkeiten.
\newblock {\em Nachr. Akad. Wiss. G\"ottingen. Math.-Phys. Kl. IIa.},
  1955:71--77, 1955.

\bibitem[Sim94a]{Simpson94I}
C.~T. Simpson.
\newblock Moduli of representations of the fundamental group of a smooth
  projective variety. {I}.
\newblock {\em Inst. Hautes \'Etudes Sci. Publ. Math.}, (79):47--129, 1994.

\bibitem[Sim94b]{Simpson94II}
C.~T. Simpson.
\newblock Moduli of representations of the fundamental group of a smooth
  projective variety. {II}.
\newblock {\em Inst. Hautes \'Etudes Sci. Publ. Math.}, (80):5--79 (1995),
  1994.

\bibitem[Thi54]{Thimm54}
W.~Thimm.
\newblock \"{U}ber meromorphe {A}bbildungen von komplexen {M}annigfaltigkeiten.
\newblock {\em Math. Ann.}, 128:1--48, 1954.

\bibitem[Wal79]{Waldschmidt79}
M.~Waldschmidt.
\newblock {\em Nombres transcendants et groupes alg\'ebriques}, volume 69-70 of
  {\em Ast\'erisque}.
\newblock Soci\'et\'e Math\'ematique de France, Paris, 1979.
\newblock With appendices by D. Bertrand and J.-P. Serre.

\bibitem[W{\"u}s84]{Wuestholz84}
G.~W{\"u}stholz.
\newblock {Zum Periodenproblem.}
\newblock {\em Invent. Math.}, 78:381--391, 1984.

\bibitem[Zar51]{Zariski51}
O.~Zariski.
\newblock Theory and applications of holomorphic functions on algebraic
  varieties over arbitrary ground fields.
\newblock {\em Mem. Amer. Math. Soc.}, 1951(5):90, 1951.

\bibitem[Zar71]{Zariski71}
O.~Zariski.
\newblock {\em Algebraic surfaces}, volume~61 of {\em Ergebnisse der Mathematik
  und ihrer Grenzgebiete}.
\newblock Springer-Verlag, New York, supplemented edition, 1971.
\newblock With appendices by S. S. Abhyankar, J. Lipman, and D. Mumford.

\end{thebibliography}

\end{document}